\documentclass{article} \textwidth=135mm \textheight= 185mm
\usepackage{amssymb,amsmath,amsthm}
\newtheorem{theorem}{Theorem}
\newtheorem{definition}{Definition}
\newtheorem{lemma}{Lemma}
\newtheorem{proposition}{Proposition}
\newtheorem{corollary}{Corollary}

\date{}
\numberwithin{equation}{section} \numberwithin{theorem}{section}
\numberwithin{lemma}{section} \numberwithin{corollary}{section}
\numberwithin{remark}{section} \numberwithin{proposition}{section}
\numberwithin{definition}{section}

\begin{document}
\newcommand{\n}{\noindent}
\newcommand{\vs}{\vskip}

\title{ On the $A$-Obstacle Problem and the Hausdorff Measure of its Free Boundary}

\author{S. Challal$^1$, A. Lyaghfouri$^1$ and  J. F. Rodrigues$^2$\\
\\
$^1$ Fields Institute, 222 College Street\\
Toronto M5T 3J1, Canada\\
\\
$^2$ University of Lisbon/CMAF\\
Av. Prof. Gama Pinto, 2\\
1649-003 Lisboa, Portugal} \maketitle

\vs 0.5cm
\begin{abstract} In this paper we prove existence and uniqueness of an entropy solution
to the  $A$-obstacle problem, for $L^1-$data. We also extend the Lewy-Stampacchia inequalities to the
general framework of $L^1-$data, and show convergence and stability results. We then prove
that the free boundary has finite $(N-1)$-Hausdorff measure, which completes previous works on this
subject by Caffarelli for the Laplace operator and by Lee and Shahgholian for the $p-$Laplace operator when $p>2$.

\end{abstract}

\vs 0.3cm

\n MSC: 35R35, 35B05, 35J60

\vs 0.3cm

\n Key words : Obstacle problem, Entropy Solution, $A$-Laplace
operator, Lewy-Stampacchia Inequalities, Stability, Free boundary,
Hausdorff measure.

\vs 0.5cm

\section{Introduction}\label{S:intro}

\vs 0.5cm

\n We consider the obstacle problem
in a bounded domain of $\mathbb{R}^N$  $(N \geqslant
2)$, associated with  data $f\in L^1(\Omega)$, an admissible obstacle $\psi$ and boundary data
$g$ in $W^{1,A}(\Omega) \cap L^\infty (\Omega)$:
\begin{equation}\label{1.2}
\begin{cases}
     &\Delta_A u=f\qquad \text{in}\quad \{u>\psi\}, \\
     &  u\geqslant \psi \quad \text{in }\quad \Omega,\\
     & u=g \quad  \text{on }\quad \partial\Omega.
  \end{cases}
\end{equation}

\n where the $A-$Laplace operator
\begin{equation}\label{1.2}
\Delta_A u= div\Big({a(|\nabla u|)\over
|\nabla u|} \nabla u \Big)=div(\nabla_A u),
\end{equation}

\n is associated with a $C^1$ function $a~:[0,+\infty)\rightarrow[0,+\infty)
$, $a(0)=0$ and with $\displaystyle{ A(t) = \int_{0}^{t} a(s)
ds}$ the $N-$function of the Orlicz-Sobolev space
$W^{1,A}(\Omega)= \{ u\in L^A(\Omega):~ \nabla u \in
 L^A(\Omega)\}$, which is the usual subspace of $W^{1,1}(\Omega)$ associated with the norm
$ \|u\|_{W^{1,A}(\Omega)} = \|u\|_A + \|\nabla u\|_A $, and where
$L^A(\Omega)$ is the Orlicz space equipped with the respective Luxembourg norm (see \cite{[A]}).

\n The obstacle  problem (1.1), as it is well known (see, for instance, \cite{[R1]}), may be formulated
as a minimization problem of the functional
\begin{equation}\label{1.3}
J(v)=\int_\Omega A(|\nabla v|)dx +\int_\Omega fvdx,
\end{equation}
for a given $f\in L^\infty(\Omega)$, in the convex set

\begin{equation}\label{1.4}
K_{\psi,g}= \{ v\in  W^{1,A}(\Omega), \quad v-g\in
W^{1,A}_0(\Omega),\quad v\geqslant \psi \quad \hbox{ a.e. in }
\Omega\},
\end{equation}

\n provided $K_{\psi,g}\neq \emptyset$, which holds if $(\psi-g)^+\in W_0^{1,A}(\Omega)$, for instance,
and where $W_0^{1,A}(\Omega)$ is the closure of ${\cal D}(\Omega)$ in $W^{1,A}(\Omega)$.

\vs 0.2cm \n In this work we intend to extend several known global and local properties of the solution
of the obstacle problem and of its free boundary $(\partial\{u>\psi\})\cap \Omega$ to a large class of
degenerate and singular elliptic operators under the natural condition which generalizes the Ladyzhenskaya-Uraltseva
operators, for some constants $a_0$ and $a_1$
\begin{equation}\label{1.5}
 0<a_0\leqslant { t a^\prime (t) \over a(t)}  \leqslant a_1  \qquad
 \forall t>0, \qquad a_0, a_1 \hbox{  are positive constants}.
\end{equation}

\n This not only includes the case of the $p-$Laplacian $\Delta_p$ (when $a_0=a_1=p-1>0$),
but also the interesting case of a variable exponent $p=p(t)>0$
\begin{equation}\label{1.6}
\Delta_A u=div\big(|\nabla u|^{p(|\nabla
u|)-2}\nabla u\big),
\end{equation}
\n corresponding to set $a(t)=t^{p(t)-1}$, for which (1.5) holds if
$a_0\leqslant t\ln(t)p'(t)+p(t)-1\leqslant a_1$, for all  $t>0.$

\n Other examples of functions satisfying (1.5) are given by $a(t)=t^\alpha\ln(\beta t+\gamma)$,
with $\alpha, \beta, \gamma>0$, or by discontinuous power transitions like
$a(t)=c_1t^{\alpha}$, if  $0\leq t<t_0$, and $a(t)=c_2t^{\beta}+c_3$, if $t\geqslant t_0$,
where $\alpha, \beta, t_0$, are positive numbers, and $c_1$,
$c_2$ and $c_3$ are real numbers such that $a(t)$ is a $C^1$
function (take $a_0=\min(\alpha,\beta)$ and $a_1=\max(\alpha,\beta)$).
Since linear combinations with positive constants, products and compositions of $C^1$
functions satisfying (1.5) also satisfy a similar condition, with different
positive constants $a_0$ and $a_1$, we conclude that there exists a large class
of $N-$functions $A$ that are included in this work. Without loss of generality, we will
assume that $a_0\leqslant 1\leqslant a_1$.

\vs0.2cm \n As it is well known, minimizing (1.3) in (1.4) is equivalent to solving the
following variational inequality
\begin{equation}\label{1.5}
u\in K_{g,\psi}~:\quad\int_\Omega \nabla_A u
 . \nabla (v-u)dx +\int_\Omega f(v-u)dx \geqslant 0,
\quad  \forall v \in K_{\psi,g},
\end{equation}

\n where we used the notation $\displaystyle{\nabla_A u={{a(|\nabla u|)}\over {|\nabla u|}}\nabla u}$, for the $A-$gradient
associated with the $\Delta_A$ introduced in (1.2).

\vs0.2cm \n Since the assumption (1.5) implies that $-\Delta_A$ is a coercive and strictly
monotone operator in $W^{1,A}(\Omega)$ (see \cite{[L1]} or \cite{[CL2]}), which is a reflexive and separable
Banach space (see \cite{[A]}), the existence and uniqueness of a variational solution
to (1.7) follows easily. In fact, since (1.5) also implies the inequality (2.4), which in turn
leads to the continuous embedding $L^A(\Omega)\subset L^{a_0+1}(\Omega)$, we have
$W^{1,A}(\Omega)\subset W^{1,a_0+1}(\Omega)$. Therefore by Sobolev embedding, we have
\begin{equation}\label{1.8}
W^{1,A}(\Omega)\subset L^q(\Omega)
\end{equation}
\n for $q=N(a_0+1)/(N-(a_0+1))$ if $1+a_0<N$, for any $q<\infty$ if $a_0+1=N$, and
for $q=\infty$ if $a_0+1>N$. Hence, in this later case, we may still solve (1.7) in
$W^{1,A}(\Omega)$ with $f$ only in $L^1(\Omega)$, and with $f\in L^r(\Omega)$ for $r>1$
if $a_0+1=N$, or any $r\geq q'=N(1+a_0)/(Na_0+a_0+1)$ if $a_0+1<N$.

\n However if $a_0+1<N$, for a general $f\in L^1(\Omega)$, the second integral in (1.7)
is not defined, and following \cite{[BBGGOV]} and \cite{[BC]}, we are led to the more
general definition of a solution to the obstacle problem, using the truncation function
$$T_s(r)=\max(-s,\min(s,r)),~~r\in \mathbb{R}:$$

\n An \emph{entropy solution} of the obstacle problem (1.1) is a
measurable function $u$ such that $u\geqslant \psi$ a.e. in
$\Omega$, and, for every $s>0$, $T_s(u)-T_s(g)\in
W^{1,A}_0(\Omega)$ and
\begin{equation}\label{1.9}
\int_\Omega \nabla_A u
 . \nabla (T_s(v-u))dx+\int_\Omega fT_s(v-u)dx\geqslant 0,
\qquad  \forall v \in K_{\psi,g}.
\end{equation}

\n We observe that no global integrability condition is required on $u$ nor on its gradient,
but as it was shown in \cite{[BBGGOV]} for Sobolev spaces, if $T_s(u)\in W^{1,A}(\Omega)$
for all $s>0$, then there exists a unique measurable vector field $U :\, \Omega\,\rightarrow\, \mathbb{R}^N$ such that
$\nabla (T_s(u))=\chi_{\{|u|<s\}}U$ a.e. in $\Omega$, $s>0$, which, in fact, coincides with the
standard distributional gradient of $\nabla u$ whenever $u\in W^{1,1}(\Omega)$.

\vs0.2cm \n Our first results concern global properties of entropy solutions in the
$L^1(\Omega)$-framework. Without loss of generality we can assume that $a_0<N-1$.

\begin{theorem}\label{t1.1} For given $f\in L^1(\Omega)$ and $\psi,g \in W^{1,A}(\Omega)\cap L^\infty(\Omega)$
such that $K_{\psi,g}\neq\emptyset$, there exists a unique entropy solution $u$ to
(1.9), which moreover satisfies
\begin{equation}\label{1.10}
u\in W^{1,p}(\Omega),\quad \text{ for all }p,\,\, 1<p<{N\over{N-1}}a_0.
\end{equation}

\n In addition, $u$ depends continuously on $f$, i.e., if $f_n\,\rightarrow\,f$ in $L^1(\Omega)$
and $u_n$ is associated with $(f_n,\psi,g)$, then
\begin{equation}\label{1.11}
u _n\,\rightarrow\,u\quad\text{in } W^{1,p}(\Omega),\quad \text{ for all } p,~ 1<p<{N\over{N-1}}a_0.
\end{equation}
\end{theorem}

\vs 0,2 cm\n We notice that existence results for the obstacle problem in more general Orlicz-Sobolev spaces
for Leray-Lions type operators with $L^1(\Omega)$ data have been obtained in \cite{[ABR]}. Also in the
$L^1$-framework for Orlicz-Sobolev spaces with a $N-$function satisfying the $\Delta_2-$condition
(see \cite{[A]}, p. 231), but not necessarily (1.5), the existence and uniqueness of an entropy solution
to the obstacle problem with homogeneous boundary data $g=0$ is given in \cite{[BB]}, and with additional assumptions
on the obstacle has been obtained earlier. We observe that also in Remark 3.2 of \cite{[BB]}, for that case,
a result similar to (1.10), improving the previous regularity result of \cite{[BC]} for the $p-$Laplacian under the unnecessary
restriction $2-1/N<p<N$.

\vs 0,2 cm\n The next theorem extends the Lewy-Stampacchia inequalities
for the entropy solution of the obstacle problem with $f$ in $L^1(\Omega)$.
This is particularly important to easily obtain the H\"{o}lder continuity of the
gradient of the solution, in the case of bounded data, as an immediate consequence
of the known results of \cite{[L1]} (see also \cite{[L2]}) for the $\Delta_A$
operator.

\begin{theorem}\label{t1.2} Assume that $(f-\Delta_A\psi)^+\in L^1(\Omega)$, in particular
if $\Delta_A\psi\in L^1(\Omega)$, the entropy solution $u$ is such that
 $\Delta_A u\in L^1(\Omega)$ and the following Lewy-Stampacchia inequalities hold
\begin{eqnarray}\label{1.12}
f-(f-\Delta_A\psi)^+\leqslant \Delta_A u\leqslant f\quad\text{a.e. in}\quad\Omega.
\end{eqnarray}
Moreover, if $f$ and $(f-\Delta_A\psi)^+\in L^\infty(\Omega)$, then
$u\in C^{1,\alpha}(\Omega)$, and also $u\in C^{1,\alpha}(\overline{\Omega})$,
provided that $g\in C^{1,\alpha}(\partial\Omega)$ and also $\partial\Omega$
is a $C^{1,\alpha}$ surface, for some $0<\alpha<1$.
\end{theorem}

\vs 0,2 cm\n We observe that the $C^{1,\alpha}$ regularity of the solution
of the one and two obstacle problems for these type of degenerate operators was obtained
in \cite{[L2]} for $C^{1,\alpha}$ obstacles.

\vs 0,1 cm\n The inequalities (1.12) and the $C^{1,\alpha}$ regularity are not
enough to prove in all generality that the entropy solution of the obstacle problem
also solves the semilinear equation
\begin{equation}\label{1.13}
\Delta_A u-(\Delta_A\psi-f)\chi_{\{u=\psi\}}=f\quad\text{a.e. in}\quad\Omega,
\end{equation}
\n as, for instance, in the case of the Laplacian (see \cite{[R1]}, for example).
Here $\chi_{\{u=\psi\}}$ denotes the characteristic function of the set $\{u=\psi\}$, i.e.
$\chi_{\{u=\psi\}}(x)=1$ if $x\in \{u=\psi\}$, and $\chi_{\{u=\psi\}}(x)=0$ if $x\in \{u>\psi\}$,.

\vs 0,1 cm\n However (1.13) holds, for example, when $\psi=0$ and $0<\lambda_0\leq f\leq \Lambda_0$,
as a consequence of the porosity of the free boundary (see \cite{[CL1]}) that yields
$(\partial\{u=\psi\})\cap\Omega$ is of Lebesgue measure zero. In fact, the weaker regularity condition
(1.13) enables us to show a first corollary on the stability of free boundaries,
extending the results of \cite{[R1]} and \cite{[R2]}.

\begin{corollary}\label{c1.1} Under the $L^1$ assumptions of Theorem 1.1 and Theorem 1.2, for
$f$, $f_n$ and $(\psi,g)$, assume also that the limit entropy solution $u$ satisfies (1.13) and suppose
that
\begin{equation}\label{1.14}
\Delta_A\psi\neq f \quad\text{a.e. in}~~ \Omega.
\end{equation}
Then, the respective characteristic functions of the coincidence sets $\{u_n=\psi\}$ converges in the
following sense
\begin{equation}\label{1.15}
\chi_{\{u_n=\psi\}}\,\rightarrow\,\chi_{\{u=\psi\}}\quad\text{in} \quad L^s(\Omega)\quad\text{for all}~~s<\infty.
\end{equation}
\end{corollary}

\vs 0,2 cm\n A third global property of entropy solutions of the obstacle problem is related to the
monotonicity and $T-$accretivity of $\Delta_A$ in $L^1(\Omega)$.

\vs 0,3 cm

\begin{theorem}\label{t1.3}
Let the $L^1$ assumptions of Theorem 1.1 and Theorem 1.2 hold for two sets
of data $f_1, \psi, g)$ and $f_2, \psi, g)$ and denote by $u_i$,
$\xi_i=\Delta_A u_i-f_i$, $i=1, 2$, their respective entropy solutions.
Assume moreover that $\psi\in W_{loc}^{2,\infty}(\Omega)$. Then we have
\begin{equation}\label{1.16}
\int_\Omega|\xi_1-\xi_2|dx\leqslant \int_\Omega|f_1-f_2|dx.
\end{equation}
\end{theorem}

\vs 0,2 cm\n Similarly to earlier works for the $p-$Laplacian (see \cite{[R1]} for $p=2$ and \cite{[R2]}),
from this $L^1$ contraction property for the obstacle problem we can easily obtain an estimate
for the stability of two coincidence sets $I_i=I(u_i)=\{u_i=\psi\}$, $i=1, 2$, in terms of the
Lebesgue measure ${\cal L}^N$ of their symmetric difference
$$I_1\div I_2=(I_1\setminus I_2)\cup(I_2\setminus I_1).$$

\vs 0,5 cm

\begin{corollary}\label{c1.2}
In addition to the assumptions of Theorem 1.3, suppose that both entropy solutions satisfy (1.13)
and the data satisfy the local non-degeneracy condition, for some constant $\lambda$,
\begin{equation}\label{e1.17}
\Delta_A\psi-f_i\leqslant -\lambda<0,\quad\text{a.e. in}\quad
\omega,\quad i=1, 2
\end{equation}
in a measurable subset $\omega\subset \Omega$. Then we have
\begin{equation}\label{e1.18}
{\cal L}^N\big((I_1\div I_2)\cap \omega\big)\leqslant {1\over\lambda}\int_\Omega|f_1-f_2|dx.
\end{equation}
\end{corollary}

\vs 0,5 cm In order to prove local properties, we shall restrict ourselves to the case where
\begin{equation}\label{e1.19}
\psi=0\quad\text{and}\quad f\in L^\infty(\Omega).
\end{equation}

\vs 0,2 cm\n We start by estimating the growth of the solution an its gradient in an open neighborhood
$\omega\subset \Omega$  of a generic free boundary point $x_0\in (\partial \{u>0\})\cap\omega$, where
we shall assume the existence of two constants $\lambda, \Lambda$, such that
\begin{equation}\label{e1.20}
0<\lambda_0\leqslant f\leqslant \Lambda_0\quad\text{a.e. in}\quad
\omega.
\end{equation}

\n We denote by $a^{-1}$ the inverse function of $a$ and we introduce the associated complementary
function $\widetilde{A}$ :

\begin{equation}\label{e-21}
\widetilde{A}(t)=\int_{0}^{t} a^{-1}(s) ds.
\end{equation}

\n Then we have

\vs 0,5cm
\begin{theorem}\label{t1.4} Under the assumptions (1.20) and (1.21), the solution $u$
to the obstacle problem, and its gradient have the following growth rates near a free boundary
point $x_0\in (\partial \{u>0\})\cap\Omega$
\begin{equation}\label{e-22}
0\leqslant u(x) \leqslant C_0 \widetilde{A} (|x-x_0|) \qquad \forall x\in \omega
\end{equation}
\begin{equation}\label{e-23}
|\nabla u(x)| \leqslant C_1 a^{-1}(|x-x_0|) \qquad \forall x\in \omega,
\end{equation}
\n for positive constants $C_0, C_1$ depending on $N$, $a_0$, $a_1$, $\lambda_0$, $\Lambda_0$,
and on the $L^\infty_{loc}$ and $C^{1,\alpha}_{loc}$ norms of $u$ respectively.
\end{theorem}

\vs 0,5 cm\n In fact (1.22) was already proved in \cite{[CL1]} and extends the well-known
quadratic growth property observed by Caffarelli in the pioneering work \cite{[C]} for the Laplacian
and the corresponding $p/(p-1)$ growth for the $p-$Laplacian obtained in \cite{[KKPS]} for
$1<p<\infty$. The natural consequence (1.23) on the growth of the gradient was also extended
to the $p-$Laplacian, but only for $p>2$ in \cite{[LS]}.

\n In order to obtain further regularity on the free boundary, we shall extend local
$L^2-$estimates for the second derivatives of the solution obtained in \cite{[CL3]}, by
assuming additionally
\begin{equation}\label{e-24}
\nabla f \in {\cal M}^N_{loc}(\Omega),
\end{equation}

\n where ${\cal M}^N_{loc}(\Omega)$ is the Morrey space (see \cite{[MZ]} for definition of this space).

\begin{equation}\label{e-25}
\displaystyle{{a(t)\over t}}~~\text{ is monotone}~~\text{ for}~~0<t<t_*,
\end{equation}

\n for some $t_*>0$, where the monotonicity near zero may be either non-increasing or
non-decreasing.

\vs 0.5cm

\begin{theorem}\label{t1.5} Under the assumptions (1.20)-(1.21) and (1.24)-(1.25), the
free boundary of the $A-$Obstacle problem is locally of finite Hausdorff measure. More
precisely, if $x_0\in (\partial \{u>0\})\cap\omega$ and $\omega\subset\Omega$, then
there exists a positive  constant $C_0$ such that
\begin{equation*}
{\cal H}^{N-1}(\partial \{u>0\}\cap B_r(x_0))\leqslant C_0 r^{n-1},\quad B_r(x_0)\subset\omega.
\end{equation*}
\end{theorem}

\vs 0,5cm\n This result is a natural extension of the same property known for the
$p-$Obstacle problem, obtained in \cite{[C]} for $p=2$ and recently in \cite{[LS]} for
$p>2$. It is new for the case $1<p<2$.

\n As a consequence of this theorem and general results
from geometric measure theory, the coincidence set $I(u)=\{u=0\}$ is a set of finite
perimeter and, hence, the free boundary consists, up to a "small" singular set, of a
union of an at most countable number of $C^1$ hypersurfaces.

\vs 0,2cm\n The arguments on the growth of the solution near the free boundary, that
are used to prove Theorem 1.5, also yield another stability result of the coincidence
set $I(u)=\{u=0\}$, in terms of $L^\infty$ variations of the solution, extending the results
of \cite{[C]} and \cite{[LS]}.

\vs 0,5cm

\begin{corollary}\label{c1.3} If $u_1$ and $u_2$ are two local solutions of the
$A-$Obstacle problem, under the conditions of Theorem 1.5, then there exists a
positive constant $C$ such that
\begin{equation}\label{e1.26}
{\cal L}^N\big((I(u_1)\div I(u_2))\cap\omega'\big)\leqslant C\widetilde{A}^{-1}\big(||u_1-u_2||_{\infty,\omega}\big),\quad
\omega'\subset\subset\omega.
\end{equation}

\n For each $\epsilon>0$, such that $||u_1-u_2||_{\infty,\omega}<\widetilde{A}^{-1}(\epsilon)$,
we have
\begin{equation}\label{e1.27}
(I(u_2)\cap\omega')_{(-C\epsilon)}\subset I(u_1)\cap\omega'\subset\{u_2<\widetilde{A}^{-1}(\epsilon)\},
\end{equation}
\n where $(I(u_2)\cap\omega')_{(-C\epsilon)}=\{x\in I(u_2)\cap\omega'~:~d(x,\{u_2>0\})$.
\end{corollary}

\vs 0.5 cm \n

\section{Some Auxiliary Results}\label{S2}

\vs 0,5cm\n We start by recalling some useful inequalities that follow from the assumption (1.5) (see
\cite{[L1]}, \cite{[MW]}) :

\begin{equation}\label{2.1}
  {t a(t) \over 1+a_1} \leqslant A(t) \leqslant t a(t)\qquad
  \forall  t\geqslant 0,
\end{equation}

\begin{equation}\label{2.2}
 sa(t) \leqslant ta(t)+sa(s) \quad \forall s,t \geqslant 0,
\end{equation}

\begin{equation}\label{3.3}
  \min ( s^{a_0},s^{a_1})
   a(t) \leqslant a(st) \leqslant \max ( s^{a_0},s^{a_1})
   a(t) \quad \forall s,t \geqslant 0,
\end{equation}

\begin{equation}\label{2.4}
  \min ( s^{1+a_0},s^{1+a_1})
   {{A(t)}\over {1+a_1}} \leqslant A(st) \leqslant (1+a_1)\max ( s^{1+a_0},s^{1+a_1})
   A(t) \quad \forall s,t \geqslant 0,
\end{equation}

\begin{equation}\label{2.5}
  \min ( s^{1/a_0},s^{1/a_1})
   a^{-1}(t) \leqslant a^{-1}(st) \leqslant \max ( s^{1/a_0},s^{1/a_1})
   a^{-1}(t)\quad \forall s,t \geqslant 0,
\end{equation}

\begin{equation}\label{2.6}
  {a_0 \over 1+a_0}  t a^{-1}(t)\leqslant  \widetilde{A}(t) \leqslant t a^{-1}(t)\qquad
  \forall  t\geqslant 0.
\end{equation}

\n We also have the following monotonicity inequality (see \cite{[CL2]} and \cite{[CL4]})

\begin{equation}\label{2.7}
 \Big({{a(|\xi|)}\over {|\xi|}}\xi-{{a(|\zeta|)}\over {|\zeta|}}\zeta\Big). (\xi-\zeta) \geqslant C(A,N)|\xi-\zeta|^2{ a\big((|\xi|^2+|\zeta|^2)^{1/2}\big)\over {(|\xi|^2+|\zeta|^2)^{1/2}}} \qquad \forall (\xi,\zeta)\in
 \mathbb{R}^{2N}\setminus\{ 0\}.
\end{equation}

\vs 0,5 cm\n We start by establishing the Lewy-Stampacchia inequalities for the variational solution.

\begin{proposition}\label{p2.1} Assume that $f,(f-\Delta_A\psi)^+\in L^\infty(\Omega)$. Then
the variational solution $u$ of the obstacle problem associated with
$(f,\psi,g)$ satisfies the Lewy-Stampacchia inequalities
\begin{eqnarray*}
f-(f-\Delta_A\psi)^+\leqslant \Delta_A u\leqslant f\quad\text{a.e.
in}\quad\Omega.
\end{eqnarray*}
\end{proposition}

\vs 0,5cm\n\emph{Proof.} Let $\epsilon\in(0,1)$, $h=(f-\Delta_A\psi)^+$ and consider the following approximated problem

\begin{equation*}
(P_\epsilon)\begin{cases} & \text{Find } u_\epsilon\in
g+W_0^{1,A}(\Omega)
\text{ such that } :\\
&  -\Delta_A u_\epsilon=-f+h(1-H_\epsilon(u_\epsilon-\psi))~~
\text{in}~~ (W^{1,A}(\Omega))',
\end{cases}
\end{equation*}

\n where $H_\epsilon(t)=(T_\epsilon(t))^+$.

\n Using the standard theory of monotone operators and Schauder
fixed point theorem, one can prove that $(P_\epsilon)$ has a
solution. Using the monotonicity of $H_\epsilon$ and the inequality (2.7), one
can also prove the uniqueness of this solution. Moreover we have

\begin{equation}\label{2.8}
|\Delta_A
u_\epsilon|_\infty=|f-(1-H_\epsilon(u_\epsilon-\psi))h|_\infty\leqslant
|f|_\infty+|h|_\infty.
\end{equation}
We claim that

\begin{equation}\label{2.9}
u_\epsilon-\epsilon\leqslant u \leqslant u_\epsilon\quad\text{a.e.
in}~~\Omega.
\end{equation}
First we show that $u_\epsilon\in K_{\psi,g}$. For this purpose, it
is enough to verify that $u_\epsilon \geqslant\psi$ a.e. in
$\Omega$. Since $(\psi-g)^+\in W_0^{1,A}(\Omega)$,
  $(\psi-u_\epsilon)^+$ is a test function for $(P_\epsilon)$. We
obtain

\begin{eqnarray*}
\int_\Omega \nabla_A u_\epsilon
 . \nabla (\psi-u_\epsilon)^+dx &=&
  \int_\Omega (-f+(1-H_\epsilon(u_\epsilon-\psi))h)(\psi-u_\epsilon)^+dx\nonumber\\
&=&\int_\Omega (-f+h)(\psi-u_\epsilon)^+dx.
\end{eqnarray*}
This leads to

\begin{eqnarray*}
\int_\Omega \big(\nabla_A \psi-\nabla_A u_\epsilon\big)
 . \nabla (\psi-u_\epsilon)^+dx &=&
  \int_\Omega (f-\Delta_A\psi-(f-\Delta_A\psi)^+)(\psi-u_\epsilon)^+dx\nonumber\\
&=&-\int_\Omega (f-\Delta_A\psi)^-(\psi-u_\epsilon)^+dx\leqslant 0.
\end{eqnarray*}

\n Using (2.7), we deduce that $\nabla (\psi-u_\epsilon)^+=0$ a.e. in
$\Omega$. But because $ (\psi-u_\epsilon)^+=(\psi-g)^+=0$ on
$\partial\Omega$, we obtain $(\psi-u_\epsilon)^+=0$ a.e. in
$\Omega$. Hence we have proved that $u_\epsilon \geqslant\psi$ a.e.
in $\Omega$.

\vs 0,2 cm \n Next we prove that $u\leqslant u_\epsilon$ a.e. in
$\Omega$. To do this we use $u-(u-u_\epsilon)^+$ as a test function
for the obstacle problem associated with $(f,\psi,g)$. We obtain

\begin{eqnarray}\label{2.10}
\int_\Omega \nabla_A u
 . \nabla (u-u_\epsilon)^+dx &\leqslant&
  \int_\Omega -f(u-u_\epsilon)^+dx.
\end{eqnarray}

\n Using $(u-u_\epsilon)^+$ as a test function for the problem
$(P_\epsilon)$, we obtain

\begin{eqnarray}\label{2.11}
\int_\Omega \nabla_A u_\epsilon
 . \nabla (u-u_\epsilon)^+dx &=&
  \int_\Omega (-f+(1-H_\epsilon(u_\epsilon-\psi))h)(u-u_\epsilon)^+dx.
\end{eqnarray}

\n Subtracting (2.10) from (2.11), we obtain

\begin{eqnarray*}
\int_\Omega \big(\nabla_A u-\nabla_A u_\epsilon\big)
 . \nabla (u-u_\epsilon)^+dx \leqslant
 -\int_\Omega h(1-H_\epsilon(u_\epsilon-\psi))(u-u_\epsilon)^+dx\leqslant 0.
\end{eqnarray*}
Using the fact that we have $ (u-u_\epsilon)^+=(g-g)^+=0$ on $\partial\Omega$, we conclude
as above that $(u-u_\epsilon)^+=0$ a.e. in $\Omega$, or equivalently
that $u_\epsilon \geqslant u$ a.e. in $\Omega$.

\vs 0,2 cm \n Now we prove that $u_\epsilon-\epsilon\leqslant u$
a.e. in $\Omega$. To do this we use $u+(u_\epsilon-u-\epsilon)^+$ as
a test function for the obstacle problem associated with
$(f,\psi,g)$ and $(u_\epsilon-u-\epsilon)^+$ as a test function for
the problem $(P_\epsilon)$. We obtain

\begin{eqnarray}\label{2.12-13}
&&-\int_\Omega \nabla_A u
 . \nabla (u_\epsilon-u-\epsilon)^+dx \leqslant
  \int_\Omega f(u_\epsilon-u-\epsilon)^+dx\\
&&\int_\Omega \nabla_A u_\epsilon
 . \nabla (u_\epsilon-u-\epsilon)^+dx =
  \int_\Omega
  (-f+(1-H_\epsilon(u_\epsilon-\psi))h)(u_\epsilon-u-\epsilon)^+dx.\nonumber\\
  &&
\end{eqnarray}

\n Adding (2.12) and (2.13), we obtain

\begin{eqnarray*}
\int_\Omega \big(\nabla_A (u_\epsilon-\epsilon)-\nabla_A u\big)
 . \nabla (u_\epsilon-\epsilon-u)^+dx &\leqslant&
 \int_\Omega h(1-H_\epsilon(u_\epsilon-\psi))(u_\epsilon-\epsilon-u)^+dx\\
 &&=0.
\end{eqnarray*}

\n Using (2.7) and the fact that we have
$(u_\epsilon-\epsilon-u)^+=(g-\epsilon-g)^+=0$ on $\partial\Omega$, we conclude as above
that $(u_\epsilon-\epsilon-u)^+=0$ a.e. in $\Omega$. This means that
$u_\epsilon-\epsilon \geqslant u$ a.e. in $\Omega$.

\vs 0,3 cm\n Now we deduce from (2.9), that $u_\epsilon$
converges to $u$ in $L^\infty(\Omega)$. We also deduce from (2.8)
that $u_\epsilon\in C_{loc}^{1,\alpha}(\Omega)$ and bounded in each
$C^{1,\alpha}(\overline{\Omega'})$ for each
$\Omega'\subset\subset\Omega$ (see \cite{[L1]}). There exists
therefore a subsequence such that $u_\epsilon \rightarrow u$ in
$C^{1,\beta}(\overline{\Omega'})$ for each
$\Omega'\subset\subset\Omega$ and each $\beta<\alpha$.

\vs 0,3 cm\n Since we have
$$f-h\leqslant \Delta_A u_\epsilon\leqslant f\quad\text{a.e.
in}\quad\Omega,$$ we get for every $\varphi\in {\cal D}(\Omega)$,
$\varphi\geqslant 0$
\begin{eqnarray*}
 \int_\Omega -f\varphi dx\leqslant \int_\Omega \nabla_A u_\epsilon
 . \nabla \varphi dx \leqslant
 \int_\Omega (-f+h)\varphi dx.
\end{eqnarray*}
Letting $\epsilon\rightarrow 0$ and using the $C_{loc}^{1,\beta}$
convergence of $u_\epsilon$ to $u$, we get
\begin{eqnarray*}
 \int_\Omega -f\varphi dx\leqslant \int_\Omega \nabla_A u
 . \nabla \varphi dx \leqslant
 \int_\Omega (-f+h)\varphi dx.
\end{eqnarray*}
Hence we obtain
$$f-(f-\Delta_A\psi)^+\leqslant \Delta_A u\leqslant f\quad\text{a.e.
in}\quad\Omega.$$ \qed

\vs 0,5 cm\n Let us now give a definition that will be needed in the next
proposition and later on.
\begin{definition}\label{d3.1}
Let $p>1$ and $w\in W^{1,p}(\Omega)$. We say that $w(x)>0$ in the
sense of $W^{1,p}(\Omega)$ if there exists a neighborhood
$\vartheta_x\subset\Omega$ of $x$ and a nonnegative function
$\zeta\in W^{1,\infty}(\Omega)$ such that $\zeta(x)>0$ and
$w\geqslant\zeta$ a.e. in $\vartheta_x$. This definition is clearly
independent of $\zeta$ and one can verify that the set $\{w>0\}$ is an
open set of $\Omega$. Similarly we define $\{w<0\}$.
 \end{definition}

\begin{proposition}\label{p2.2} Let $u$ be the entropy solution of
the obstacle problem associated with $(f,\psi,g)$. Assume that we have
$u\in W^{1,p}(\Omega)$ for some $p>1$. Then we have
$$\Delta_Au =f\quad\text{a.e.
in}\quad \{u>\psi\}.$$
\end{proposition}

\vs 0,5cm\n\emph{Proof.} According to Definition 2.1, we denote by $\Lambda$ the open set
$\{u>\psi\}$ and we argue as in \cite{[RSU]}. Let $\varphi\in {\cal D}(\Lambda)$, $s>|g|_\infty,
|\psi|_\infty$, and $\epsilon>0$ small enough so that
$v=T_s(u)\pm\epsilon\varphi\in K_{\psi, g}\cap L^\infty(\Omega)$.
Using $v$ as a test function for the obstacle problem associated
with $(f,\psi,g)$, we obtain
\begin{eqnarray*}
\int_\Omega \nabla_A u
 . \nabla (T_t(T_s(u)\pm\epsilon\varphi-u))dx \geqslant
 \int_\Omega -f T_t(T_s(u)\pm\epsilon\varphi-u)dx
\end{eqnarray*}
which can be written as
\begin{eqnarray*}
\pm\epsilon\int_{\{|T_s(u)\pm\epsilon\varphi-u|<t\}} \nabla_A u
 . \nabla \varphi dx \geqslant
 \int_\Omega -f T_t(T_s(u)\pm\epsilon\varphi-u)dx.
\end{eqnarray*}
Choosing $t>\epsilon|\varphi|_\infty$ and letting
$s\rightarrow\infty$, we get
\begin{eqnarray*}
\pm\epsilon\int_{\{|\epsilon\varphi|<t\}} \nabla_A u
 . \nabla \varphi dx \geqslant
 \int_\Omega -f T_t(\pm\epsilon\varphi)dx=\pm\epsilon\int_\Omega -f \varphi dx.
\end{eqnarray*}
Hence we get
\begin{eqnarray*}
\int_\Lambda \nabla_A u
 . \nabla \varphi dx=
 \int_\Omega -f \varphi dx
\end{eqnarray*}
or $\Delta_A u=f\quad\text{a.e. in}\quad \Lambda.$ \qed
 \vs 0,5 cm

\begin{proposition}\label{p2.3}
Let the assumptions of Proposition 2.1  hold for two sets
of data $f_1, \psi, g)$ and $f_2, \psi, g)$ and denote by $u_i$,
$\xi_i=\Delta_A u_i-f_i$, $i=1, 2$, their respective variational solutions.
Assume moreover that $\psi\in W^{2,\infty}_{loc}(\Omega)$.
Then we have
\begin{equation}\label{2.16}
\int_\Omega|\xi_1-\xi_2|dx\leqslant \int_\Omega|f_1-f_2|dx.
\end{equation}
\end{proposition}

\vs 0,5 cm\n We need two lemmas for the variational solution $u_\epsilon$
of the $A_\epsilon$-obstacle problem associated with
$(f,\psi,g)$, where $A_\epsilon$ is associated to the function
$a_\epsilon(t)=\displaystyle{{ta(\sqrt{\epsilon+t^2})}\over{\sqrt{\epsilon+t^2}}}$.
Note that $a_\epsilon$ satisfies (1.5) with the same constants $a_0$, $a_1$.
In the first lemma, we establish a regularity result for $u_\epsilon$.

\begin{lemma}\label{l2.1} Assume that $f\in L_{loc}^\infty(\Omega)$ and
$\psi\in W^{2,\infty}_{loc}(\Omega)$. If $u_\epsilon$ is
the variational solution of the $A_\epsilon$-obstacle problem associated with
$(f,\psi,g)$, then we have $u_\epsilon\in W^{2,2}_{loc}(\Omega)$. In particular we have
$\nabla_{A_\epsilon} u_\epsilon\in W_{loc}^{1,2}(\Omega)$.
\end{lemma}

\vs 0,5cm\n\emph{Proof.} Note that $u_\epsilon$ satisfies the Lewy-Stampacchia inequalities
proven in Proposition 2.1
\begin{equation}\label{e2.15}
f-(f-\Delta_{A_\epsilon}\psi)^+\leqslant \Delta_{A_\epsilon} u_\epsilon\leqslant f\quad\text{a.e.
in}\quad\Omega.
\end{equation}

\n Moreover, since $\psi\in W^{2,\infty}_{loc}(\Omega)$, one  has

\begin{equation*}
\triangle_{A_\epsilon}\psi= {{a_\epsilon(|\nabla\psi|)}\over{|\nabla\psi|}}\Big[\triangle
\psi+\Big({{|\nabla\psi|a_\epsilon'(|\nabla\psi|)}\over{a_\epsilon(|\nabla\psi|)}}-1\Big){{D^2
\psi. \nabla \psi}\over{|\nabla\psi|^2}}.\nabla \psi\Big].
\end{equation*}

\n Using (1.5), we obtain for any ball $B\subset\subset\Omega$ and $x\in B$
\begin{eqnarray*}
|\triangle_{A_\epsilon}\psi (x)|&\leqslant& {{a((\epsilon+|\nabla\psi(x)|^2)^{1/2})}\over{(\epsilon+|\nabla\psi(x)|^2)^{1/2}}}[|\triangle
\psi(x)|+(a_1-1)|D^2 \psi(x)|]\nonumber\\
&\leqslant& a_1{{a((\epsilon+|\nabla\psi|_{\infty,B}^2)^{1/2})}\over{\epsilon}}|D^2 \psi|_{\infty,B}.
\end{eqnarray*}

\n It follows that $\triangle_{A_\epsilon} \psi\in L^\infty_{loc}(\Omega)$. Taking into account
(2.15), we deduce that $\triangle_{A_\epsilon} u_\epsilon\in L^\infty_{loc}(\Omega)$, which
in turn leads to $u_\epsilon\in W^{2,2}_{loc}(\Omega)$ (see \cite{[CL3]}).
Now since
\begin{equation*}
D_i(\nabla_{A_\epsilon} u_\epsilon)= {{a_\epsilon(|\nabla u_\epsilon|)}\over{|\nabla u_\epsilon|}}\Big[\nabla
u_{\epsilon x_i}+\Big({{|\nabla u_\epsilon|a_\epsilon'(|\nabla u_\epsilon|)}\over{a_\epsilon(|\nabla u_\epsilon|)}}-1\Big){{\nabla
u_{\epsilon x_i}.  \nabla u_\epsilon}\over{|\nabla u_\epsilon|^2}}.\nabla u_\epsilon\Big],
\end{equation*}
\n we get by (1.5)
\begin{equation*}
|D_i(\nabla_{A_\epsilon} u_\epsilon)|\leqslant a_1{{a((\epsilon+|\nabla u_\epsilon(x)|^2)^{1/2})}\over{(\epsilon+|\nabla u_\epsilon(x)|^2)^{1/2}}}|\nabla
u_{\epsilon x_i}|.
\end{equation*}
\n Hence $\nabla_{A_\epsilon} u_\epsilon\in W_{loc}^{1,2}(\Omega)$.
\qed

\vs 0,5 cm

\begin{lemma}\label{l2.2}
Let the assumptions of Lemma 2.1  hold for two sets
of data $f_1, \psi, g)$ and $f_2, \psi, g)$ and denote by $u_\epsilon^i$,
$\xi_\epsilon^i=\Delta_A u_\epsilon^i-f_i$, $i=1, 2$, their respective variational solutions.
Then we have
\begin{equation}\label{2.16}
\int_\Omega|\xi_\epsilon^1-\xi_\epsilon^2|dx\leqslant \int_\Omega|f_1-f_2|dx.
\end{equation}
\end{lemma}

\vs 0,5cm\n\emph{Proof.}  By Lemma 2.1, we have $\Delta_{A_\epsilon} u_\epsilon^i\in L^\infty(\Omega)$.
Integrating by parts and using the monotonicity of $\displaystyle{{{a((\epsilon+|\xi|^2)^{1/2})}\over{(\epsilon+|\xi|^2)^{1/2}}}\xi}$, we get
for $\delta>0$
\begin{eqnarray*}
&&\int_\Omega (\Delta_{A_\epsilon} u_\epsilon^1-\Delta_{A_\epsilon} u_\epsilon^2).(T_\delta(u_\epsilon^1-u_\epsilon^2))dx
= -\int_\Omega (\nabla_{A_\epsilon} u_\epsilon^1-\nabla_{A_\epsilon} u_\epsilon^2).\nabla(T_\epsilon(u_\epsilon^1-u_\epsilon^2))dx\nonumber\\
&&\hskip 4cm =-\int_\Omega T'_\delta(u_1-u_2)(\nabla_{A_\epsilon} u_\epsilon^1-\nabla_{A_\epsilon} u_\epsilon^2).(\nabla u_\epsilon^1-\nabla u_\epsilon^2)dx\leqslant 0.
\end{eqnarray*}

\n Letting $\delta\rightarrow 0$, we get for $T_0(s)=sign(s)=s/|s|$ if $s\neq0$, and $T_0(0)=0$
\begin{equation}\label{2.17}
\int_\Omega (\Delta_{A_\epsilon} u_\epsilon^1-\Delta_{A_\epsilon} u_\epsilon^2).(T_0(u_\epsilon^1-u_\epsilon^2))dx\leqslant 0.
\end{equation}

\n Note that
\begin{equation}\label{2.18}
(\xi_\epsilon^2-\xi_\epsilon^1)(u_\epsilon^2-u_\epsilon^1)\geqslant 0\quad\text{a.e. in}~~\Omega.
\end{equation}

\n Indeed since $\{u_\epsilon^2>u_\epsilon^1\}\subset\{u_\epsilon^2>\psi\}$, we have
$\xi_\epsilon^2=0$ a.e. in $\{u_\epsilon^2>u_\epsilon^1\}$. It follows that
$(\xi_\epsilon^2-\xi_\epsilon^1)(u_\epsilon^2-u_\epsilon^1)=-\xi_\epsilon^1(u_\epsilon^2-u_\epsilon^1)\geqslant 0$ a.e. in
$\{u_\epsilon^1>u_\epsilon^2\}$, since $\xi_\epsilon^1\leqslant 0$ a.e. in $\Omega$ by Proposition 2.1.
 Similarly, we show that  $(\xi_\epsilon^2-\xi_\epsilon^1)(u_\epsilon^2-u_\epsilon^1)\geqslant 0$ a.e. in
$\{u_\epsilon^2<u_\epsilon^1\}$. Hence (2.18) holds.

\n We know from Lemma 2.1 that we have $\nabla_{A_\epsilon} u_\epsilon^i\in W_{loc}^{1,2}(\Omega)$.
Since $\nabla_{A_\epsilon} u_\epsilon^2=\nabla_{A_\epsilon} u_\epsilon^1$ a.e. in $\{u_\epsilon^2=u_\epsilon^1\}$, we obtain
$div(\nabla_{A_\epsilon} u_\epsilon^2)=div(\nabla_{A_\epsilon} u_\epsilon^1)$ a.e. in $\{u_\epsilon^2=u_\epsilon^1\}$. Hence we
have
\begin{equation*}
\xi_\epsilon^2-\xi_\epsilon^1=f_1-f_2+\Delta_{A_\epsilon} u_\epsilon^2-\Delta_{A_\epsilon} u_\epsilon^1=f_1-f_2\quad\text{a.e. in}~~ \{u_\epsilon^2=u_\epsilon^1\}.
\end{equation*}

\n Taking into account (2.18), we obtain
\begin{eqnarray}\label{2.19}
&&\int_\Omega|\xi_\epsilon^1-\xi_\epsilon^2|dx=\int_{\{u_\epsilon^2=u_\epsilon^1\}}|\xi_\epsilon^1-\xi_\epsilon^2|dx
+\int_{\{u_\epsilon^2>u_\epsilon^1\}}|\xi_\epsilon^1-\xi_\epsilon^2|dx
+\int_{\{u_\epsilon^2<u_\epsilon^1\}}|\xi_\epsilon^1-\xi_\epsilon^2|dx\nonumber\\
&&\qquad=\int_{\{u_\epsilon^2=u_\epsilon^1\}}|f_1-f_2|dx+\int_{\{u_\epsilon^2>u_\epsilon^1\}}(\xi_2-\xi_1)dx
+\int_{\{u_\epsilon^2<u_\epsilon^1\}}(\xi_1-\xi_2)dx\nonumber\\
&&\qquad=\int_{\{u_\epsilon^2=u_\epsilon^1\}}|f_1-f_2|dx+\int_{\{u_\epsilon^2>u_\epsilon^1\}}(f_1-f_2+\Delta_{A_\epsilon} u_\epsilon^2-\Delta_{A_\epsilon} u_\epsilon^1)dx\nonumber\\
&&\qquad+\int_{\{u_\epsilon^2<u_\epsilon^1\}}(f_2-f_1+\Delta_{A_\epsilon} u_\epsilon^1-\Delta_{A_\epsilon} u_\epsilon^2)dx\nonumber\\
&&\qquad=\int_{\{u_\epsilon^2=u_\epsilon^1\}}|f_1-f_2|dx+\int_{\{u_\epsilon^2>u_\epsilon^1\}}(f_1-f_2)dx
+\int_{\{u_\epsilon^2<u_\epsilon^1\}}(f_2-f_1)dx\nonumber\\
&&\qquad+\int_{\{u_\epsilon^2>u_\epsilon^1\}}(\Delta_{A_\epsilon} u_\epsilon^2-\Delta_{A_\epsilon} u_\epsilon^1)dx
+\int_{\{u_\epsilon^2<u_\epsilon^1\}}(\Delta_{A_\epsilon} u_\epsilon^1-\Delta_{A_\epsilon} u_\epsilon^2)dx.
\end{eqnarray}

\n Taking into account (2.17) and (2.19), we obtain
\begin{eqnarray*}
\int_\Omega|\xi_\epsilon^1-\xi_\epsilon^2|dx&\leqslant& \int_\Omega |f_1-f_2|dx+\int_\Omega (\Delta_{A_\epsilon} u_\epsilon^1-\Delta_{A_\epsilon} u_\epsilon^2).(T_0(u_\epsilon^1-u_\epsilon^2))dx\\
&&\leqslant \int_\Omega |f_1-f_2|dx.
\end{eqnarray*}
\qed

\vs 0,5cm\n\emph{Proof of Proposition 2.3.} Let $u_1, u_2$ be two variational
solutions of the $A-$obstacle problem associated with $(f_1,\psi, g)$
and $(f_2,\psi, g)$. We denote by $u_{\epsilon}^1, u_{\epsilon}^2$ the variational solutions
of the $A_{\epsilon}-$obstacle problem associated with $(f_1,\psi, g)$
and $(f_2,\psi, g)$ respectively. We observe that we can prove, as in the proof of Lemma 2.3 of
\cite{[CL3]}, that there exists a subsequence $(\epsilon_n)$ such that $u_{\epsilon_n}^1$
and $u_{\epsilon_n}^2$ converge in $W^{1,A}(\Omega)$ respectively to
$u_1$ and $u_2$.

\n Let $\xi_{\epsilon_n}^i=\Delta_{A_{\epsilon_n}} u_{\epsilon_n}^i-f_i$, $i=1, 2$.
Then we have by Lemma 2.2
\begin{equation}\label{2.20}
\int_\Omega|\xi_{\epsilon_n}^1-\xi_{\epsilon_n}^2|dx\leqslant \int_\Omega|f_1-f_2|dx.
\end{equation}

\n Let $O$ be an open set such that $O\subset\subset\Omega$. Let $\delta>0$,
and consider the open set $O_\delta=\{x\in O~:~|\nabla\psi(x)|>\delta\}$.
Note that we can assume, without loss of generality, that we have $\epsilon_n<\delta$ for all $n\geqslant 1$.
Using (1.5) and arguing as in the proof of Lemma 2.1, we obtain for $n\geqslant 1$ and $x\in O_\delta$
\begin{eqnarray*}
|\triangle_{A_{\epsilon_n}}\psi (x)|&\leqslant& {{a((\epsilon_n+|\nabla\psi(x)|^2)^{1/2})}\over{(\epsilon_n+|\nabla\psi(x)|^2)^{1/2}}}[|\triangle
\psi(x)|+(a_1-1)|D^2 \psi(x)|]\nonumber\\
&\leqslant& a_1{{a((\delta+|\nabla\psi|^2)^{1/2})}\over{|\nabla\psi(x)|}}|D^2 \psi|_{\infty,O}\nonumber\\
&\leqslant& a_1{{a(\sqrt{2}|\nabla\psi|_{\infty,B})}\over{\delta}}|D^2 \psi|_{\infty,O}.
\end{eqnarray*}

\n It follows that $\triangle_{A_{\epsilon_n}} \psi$ is bounded uniformly in $L^\infty(O_\delta)$.
Taking into account (2.15), we deduce that $\triangle_{A_{\epsilon_n}} u_{\epsilon_n}^i$ are bounded
uniformly in $L^\infty(O_\delta)$.

\n Consider now the measurable set $O_0=\{x\in O~:~\nabla\psi(x)=0\}$. Since $\psi\in W^{2,\infty}(O)$, we have
$\Delta_{A_{\epsilon_n}}\psi=0$ a.e. in $O_0$. We deduce from (2.15) that
$-f_i^-\leqslant\Delta_{A_{\epsilon_n}} u_{\epsilon_n}^i\leqslant f_i$ a.e. in $O_0$.
Therefore $\triangle_{A_{\epsilon_n}} u_{\epsilon_n}^i$ are bounded uniformly
in $L^\infty(O_0)$.

\n We conclude that $\triangle_{A_{\epsilon_n}} u_{\epsilon_n}^i$ are bounded uniformly
in $L^\infty(O_0\cup O_\delta)$. Since moreover $\nabla_{A_{\epsilon_n}} u_{\epsilon_n}^i$ converges
to $\nabla_A u_i$ in $L^{\widetilde{A}}(\Omega)$, we deduce that we have up to a subsequence
\begin{eqnarray*}
\Delta_{A_{\epsilon_n}} u_{\epsilon_n}^i~\rightharpoonup~ \Delta_A u^i\quad\text{weakly-*} \text{ in}~~L^\infty(O_0\cup O_\delta),
\end{eqnarray*}
and therefore
\begin{equation*}
\Delta_{A_{\epsilon_n}} u_{\epsilon_n}^i~\rightharpoonup~ \Delta_A u^i\quad\text{in}~~L^1(O_0\cup O_\delta).
\end{equation*}

\n It follows that
\begin{equation}\label{2.21}
\xi_{\epsilon_n}^i~\rightharpoonup~ \xi_i\quad\text{in}~~L^1(O_0\cup O_\delta).
\end{equation}

\n Using the semicontinuity of the norm $||.||_{L^1(O_0\cup O_\delta)}$, we get from (2.20)-(2.21), since
$\displaystyle{\int_{O_0\cup O_\delta}|\xi_{\epsilon_n}^1-\xi_{\epsilon_n}^2|dx\leqslant \int_\Omega |\xi_{\epsilon_n}^1-\xi_{\epsilon_n}^2|dx}$
\begin{eqnarray*}
\int_{O_0\cup O_\delta}|\xi_1-\xi_2|dx\leqslant \liminf_{n\rightarrow\infty}\int_{O_0\cup O_\delta}|\xi_{\epsilon_n}^1-\xi_{\epsilon_n}^2|dx\leqslant \int_\Omega|f_1-f_2|dx.
\end{eqnarray*}

\n Given that $O$ and $\delta$ are arbitrary, we obtain
\begin{eqnarray*}
\int_\Omega|\xi_1-\xi_2|dx\leqslant \int_\Omega|f_1-f_2|dx.
\end{eqnarray*}

\qed

\section{Lemmas on Entropy Solutions}\label{S3}

\vs 0,5 cm\n Since variational solutions, in particular those obtained with
$f\in L^\infty(\Omega)$, are also entropy solutions, in this section we establish
several auxiliary results on convergence of sequences of entropy solutions when
$f_n\rightarrow f$ in $L^1(\Omega)$. Without loss of generality, we assume that
$||f_n||_1\leqslant ||f||_1+1. $We start with an a priori estimate.

\begin{lemma}\label{l3.1} Let $v_0\in K_{\psi,g}\cap L^\infty(\Omega)$, and
let $u$ be an entropy solution of the obstacle problem associated
with $(f,\psi,g)$. Then we have
$$ \int_\Omega A(|\nabla T_t(u)|)dx \leqslant (1+a_1)2^{1+a_1}\int_\Omega A(|\nabla v_0|)dx+2t||f||_1,\quad\forall t>0. $$
\end{lemma}

\vs 0,5cm\n \emph{Proof.} Using $v_0$ as a test function for (1.9), we get for $t>0$
\begin{equation*}
\int_{\{|v_0-u|<t\}} \nabla_A u
 . (\nabla v_0-\nabla u))dx \geqslant \int_\Omega  -fT_t(v_0-u)dx
\end{equation*}
or
\begin{eqnarray}\label{3.1}
\int_{\{|v_0-u|<t\}} |\nabla u|a(|\nabla u|)dx&\leqslant&
\int_{\{|v_0-u|<t\}} \nabla_A u
.\nabla v_0 dx +\int_\Omega  fT_t(v_0-u)dx\nonumber\\
&\leqslant& \int_{\{|v_0-u|<t\}} a(|\nabla u|).|\nabla v_0| dx.
+t||f||_1.
\end{eqnarray}

\n Using inequalities (2.2)-(2.3), we get
\begin{eqnarray}\label{3.2}
\int_{\{|v_0-u|<t\}} &&a(|\nabla u|).|\nabla v_0| dx\leqslant{1\over
2}\int_{\{|v_0-u|<t\}} |\nabla u|a(|\nabla u|)dx+ \int_{\{|v_0-u|<t\}}
|\nabla v_0|a(2|\nabla v_0|)dx\nonumber\\
&&\leqslant {1\over 2}\int_{\{|v_0-u|<t\}} |\nabla u|a(|\nabla u|)dx+
2^{a_1}\int_{\{|v_0-u|<t\}} |\nabla v_0|a(|\nabla v_0|)dx.
\end{eqnarray}

\n Combining (3.1) and (3.2), we obtain
\begin{equation}\label{3.3}
\int_{\{|v_0-u|<t\}} |\nabla u|a(|\nabla u|) dx\leqslant
2^{1+a_1}\int_{\{|v_0-u|<t\}} |\nabla v_0|a(|\nabla v_0|)dx+2t||f||_1.
\end{equation}

\n Replacing $t$ by $t+|v_0|_{\infty}$ in (3.3) and using the fact
that $\{|u|<t\} \subset\{|v_0-u|<t+|v_0|_{\infty}\}$, we obtain
\begin{equation*}
\int_{\{|u|<t\}} |\nabla u|a(|\nabla u|) dx\leqslant
2^{1+a_1}\int_\Omega |\nabla v_0|a(|\nabla v_0|)dx+2t||f||_1.
\end{equation*}

\n Taking into account (2.1), we get the result.
\qed

\vs 0,5 cm\n In the rest of this section, let $(u_n)$ be a sequence
of entropy solutions of the obstacle problem associated with $(f_n,\psi,g)$
and assume that
\begin{equation}\label{e3.4}
f_n\rightarrow f\quad \text{in}~~L^1(\Omega).
\end{equation}

\vs 0,5 cm

\begin{lemma}\label{l3.2} There exists a measurable function $u$ such that
\begin{eqnarray*}u _n &\rightarrow &u\quad\text{in measure,}\\
T_k(u_n) &\rightharpoonup &T_k(u)\quad\text{weakly in } W^{1,A}(\Omega)\\
T_k(u _n) &\rightarrow &T_k(u)\quad\text{strongly in }
L^A(\Omega)\text{ and a.e. in }\Omega.
\end{eqnarray*}
\end{lemma}

\vs 0,5cm\n \emph{Proof.} Let $s, t$ and $\epsilon$ be positive
numbers. It is easy to verify that we have for every $n, m \geqslant
1$
\begin{eqnarray}\label{3.5}
{\cal L}^N(\{|u_n-u_m|>s\})&\leqslant & {\cal L}^N(\{|u_n|>t\})+{\cal L}^N(\{|u_m|>t\})\nonumber\\
&&+{\cal L}^N(\{|T_t(u_n)-T_t(u_m)|>s\}).
\end{eqnarray}

\n Moreover we have
\begin{eqnarray}\label{3.6}
{\cal L}^N(\{|u_n|>t\})&=& {1\over{A(t)}}\int_{\{|u_n|>t\}}A(t)dx\nonumber\\
&\leqslant&{1\over{A(t)}}\int_\Omega A(|T_t(u_n)|)dx.
\end{eqnarray}

\n Since $v_0=g+(\psi-g)^+\in K_{\psi,g}\cap L^\infty(\Omega)$, we
deduce from Lemma 3.1 that
\begin{equation}\label{3.7}
\int_\Omega A(|\nabla T_t(u_n)|)dx=\int_{\{|u|<t\}}A(|\nabla u_n|)dx
\leqslant (1+a_1)2^{1+a_1}\int_\Omega A(|\nabla v_0|)dx+2t||f||_1.
\end{equation}

\n Using the convexity of $A$, Poincar\'{e}'s inequality, the fact
that $T_t(u_n)-T_t(g)\in W_0^{1,A}(\Omega)$, and (2.4), we deduce
from (3.6)-(3.7) for every $t>|g|_\infty$ and for some positive constant
$C$ independent of $n$ and $t$
\begin{eqnarray*}
{\cal L}^N(\{|u_n|>t\})&\leqslant&
{C\over{A(t)}}\int_\Omega (A(|\nabla g|)+A(|\nabla \psi|))dx+{{2t}\over{A(t)}}||f||_1\nonumber\\
&\leqslant&{{Ct+C}\over{A(t)}}.
\end{eqnarray*}

\n We deduce that there exists $t_\epsilon>0$ such that
\begin{equation}\label{3.8}
{\cal L}^N(\{|u_n|>t\})<{\epsilon\over 3}\qquad \forall t\geqslant
t_\epsilon,~~\forall n\geqslant 1.
\end{equation}

\n Now we have as in (3.6)
\begin{equation}\label{3.9}
{\cal L}^N(\{|T_{t_\epsilon}(u_n)-T_{t_\epsilon}(u_m)|>s\})\leqslant{1\over{A(s)}}\int_\Omega
A(|T_{t_\epsilon}(u_n)-T_{t_\epsilon}(u_m)|)dx.
\end{equation}

\n Using (3.7) again, we see that $(T_{t_{\epsilon}}(u_n))$ is a
bounded sequence in $W^{1,A}(\Omega)$. It follows that up to a
subsequence $(T_{t_{\epsilon}}(u_n))$ converges strongly in
$L^A(\Omega)$. Taking into account (3.9), there exists
$n_0=n_0(t_\epsilon,s)\geqslant 1$ such that
\begin{equation}\label{3.10}
{\cal L}^N(\{|T_{t_\epsilon}(u_n)-T_{t_\epsilon}(u_m)|>s\})<{\epsilon\over 3}\quad \forall n,
m\geqslant n_0.
\end{equation}

\n Combining (3.5), (3.8) and (3.10), we obtain
\begin{equation*}
{\cal L}^N(\{|u_n-u_m|>s\})<\epsilon\qquad \forall n, m\geqslant n_0.
\end{equation*}
\n Hence $(u_n)$ is a Cauchy sequence in measure, and therefore
there exists a measurable function $u$ such that $u
_n\,\rightarrow\,u$ in measure. The remainder of the lemma is a
consequence of the fact that $(T_k(u_n))$ is a bounded sequence in
$W^{1,A}(\Omega)$.  \qed

\vs 0,5 cm
\begin{lemma}\label{l3.3} Let $u$ be an entropy solution of the
obstacle problem associated with $(f,\psi,g)$. Then we have for every $t>0$
$$\lim_{h\rightarrow \infty}\int_{\{h\leqslant|u|\leqslant h+t\}}A(|\nabla u|)dx =0. $$
\end{lemma}

\vs 0,5cm\n \emph{Proof.} Let $t, h>0$ and $u$ be an entropy
solution of the obstacle problem associated with $(f,\psi,g)$. For
$h\geqslant |g|_\infty, |\psi|_\infty$, it is easy to check that
$T_h(u)$ is a test function for (1.9). We obtain
\begin{equation*}
\int_\Omega \nabla_A u
 . \nabla (T_t(u-T_h(u)))dx \leqslant \int_\Omega  -fT_t(u-T_h(u))dx
\end{equation*}
or
\begin{eqnarray*}
\int_{\{h\leqslant|u|\leqslant h+t\}}  |\nabla u|a(|\nabla
u|)dx&\leqslant& t\int_{\{|u|>h\}}  |f|dx~~\rightarrow~~0,\\
&&~~\text{as}~~h\rightarrow \infty.
\end{eqnarray*}

\n We conclude by taking into account (2.1).
\qed

\vs 0,5cm
\begin{proposition}\label{p3.1} There exists a subsequence of $(u_n)$ and a
measurable function $u$ such that for each $p\in\big(1,{N\over{N-1}}a_0\big)$,
we have
\begin{eqnarray*}u _n  &\rightarrow & u\quad\text{strongly in}\quad W^{1,p}(\Omega).
\end{eqnarray*}

\n If moreover $a_1<{N\over{N-1}}a_0$, then
\begin{eqnarray*}
\nabla_A u _n  &\rightarrow & \nabla_A u\quad\text{strongly in}\quad L^1(\Omega).
\end{eqnarray*}
\end{proposition}

\vs 0,5cm\n To prove Proposition 3.1, we need two preliminary lemmas.

\vs 0.5cm
\begin{lemma}\label{l3.4} There exists a subsequence of $(u_n)$ such that for each
$p\in\big(1,{N\over{N-1}}a_0\big)$, we have
\begin{eqnarray*}u _n  &\rightharpoonup & u\quad\text{weakly in } W^{1,p}(\Omega)\\
u _n  &\rightarrow & u\quad\text{strongly in } L^p(\Omega).
\end{eqnarray*}
\end{lemma}
\vs 0,5cm\n\emph{Proof.} Let $k>0$ and $n\geqslant 1$ and define $D_k=\{|u_n|\leqslant k\}$ and $B_k=\{k\leqslant|u_n|<k+1\}$. Using
Lemma 3.1 with $v_0=g+(\psi-g)^+$, we obtain
\begin{equation}\label{3.11}
 \int_{D_k} A(|\nabla u_n|)dx \leqslant 2^{1+a_1}(1+a_1)\int_\Omega A(|\nabla
 v_0|)dx+2k||f||_1.
\end{equation}

\n Using the function $T_k(u_n)$ for $k>|g|_\infty, |\psi|_\infty$, as a test function for the problem
associated with $(f_n,\psi,g)$, we obtain
\begin{equation*}
\int_\Omega\nabla_A u_n
 . \nabla(T_1(u_n-T_k(u_n)))dx \leqslant \int_\Omega  -f_nT_1(u_n-T_k(u_n))dx
\end{equation*}
or
\begin{eqnarray*}
\int_{B_k} |\nabla u_n|a(|\nabla u_n|)dx\leqslant|f_n|_1\leqslant||f||_1+1.
\end{eqnarray*}

\n Using inequality (2.1), we get
\begin{eqnarray}\label{3.12}
\int_{B_k} A(|\nabla u_n|)dx\leqslant||f||_1+1.
\end{eqnarray}

\n Now let $p\in\big(1,{N\over{N-1}}a_0\big)$ and consider the
function $A_p(t)=A(t^{1/p})$. It is easy to verify that the function
$a_p(t)=A_p'(t)={1\over p}t^{{1\over p}-1}a(t^{1/p})$ satisfies
the inequalities (1.5) with constants $b_0, b_1$, i.e.
\begin{equation*}
 b_0\leqslant { t a_p^\prime (t) \over a_p(t)}  \leqslant b_1  \qquad
 \forall t>0, \quad \text{with}\quad b_0={{a_0+1-p}\over p},~~ b_1={{a_1+1-p}\over p}.
\end{equation*}

\n Note that because $a_0<N-1$, we have $p<{N\over{N-1}}a_0=a_0+{1\over{N-1}}a_0<a_0+1$.
We also have $p<{N\over{N-1}}a_0<N$, since $a_0<N-1$.

\n Applying H\"{o}lder's inequality (see \cite{[A]}), we obtain
\begin{equation}\label{3.13}
 \int_{B_k} |\nabla u_n|^p dx \leqslant 2 \big||\nabla u_n|^p\big|_{L^{A_p}(B_k)}.
 \big|1\big|_{L^{\widetilde{A}_p}(B_k)}.
\end{equation}
From (3.12), we deduce that
\begin{equation}\label{3.14}
 \int_{B_k} A_p(|\nabla u_n|^p) dx= \int_{B_k} A(|\nabla u_n|)dx \leqslant ||f||_1+1.
\end{equation}
Let $\lambda_*=\max\big(\big((1+b_1)(1+||f||_1)\big)^{1\over{1+b_0}},
\big((1+b_1)(1+||f||_1)\big)^{1\over{1+b_1}}\big)$. Using the
inequality (2.4) for $A_p$  and (3.14), we obtain
\begin{eqnarray*}
 \int_{B_k} A_p\big({{|\nabla u_n|^p}\over \lambda_*}\big) dx &\leqslant&
 (1+b_1)\max\big({1\over\lambda_*^{1+b_0}},{1\over\lambda_*^{1+b_1}}\big)
\int_{B_k} A_p(|\nabla u_n|^p) dx\\
&\leqslant&(1+b_1)(1+||f||_1)
\max\big({1\over\lambda_*^{1+b_0}},{1\over\lambda_*^{1+b_1}}\big)\leqslant 1.
\end{eqnarray*}
This leads to
\begin{equation}\label{3.15}
\big||\nabla u_n|^p\big|_{L^{A_p}(B_k)}\leqslant
\max\big(\big((1+b_1)(1+||f||_1)\big)^{1\over{1+b_0}},
\big((1+b_1)(1+||f||_1)\big)^{1\over{1+b_1}}\big).
\end{equation}
Next we evaluate
\begin{eqnarray}\label{3.16}
\big|1\big|_{L^{\widetilde{A}_p}(B_k)}&=&\inf\big\{\lambda>0~:~\int_{B_k}
\widetilde{A}_p\big({1\over \lambda}\big) dx \leqslant 1~\big\}\nonumber\\
 &=&\inf\big\{\lambda>0~:~
\widetilde{A}_p\big({1\over \lambda}\big) \leqslant {1\over {|B_k|}}~\big\}\nonumber\\
 &=& {1\over {\widetilde{A}^{-1}_p\big({1\over {|B_k|}}\big)}}.
\end{eqnarray}
Taking into account (3.13), (3.15) and (3.16), we obtain
\begin{equation*}
 \int_{B_k} |\nabla u_n|^p dx \leqslant {C_1\over {\widetilde{A}^{-1}_p\big({1\over {|B_k|}}\big)}}.
\end{equation*}
or
\begin{equation}\label{3.17}
{1\over {|B_k|}}\leqslant \widetilde{A}_p\Bigl({C_1\over{\int_{B_k}
|\nabla u_n|^p dx}}\Bigl).
\end{equation}
Let $p_*={{pN}\over{N-p}}$. Since $\lim_{k\rightarrow\infty}|B_k|=0,$ we deduce from (3.17) that
$\displaystyle{\lim_{k\rightarrow\infty}\int_{B_k} |\nabla
u_n|^p dx=0}$. It follows that there exists $k_0\geqslant 1$ such that
$$\int_{B_k} |\nabla
u_n|^p dx\leqslant 1,~~\forall k\geqslant k_0,~~\forall n\geqslant 1.$$

\n Since $\displaystyle{|B_k|\leqslant
{1\over{k^{p_*}}} \int_{B_k} |u_n|^{p_*} dx}$, we obtain from (3.17) and the
inequalities (2.5)-(2.6) applied to $a_p$, that
for all $k\geqslant k_0$
\begin{eqnarray*}
&&\Bigl({1\over{k^{p_*}}} \int_{B_k} |u_n|^{p_*}
dx\Bigl)^{-1}\leqslant\widetilde{A}_p\Bigl({C_1\over{\int_{B_k}
|\nabla u_n|^p dx}}\Bigl)\\
&&\qquad\leqslant\big(1+{1\over{b_0}}\big)\max\Biggl(\Bigl({1\over{\int_{B_k} |\nabla
u_n|^p dx}}\Bigl)^{{1\over b_0}+1},\Bigl({1\over{\int_{B_k}
|\nabla u_n|^p dx}}\Bigl)^{{1\over b_1}+1}\Biggl)\widetilde{A}_p(C_1)\\
&&\qquad\leqslant\big(1+{1\over{b_0}}\big)\Bigl({1\over{\int_{B_k} |\nabla
u_n|^p dx}}\Bigl)^{{1\over b_0}+1}\widetilde{A}_p(C_1).
\end{eqnarray*}

\n This leads for $\displaystyle{C_2=\Bigl(\big(1+{1\over{b_0}}\big)\widetilde{A}_p(C_1)\Bigl)^{{b_0}\over{1+b_0}}}$
to
\begin{eqnarray*}
\int_{B_k} |\nabla u_n|^p dx\leqslant
C_2
{1\over{k^{{p_*b_0}\over{1+b_0}}}}
\Bigl(\int_{B_k}
|u_n|^{p_*}dx\Bigl)^{{b_0}\over{1+b_0}}.
\end{eqnarray*}
Summing up over from $k=k_0$ to $k=K$ and using H\"{o}lder's inequality, we obtain
\begin{eqnarray}\label{3.18}
\sum_{k=k_0}^K\int_{B_k} |\nabla u_n|^p dx&\leqslant& C_2
\Bigl(\sum_{k=k_0}^K\int_{B_k}
|u_n|^{p_*}dx\Biggl)^{{b_0}\over{1+b_0}}
\Bigl(\sum_{k=k_0}^K{1\over{k^{p_*b_0}}}\Biggl)^{1\over{1+b_0}}.
\end{eqnarray}
Note that
\begin{equation}\label{3.19}
\int_{\{|u_n|\leqslant K\}}|\nabla u_n|^p dx=\int_{D_{k_0}} |\nabla u_n|^p dx+
\sum_{k=k_0}^K\int_{B_k} |\nabla u_n|^p dx.
\end{equation}
To estimate the first integral in the right hand side of (3.19), we
argue as above by using H\"{o}lder's inequality. We obtain
\begin{eqnarray}\label{3.20}
\int_{D_{k_0}} |\nabla u_n|^p dx&\leqslant& 2 \big||\nabla
u_n|^p\big|_{L^{A_p}(D_{k_0})}.
 \big|1\big|_{L^{\widetilde{A}_p}(D_{k_0})}\leqslant 2\lambda_*.{1\over
{\widetilde{A}^{-1}_p\big({1\over {|D_{k_0}|}}\big)}}\nonumber\\
&\leqslant& 2\lambda_*.{1\over {\widetilde{A}^{-1}_p\big({1\over
{|\Omega|}}\big)}}=C_0.
\end{eqnarray}

\n Note that the series
$\displaystyle{\sum_{k=k_0}^K{1\over{k^{p_*b_0}}}}$ converges since
$a_0>{N\over{N-1}}p$, and
$$p_*b_0={{pN}\over{N-p}}{{a_0+1-p}\over p}
={N\over{N-p}}(a_0+1-p)>{N\over{N-p}}({N\over{N-1}}p+1-p)>1.$$
Combining (3.18)-(3.20), we get for $k_0$ large enough
\begin{eqnarray}\label{3.21}
\int_{\{|u_n|\leqslant K\}}|\nabla u_n|^p dx&\leqslant&C_0+ C_2\Bigl(
\int_{\{|u_n|\leqslant K\}}|u_n|^{p_*} dx\Bigl)^{{b_0}\over{1+b_0}}\Biggl
(\sum_{k=k_0}^\infty{1\over{k^{p_*b_0}}}\Biggl)^{1\over{1+b_0}}\nonumber\\
&=&C_0+ C_3\Bigl(
\int_{\{|u_n|\leqslant K\}}|u_n|^{p_*} dx\Bigl)^{{b_0}\over{1+b_0}}.
\end{eqnarray}

\n Note that $T_K(u_n)\in W^{1,p}(\Omega)$. Indeed we have by
H\"{o}lder's inequality
\begin{equation*}
 \int_\Omega |\nabla T_K(u_n)|^p dx \leqslant
 2 \big||\nabla T_K(u_n)|^p\big|_{L^{A_p}(\Omega)}.
 \big|1\big|_{L^{\widetilde{A}_p}(\Omega)}
 ={2\over {\widetilde{A}^{-1}_p\big({1\over {|\Omega|}}\big)}}.
 \big||\nabla T_K(u_n)|^p\big|_{L^{A_p}(\Omega)}.
\end{equation*}
and
\begin{equation*}
 \int_\Omega A_p(|\nabla T_K(u_n)|^p) dx= \int_\Omega A(|\nabla T_K(u_n)|)dx <\infty.
\end{equation*}

\n Similarly we verify that $T_K(g)=g\in W^{1,p}(\Omega)$ for
$K>|g|_\infty$. Hence we have $T_K(u_n)-g\in W_0^{1,p}(\Omega)$.
Using the Sobolev embedding $W_0^{1,p}(\Omega)\subset
L^{p_*}(\Omega)$, Poincar\'{e}'s inequality, we obtain
\begin{eqnarray}\label{3.22}
|T_K(u_n)|_{L^{p_*}(\Omega)}^p&\leqslant&
2^{p-1}\big(|T_K(u_n)-g|_{L^{p_*}(\Omega)}^p+|g|_{L^{p_*}(\Omega)}^p\big)\nonumber\\
&\leqslant&
2^{p-1}\big(C|\nabla(T_K(u_n)-g)|_{L^p(\Omega)}^p+|g|_{L^{p_*}(\Omega)}^p\big)\nonumber\\
&\leqslant& 2^{2(p-1)}C|\nabla
T_K(u_n)|_{L^p(\Omega)}^p+2^{2(p-1)}C|\nabla g|_{L^p(\Omega)}^p
+2^{p-1}|g|_{L^{p_*}(\Omega)}^p\nonumber\\
&=&C_4\Big(\int_{\{|u_n|\leqslant K\}}|\nabla u_n|^p dx+1\Big).
\end{eqnarray}

\n Using the fact that
\begin{eqnarray}\label{3.23}
\int_{\{|u_n|\leqslant K\}}|u_n|^{p_*} dx \leqslant
\int_{\{|u_n|\leqslant K\}}|T_K(u_n)|^{p_*}
dx\leqslant|T_K(u_n)|_{L^{p_*}(\Omega)}^{p_*}
\end{eqnarray}
we obtain from (3.21)-(3.22), for $C_5=C_3C_4^{{p_*b_0}\over{p(1+b_0)}}$

\begin{eqnarray}\label{3.24}
\int_{\{|u_n|\leqslant K\}}|\nabla u_n|^p dx&\leqslant&
C_0+ C_3|T_K(u_n)|_{L^{p_*}(\Omega)}^{{p_*b_0}\over{1+b_0}}\nonumber\\
&\leqslant&
C_0+ C_5\Big(\int_{\{|u_n|\leqslant K\}}|\nabla u_n|^p dx+1\Big)^{{p_*b_0}\over{p(1+b_0)}}.
\end{eqnarray}

\n Now it is easy to verify that
$${{p_*b_0}\over{p(1+b_0)}}<1~\Leftrightarrow~a_0+1<N.$$

\n It follows from (3.24) that for $k_0$ large enough, the integral
$\displaystyle{\int_{\{|u_n|\leqslant K\}}|\nabla u_n|^p dx}$ is
bounded independently of $n$ and $K$. Using (3.22)-(3.23),
we deduce that $\displaystyle{\int_{\{|u_n|\leqslant K\}}| u_n|^{p_*}
dx}$ is also bounded independently of $n$ and $K$. Letting
$K\rightarrow\infty$, we deduce that $|\nabla u_n|_{L^{p}(\Omega)}$
and $|u_n|_{L^{p_*}(\Omega)}$ are uniformly bounded independently of
$n$. In particular $u_n$ is bounded in $W^{1,p}(\Omega)$. Therefore
there exists a subsequence of $(u_n)$ and a function $v\in
W^{1,p}(\Omega)$ such that
\begin{eqnarray*}
u_n &\rightharpoonup &v\quad\text{weakly in } W^{1,p}(\Omega)\\
u_n &\rightarrow &v\quad\text{strongly in } L^p(\Omega)\text{ and
a.e. in }\Omega.
\end{eqnarray*}
But since $u_n \rightarrow u$ in measure in $\Omega$, we have necessarily
$u=v$ and $u\in W^{1,p}(\Omega)$.

\qed

\vs 0.5cm
\begin{lemma}\label{l3.5} There exists a subsequence of $(u_n)$ such that for each $k>0$
\begin{eqnarray*}
\lim_{n,m\rightarrow\infty} \int_{\{|u_n|\leqslant k, |u_m|\leqslant k\}} \big(\nabla_A u_n-\nabla_A u_m\big)
 . \nabla (u_n-u_m)dx= 0.
\end{eqnarray*}
\end{lemma}
\vs 0,5cm\n\emph{Proof.} We argue as in \cite{[BGO]} for the proof of the uniqueness of the entropy
solution. So let $n, m \geqslant1$, $k>0$ and
$h\geqslant |g|_\infty, |\psi|_\infty$, with $k<h$. Let $u_n$ and $u_m$ be two entropy solutions
of the obstacle problem associated with $(f_n,\psi,g)$ and
$(f_m,\psi,g)$ respectively.  It is easy to check that
$T_h(u_m)$ and $T_h(u_n)$ are test functions for (1.9)
associated with $(f_n,\psi,g)$ and $(f_m,\psi,g)$ respectively. We obtain for $t>0$
\begin{equation*}
\int_\Omega \nabla_A u_n
 . \nabla (T_t(u_n-T_h(u_m)))dx \leqslant \int_\Omega  -f_nT_t(u_n-T_h(u_m))dx
\end{equation*}
\begin{equation*}
\int_\Omega \nabla_A u_m
 . \nabla (T_t(u_m-T_h(u_n)))dx \leqslant \int_\Omega  -f_mT_t(u_m-T_h(u_n))dx
\end{equation*}
Adding the two inequalities, we get
\begin{eqnarray}\label{3.25}
&&\int_{\{|u_n-T_h(u_m)|\leqslant t\}}  \nabla_A u_n
 . \nabla (u_n-T_h(u_m))dx\nonumber\\
&&\qquad+ \int_{\{|u_m-T_h(u_n)|\leqslant t\}} \nabla_A u_m
 . \nabla (u_m-T_h(u_n))dx\nonumber\\
&&\qquad\leqslant \int_\Omega -(f_nT_t(u_n-T_h(u_m))+f_mT_t(u_m-T_h(u_n)))dx.
\end{eqnarray}
Let us introduce the following sets as in \cite{[BGO]}
\begin{eqnarray*}
&&C_0^h=\{|u_n-u_m|\leqslant t,~|u_n|\leqslant h,~|u_m|\leqslant h\}\\
&&C_1^h=\{|u_n-T_h(u_m)|\leqslant t,~|u_m|> h\}\\
&&C_2^h=\{|u_n-T_h(u_m)|\leqslant t,~|u_m|\leqslant h,~|u_n|> h\}\\
&&C_3^h=\{|u_m-T_h(u_n)|\leqslant t,~|u_n|> h\}\\
&&C_4^h=\{|u_m-T_h(u_n)|\leqslant t,~|u_n|\leqslant h,~|u_m|> h\}.
\end{eqnarray*}
Then it is easy to see that $\{|u_n-T_h(u_m)|\leqslant t\}=C_0^h\cup
C_1^h\cup C_2^h$ and $\{|u_m-T_h(u_n)|\leqslant t\}=C_0^h\cup C_3^h\cup C_4^h$.
Moreover we have
\begin{eqnarray}\label{3.26}
&&\int_{C_0^h}  \nabla_A u_n
 . \nabla (u_n-T_h(u_m))dx+ \int_{C_0^h} \nabla_A u_m
 . \nabla (u_m-T_h(u_n))dx\nonumber\\
&&\qquad\qquad=\int_{C_0^h}  \big(\nabla_A u_n-\nabla_A u_m\big). \nabla (u_n-u_m)dx.
\end{eqnarray}

\begin{eqnarray}\label{3.27}
&&\int_{C_1^h}  \nabla_A u_n
 . \nabla (u_n-T_h(u_m))dx=\int_{C_1^h} \nabla_A u_n. \nabla u_n dx\geqslant 0\nonumber\\
&&\int_{C_3^h}  \nabla_A u_m
 . \nabla (u_m-T_h(u_n))dx=\int_{C_3^h} \nabla_A u_m. \nabla u_m dx\geqslant 0\nonumber\\
&&\int_{C_2^h}  \nabla_A u_n
 . \nabla (u_n-T_h(u_m))dx=\int_{C_2^h} \nabla_A u_n. \nabla(u_n-u_m) dx
 \geqslant -\int_{C_2^h} \nabla_A u_n. \nabla u_m dx\nonumber\\
&&\int_{C_4^h}  \nabla_A u_m
 . \nabla (u_m-T_h(u_n))dx=\int_{C_4^h} \nabla_A u_m. \nabla(u_m-u_n) dx
 \geqslant -\int_{C_4^h} \nabla_A u_m. \nabla u_n dx.\nonumber\\
 &&
\end{eqnarray}
We deduce from (3.25)-(3.27) that
\begin{eqnarray}\label{3.28}
&&\int_{C_0^h}  \big(\nabla_A u_n-\nabla_A u_m\big). \nabla (u_n-u_m)dx\leqslant
\int_\Omega  -(f_nT_t(u_n-T_h(u_m))+f_mT_t(u_m-T_h(u_n)))dx\nonumber\\
&&\qquad +\int_{C_2^h} |\nabla u_n| a(|\nabla u_m|) dx+\int_{C_4^h}
|\nabla u_m| a(|\nabla u_n|) dx.
\end{eqnarray}
Using  inequality (2.2), we obtain
\begin{eqnarray}\label{3.29}
&&\int_{C_2^h} |\nabla u_n| a(|\nabla u_m|) dx\leqslant \int_{C_2^h}
|\nabla u_n| a(|\nabla u_n|) dx
+\int_{C_2^h} |\nabla u_m| a(|\nabla u_m|) dx\nonumber\\
&&\leqslant \int_{\{h<|u_n|\leqslant h+t\}} |\nabla u_n| a(|\nabla
u_n|) dx +\int_{\{h-t<|u_m|\leqslant h\}} |\nabla u_m| a(|\nabla u_m|)
dx.
\end{eqnarray}
Similarly we have
\begin{eqnarray}\label{3.30}
&&\int_{C_4^h} |\nabla u_m| a(|\nabla u_n|) dx \leqslant \int_{C_4^h}
|\nabla u_n| a(|\nabla u_n|) dx
+\int_{C_4^h} |\nabla u_m| a(|\nabla u_m|) dx\nonumber\\
&&\quad\leqslant \int_{\{h<|u_m|\leqslant h+t\}} |\nabla u_m| a(|\nabla u_m|) dx
+\int_{\{h-t<|u_n|\leqslant h\}} |\nabla u_n| a(|\nabla u_n|) dx.
\end{eqnarray}
Now since $k<h$, we have $C_0^k\subset C_0^h$. Choosing $t=2k$, we obtain from (3.28),
by using (2.7) and taking into account (3.29)-(3.30)
\begin{eqnarray*}&&\int_{\{|u_n|\leqslant k, |u_m|
\leqslant k\}}\big(\nabla_A u_n-\nabla_A u_m\big). \nabla (u_n-u_m)dx\nonumber\\
&&\qquad\leqslant  \int_\Omega -(f_nT_t(u_n-T_h(u_m))+f_mT_t(u_m-T_h(u_n)))dx\nonumber\\
&&\qquad +\int_{\{h<|u_n|\leqslant h+t\}} |\nabla u_n| a(|\nabla u_n|) dx
+\int_{\{h-t<|u_n|\leqslant h\}} |\nabla u_n| a(|\nabla u_n|) dx\nonumber\\
&&\qquad +\int_{\{h<|u_m|\leqslant h+t\}} |\nabla u_m| a(|\nabla u_m|) dx
+\int_{\{h-t<|u_m|\leqslant h\}} |\nabla u_m| a(|\nabla u_m|) dx.
\end{eqnarray*}
Letting $h\rightarrow\infty$, and using Lemma 3.3 and Lebesgue's theorem, we get
\begin{eqnarray*}\int_{\{|u_n|\leqslant k, |u_m|\leqslant k\}}\big(\nabla_A u_n-\nabla_A u_m\big). \nabla (u_n-u_m)dx&\leqslant& \int_\Omega -(f_n-f_m)T_t(u_n-u_m)dx\\
&\leqslant& t|f_n-f_m|_1.
\end{eqnarray*}
We conclude by using the monotonicity inequality (2.7) and the convergence of $(f_n)$ to $f$ in $L^1(\Omega)$.

\qed

\vs 0.5cm\n \emph{Proof of Proposition 3.1.} We argue as in the proof of Lemma 3.2.
So let $s, k$ and $\epsilon$ be positive numbers. Then we have for every $n, m \geqslant
1$, $\{|\nabla u_n-\nabla u_m|>s\}\subset E_n^1\cup E_m^1 \cup E_n^2\cup E_m^2\cup E_{m,n}$
where $E_n^1=\{|u_n|>k\}$, $E_n^2=\{|\nabla u_n|>k\}$ and
$E_{n,m}=\{||\nabla u_n-\nabla u_m|>s, |u_n|\leqslant k, |u_m|\leqslant k, |\nabla u_n|\leqslant k,  |\nabla u_m|\leqslant k\}$.
\begin{eqnarray}\label{3.31}
&{\cal L}^N(\{|\nabla u_n-\nabla u_m|>s\})\leqslant  {\cal L}^N(E_n^1)+{\cal L}^N(E_m^1)+{\cal L}^N(E_n^2)+{\cal L}^N(E_m^2)+{\cal L}^N(E_{m,n}).\nonumber\\
&
\end{eqnarray}

\n Using the fact that by Lemma 3.4, the sequence  $(u_n)$ is uniformly bounded in $W^{1.p}(\Omega)$, we obtain

\begin{eqnarray*}
{\cal L}^N(E_n^1)\leqslant {1\over{k^p}}\int_\Omega |u_n|^pdx\leqslant
{C\over{k^p}}.
\end{eqnarray*}
\begin{eqnarray*}
{\cal L}^N(E_n^2)\leqslant {1\over{k^p}}\int_\Omega |\nabla u_n|^pdx\leqslant {C\over{k^p}}.
\end{eqnarray*}

\n We deduce that there exists $k_\epsilon>0$ such that
\begin{equation}\label{3.32}
\forall k\geqslant
k_\epsilon,~~\forall n, m\geqslant 1,  \quad {\cal L}^N(E_n^1)+{\cal L}^N(E_m^1)+{\cal L}^N(E_n^2)+{\cal L}^N(E_m^2)<{\epsilon\over 2}.
\end{equation}

\n Now using the inequality (2.7), we get for all $n,m \geqslant 1$

\begin{eqnarray}\label{e3.33}
&&\int_{E_{m,n}}|\nabla u_n-\nabla u_m|^2{ a\big((|\nabla u_n|^2+|\nabla u_m|^2)^{1/2}\big)\over {(|\nabla u_n|^2+|\nabla u_m|^2)^{1/2}}}dx\nonumber\\
&&\leqslant{1\over { C(A,N)}} \int_{E_{m,n}} \big(\nabla_A u_n-\nabla_A u_m\big)
 . (\nabla u_n-\nabla u_m)dx\nonumber\\
&&\leqslant{1\over { C(A,N)}} \int_{\{|u_n|\leqslant k, |u_m|\leqslant k\}} \big(\nabla_A u_n-\nabla_A u_m\big)
 . (\nabla u_n-\nabla u_m)dx.
\end{eqnarray}

\n Using the fact that we have in $E_{m,n}$
\begin{eqnarray*}
&&s<|\nabla u_n-\nabla u_m|\leqslant|\nabla u_n|+|\nabla u_m|\leqslant \sqrt{2}(|\nabla u_n|^2+|\nabla u_m|^2)^{1/2}\\
&&(|\nabla u_n|^2+|\nabla u_m|^2)^{1/2}\leqslant k\sqrt{2},
\end{eqnarray*}
we deduce from (3.33)

\begin{eqnarray*}
{\cal L}^N(E_{m,n})\leqslant{{k\sqrt{2}\over{s^2a(s/\sqrt{2})}C(A,N)}}\int_{\{|u_n|\leqslant k, |u_m|\leqslant k\}} \big(\nabla_A u_n-\nabla_A u_m\big)
 . (\nabla u_n-\nabla u_m)dx.
\end{eqnarray*}

\n Using Lemma 3.5, for $k=k_\epsilon$, there exists $n_0=n_0(s,k,\epsilon)\geqslant 1$ such that
\begin{equation}\label{3.34}
{\cal L}^N(E_{m,n})<{\epsilon\over 2}\quad \forall n,
m\geqslant n_0.
\end{equation}

\n Combining (3.31), (3.32) and (3.34), we obtain
\begin{equation*}
{\cal L}^N(\{|\nabla u_n-\nabla u_m|>s\})<\epsilon\qquad \forall n, m\geqslant n_0.
\end{equation*}
\n We deduce that  the sequence  $(\nabla u_n)$ is a Cauchy sequence in measure, and therefore
there exists a measurable function $V$ such that $(\nabla u_n)$
converges in measure to the measurable function $V$.
 We shall prove that $V=\nabla u$ and that $(\nabla u_n)$ converges
 strongly to $\nabla u$ in $L^p(\Omega)$ for each $p\in\big(1,{N\over{N-1}}a_0\big)$.
 To do that, we will apply Vitalli's Theorem, using the fact that by Lemma 3.4,
  $\nabla u_n $  is bounded in $L^p(\Omega)$ for each $p\in\big(1,{N\over{N-1}}a_0\big)$.

\n So let $r\in\big(p,{N\over{N-1}}a_0\big)$ and $O\subset\Omega$ be a measurable set. Then we have by H\"{o}lder's inequality
\begin{equation*}
 \int_O |\nabla u_n|^p dx \leqslant \Big(\int_O |\nabla u_n|^r dx\Big)^{p/r}.
 |O\big|^{{r-p}\over r}\leqslant C|O\big|^{{r-p}\over r}\rightarrow 0
\end{equation*}
uniformly in $n$, as $|O\big|\rightarrow 0$. We deduce that $\nabla
u_n$ converges strongly to $V$ in $L^p(\Omega)$. Using the fact
that $\nabla u_n$ converges weakly to $\nabla u$ in $L^p(\Omega)$,
we obtain $V=\nabla u$ and also that $(\nabla u_n)$ converges
strongly to $\nabla u$ in $L^p(\Omega)$ for each
$p\in\big(1,{N\over{N-1}}a_0\big)$.

\vs 0,2cm \n Now assume that $a_1<{N\over{N-1}}a_0$ and let $p\in\big(a_1,{N\over{N-1}}a_0\big)$.
Note that since $\nabla u_n$ converges to $\nabla u$ a.e. in $\Omega$, to
prove the convergence
$$ \nabla_A u_n\,\rightarrow\, \nabla_A u\quad\text{in } L^1(\Omega)$$
it suffices, Thanks to Vitalli's Theorem, to show that for every
measurable subset $O$ of $\Omega$, $\displaystyle{\int_O
|\nabla_A u_n|dx}$ converges to 0 uniformly in $n$, as
$|O\big|\rightarrow 0$. Note that $|\nabla_A u_n|\leqslant
a(|\nabla u_n|)$ and let $D(t)=(a^{-1}(t))^p$.  The function
$d(t)={{(a^{-1}(t))^p}\over t}$ satisfies the inequalities
\begin{equation*}
 d_0\leqslant { t d^\prime (t) \over d(t)}  \leqslant d_1  \qquad
 \forall t>0, \quad \text{with}\quad d_0={{p-a_1}\over a_1},~~ d_1={{p-a_0}\over a_0}.
\end{equation*}

\n It follows that we have by H\"{o}lder's inequality

\begin{equation*}
 \int_O a(|\nabla u_n|) dx \leqslant 2 \big|a(|\nabla u_n|)|\big|_{L^{D}(O)}.
 \big|1\big|_{L^{\widetilde{D}}(O)}.
\end{equation*}
Since we have also
$$\int_O D(a(|\nabla
u_n|))dx=\int_O |\nabla u_n|^p dx\leqslant C,$$
and
$$\big|1\big|_{L^{\widetilde{D}}(O)}={1\over {\widetilde{D}^{-1}\big({1\over {|O|}}\big)}}$$
we obtain
\begin{equation*}
 \int_O a(|\nabla u_n|) dx \leqslant {C\over {\widetilde{D}^{-1}\big({1\over {|O|}}\big)}}\rightarrow 0
\end{equation*}
uniformly in $n$, as $|O\big|\rightarrow 0$.

\qed

\vs 0,5cm

\begin{lemma}\label{l3.6} Let $(v_n)$ be a sequence of  measurable functions
uniformly bounded in $L^{\widetilde{A}}(\Omega)$ and converging in measure to a measurable function
$v$ in $\Omega$. Then $(v_n)$ converges strongly to $v$ in $L^1(\Omega)$.
\end{lemma}
\vs 0,5cm \n\emph{Proof.} Let $s, \epsilon$ be two positive numbers. First note that we have for each measurable subset $D$
of $\Omega$ (see the proof of Lemma 3.4) $\displaystyle{|\chi_{D}|_{L^A(\Omega)}={1\over{A^{-1}\big({1\over D}\big)}}}$.
Next we have by H\"{o}lder's inequality
\begin{eqnarray}\label{e3.35}
|v_n-v_m|_1&\leqslant&\int_{\{|v_n-v_m|\leqslant s\}} |v_n-v_m|dx+\int_{\{|v_n-v_m|>s\}} |v_n-v_m|dx\nonumber\\
&\leqslant& s|\Omega|+2|v_n-v_m|_{L^{\widetilde{A}}(\Omega)}.|\chi_{\{|v_n-v_m|>s\}}|_{L^A(\Omega)}\nonumber\\
&\leqslant& s|\Omega|+4C.|\chi_{\{|v_n-v_m|>s\}}|_{L^A(\Omega)}\nonumber\\
&\leqslant& s|\Omega|+{4C\over{A^{-1}\big({1\over |v_n-v_m|>s\}}\big)}}
\end{eqnarray}

\n Choosing $s={\epsilon\over {2|\Omega|}}$ and using the convergence in measure of  $(v_n)$ to $v$ in $\Omega$, we deduce that there exists $n_0\geqslant 1$ such that
\begin{equation}\label{e3.36}
{\cal L}^N(\{|v_n-v_m|>s\})<{1\over{A\big({{8C}\over \epsilon}\big)}}\qquad \forall n, m\geqslant n_0.
\end{equation}

\n Combining (3.35) and (3.36), we obtain $|v_n-v_m|_1<\epsilon$ for all $n, m\geqslant n_0$. Hence $(v_n)$ is a Cauchy sequence in $L^1(\Omega)$, and therefore it converges in $L^1(\Omega)$.
\qed

\vs 0,5cm

\begin{lemma}\label{l3.7} There exists a subsequence of $u_n$ such that for all $k>0$
\begin{eqnarray*}
&&\nabla_A T_k(u_n)~~ \rightarrow ~~ \nabla_A T_k(u)
\quad\text{strongly in } L^1(\Omega).
\end{eqnarray*}
\end{lemma}
\vs 0,5cm\n\emph{Proof.} Let $k$ be a positive number. First note that the sequence $(\nabla_A T_k(u_n))$
is uniformly bounded in $L^{\widetilde{A}}(\Omega)$. Indeed we have by the inequalities (2.1), (2.6), and
Lemma 3.1
\begin{eqnarray*}
\int_\Omega \widetilde{A}(|\nabla_A T_k(u_n)|)dx &\leqslant & \int_\Omega \widetilde{A}(|a(|\nabla T_k(u_n))|)dx\nonumber\\
&\leqslant&\int_\Omega |\nabla T_k(u_n)|a(|\nabla T_k(u_n)|)dx\nonumber\\
&\leqslant&(1+a_1)\int_\Omega A(|\nabla T_k(u_n)|)dx\nonumber\\
&\leqslant&C.
\end{eqnarray*}

\n Next it is enough, thanks to Lemma 3.6, to show that there exists a subsequence of $(u_n)$ such that
\begin{eqnarray}\label{e3.37}
\nabla_A T_k(u_n) &\rightarrow & \nabla_A T_k(u) \quad\text{ in
measure in }\Omega.
\end{eqnarray}
\n Let us first establish the following convergence
\begin{eqnarray}\label{e3.38}
\nabla T_k(u_n) &\rightarrow & \nabla T_k(u) \quad\text{ in
measure in }\Omega.
\end{eqnarray}

\n Let $s$ be a positive number and set $F_n=\{|\nabla T_k(u_n)-\nabla T_k(u)|>s\}$. Note that we have
$F_n\subset F_n^1\cup F_n^2\cup F_n^3$, with
\begin{eqnarray*}
&&F_n^1=\{|\nabla u_n-\nabla u|>s\}\cap\{|u_n|\leqslant k\}\cap\{|u|\leqslant k\}\\
&&F_n^2=\{|\nabla u|>s\}\cap\{|u_n|>k\geqslant |u|\}\\
&&F_n^3=\{|\nabla u_n|>s\}\cap\{|u|>k\geqslant |u_n|\}.
\end{eqnarray*}

\n Taking into account that  $(\nabla u_n)$ converges in
measure to $\nabla u$ in $\Omega$, and the fact that $F_n^1\subset \{|\nabla u_n-\nabla u|>s\}$, we see that
$$\lim_{n\rightarrow\infty}{\cal L}^N(F_n^1)=0.$$

\n Since $u, u_n\in W^{1.p}(\Omega)$, we have $\nabla u=0$ a.e. in $\{u=k\}$ and $\nabla u_n=0$ a.e. in $\{u_n=k\}$. Therefore we have
$${\cal L}^N(F_n^2)\leqslant {\cal L}^N(\{|u_n|>k>|u|\})\quad\text{and}\quad {\cal L}^N(F_n^3)\leqslant {\cal L}^N(\{|u|>k>|u_n|\})$$

\n Using the a.e. convergence of $(u_n)$ to $u$ in $\Omega$, we see that
$$\lim_{n\rightarrow\infty}{\cal L}^N(\{|u_n|>k>|u|\})=\lim_{n\rightarrow\infty}{\cal L}^N(\{|u|>k>|u_n|\})=0.$$

\n This leads to $\displaystyle{\lim_{n\rightarrow\infty}{\cal L}^N(F_n^2)=\lim_{n\rightarrow\infty}{\cal L}^N(F_n^3)=0.}$
\n Finally (3.38) follows, since we have ${\cal L}^N(F_n)\leqslant {\cal L}^N(F_n^1)+{\cal L}^N(F_n^2)+{\cal L}^N(F_n^3)$.

\n Let us now establish (3.37). So let $s$ be a positive number and set
$$E_n=\{|\nabla_A T_k(u_n)-\nabla_A T_k(u)|>s\}.$$

\n Note that we have for each positive number $t$ and
each $n\geqslant 1$
$E_n\subset E_n^1\cup E_n^2\cup E_n^3$, with
\begin{eqnarray*}
&&E_n^1=\{|\nabla T_k(u_n)|> t\},~~E_n^2=\{|\nabla T_k(u)|>t\}\\
&&E_n^3=E_n\cap\{|\nabla T_k(u_n)|\leqslant t\}\cap\{|\nabla T_k(u)|\leqslant t\}.\\
\end{eqnarray*}

\n Using the fact that by Lemma 3.4, the sequence $(u_n)$ and the function $u$ are uniformly bounded in $W^{1.p}(\Omega)$, we obtain

\begin{eqnarray*}
{\cal L}^N(E_n^1)\leqslant {1\over{t^p}}\int_\Omega |\nabla T_k(u_n)|^pdx\leqslant{1\over{t^p}}\int_\Omega |\nabla u_n|^pdx\leqslant {C\over{t^p}}.
\end{eqnarray*}
\begin{eqnarray*}
{\cal L}^N(E_n^2)\leqslant {1\over{t^p}}\int_\Omega |\nabla T_k(u)|^pdx\leqslant{1\over{t^p}}\int_\Omega |\nabla u|^pdx\leqslant {C\over{t^p}}.
\end{eqnarray*}

\n We deduce that there exists $t_\epsilon>0$ such that
\begin{equation}\label{3.39}
\forall t\geqslant
t_\epsilon,~~\forall n\geqslant 1,  \quad {\cal L}^N(E_n^1)+{\cal L}^N(E_n^2)<{\epsilon\over 2}.
\end{equation}

\n Now we claim that there exists a positive constant $C(A,N)$ depending only on $A$ and $N$ such that
we have the following inequality (see \cite{[CL4]})
\begin{equation}\label{3.40}
\Big|{{a(|\xi|)}\over{|\xi|}}\xi-{{a(|\zeta|)}\over{|\zeta|}}\zeta\Big|\leqslant C(A,N)|\xi-\zeta|{{a\big((|\xi|^2+|\zeta|^2)^{1/2}\big)}\over{(|\xi|^2+|\zeta|^2)^{1/2}}}~\forall \xi, \zeta\in
\mathbb{R}^N\times\mathbb{R}^N\setminus\{(0,0)\}.
\end{equation}

\n Since $$|\nabla_A T_k(u_n)-\nabla_A T_k(u)|\leqslant|a(\nabla T_k(u_n))|+|a(\nabla T_k(u))|\leqslant
2a\big((|\nabla T_k(u_n)|^2+|\nabla T_k(u)|^2)^{1/2}\big),$$

\n we deduce that we have in $E_n^3$
\begin{eqnarray*}
&&s<2a\big((|\nabla T_k(u_n)|^2+|\nabla T_k(u)|^2)^{1/2}\big)\\
&&(|\nabla T_k(u_n)|^2+|\nabla T_k(u)|^2)^{1/2}\leqslant t\sqrt{2},
\end{eqnarray*}

\n Applying the inequality (3.40), with $\xi=\nabla T_k(u_n)$ and $\zeta=\nabla T_k(u)$, we deduce that we have in $E_n^3$
\begin{eqnarray*}
s<|\nabla_A T_k(u_n)-\nabla_A T_k(u)|&\leqslant& C|\nabla T_k(u_n)-\nabla T_k(u)|
{{a\big((|\nabla T_k(u_n)|^2+|\nabla T_k(u)|^2)^{1/2}\big)}\over{(|\nabla T_k(u_n)|^2+|\nabla T_k(u)|^2)^{1/2}}}\\
&\leqslant&C|\nabla T_k(u_n)-\nabla T_k(u)|{{a(t\sqrt{2})}\over{a^{-1}(s/2)}},
\end{eqnarray*}
\n which leads to
\begin{eqnarray*}
E_n^3\subset \big\{|\nabla T_k(u_n)-\nabla
T_k(u)|>{{sa^{-1}(s/2)}\over{Ca(t\sqrt{2})}}\big\}.
\end{eqnarray*}

\n Using the convergence in measure of $\nabla T_k(u_n)$ to $\nabla T_k(u)$, for $t=t_\epsilon$,
we obtain the existence of $n_0=n_0(s,\epsilon)\geqslant 1$ such that
\begin{equation}\label{3.41}
{\cal L}^N(E_n^3)<{\epsilon\over 2}\quad \forall n,
m\geqslant n_0.
\end{equation}

\n Combining (3.39) and (3.41), we obtain
\begin{equation*}
{\cal L}^N(\{|\nabla_A T_k(u_n)-\nabla_A T_k(u)|>s\})<\epsilon\quad \forall n\geqslant n_0.
\end{equation*}
\n Hence the sequence  $(\nabla_A T_k(u_n))$ converges in measure to $\nabla_A T_k(u)$
and the Lemma follows.

\qed

\vs 0,5 cm

\begin{lemma}\label{l2.9} There exists a subsequence of $u_n$ such that for all $k>0$
\begin{eqnarray*}
&&\nabla_A T_k(u_n)~~ \rightharpoonup ~~ \nabla_A T_k(u)
\quad\text{
weakly in } L^{\widetilde{A}}(\{|u|<k\})\\
&&\lim_{n\rightarrow\infty} \int_\Omega \varphi\nabla_A T_k(u_n).
\nabla (T_k(u_n))dx= \int_\Omega\varphi\nabla_A u .\nabla u
 dx~~\forall\varphi\in {\cal D}(\{|u|<k\}),~\varphi \geqslant 0.
\end{eqnarray*}
\end{lemma}
\vs 0,5cm\n\emph{Proof.} Let us write
\begin{eqnarray}\label{3.42}
&&\int_\Omega \varphi\nabla_A T_k(u_n). \nabla (T_k(u_n))dx=
\int_\Omega\varphi\nabla_A T_k(u_n) .\nabla (T_k(u)) dx\nonumber\\
&&\qquad~+\int_\Omega
\varphi \big(\nabla_A T_k(u_n)-\nabla_A T_k(u)\big)
 . \nabla (T_k(u_n)-T_k(u)) dx\nonumber\\
&&\qquad~+\int_\Omega\varphi\nabla_A T_k(u) .\nabla (T_k(u_n)-T_k(u)) dx\nonumber\\
&&\qquad =I_1^n+I_2^n+I_3^n.
\end{eqnarray}

\n Since $\nabla_A T_k(u_n)$ is bounded
in $L^{\widetilde{A}}(\Omega)$, there exists a subsequence of $u_n$
and an element $\xi\in L^{\widetilde{A}}(\Omega)$ such that
\begin{eqnarray}\label{3.43}
\nabla_A T_k(u_n)~~ \rightharpoonup ~~ \xi \quad\text{ weakly
in } L^{\widetilde{A}}(\Omega).
\end{eqnarray}
\n From Lemma 3.2 and (3.43), we deduce that
\begin{eqnarray}\label{3.44}
\lim_{n\rightarrow\infty} I_1^n=\int_\Omega \varphi\xi.\nabla T_k(u) dx.
\end{eqnarray}

\n For $I_2^n$, we have
\begin{eqnarray}\label{3.45}
I_2^n&=&\int_{|u_n|\leqslant k, |u|\leqslant k}
\varphi \big(\nabla_A u_n-\nabla_A u\big)
 . \nabla (u_n-u) dx +\int_{|u|\leqslant k<|u_n|}\varphi\nabla_A u.\nabla u dx\nonumber\\
&=& I_{2,1}^n+I_{2,2}^n.
\end{eqnarray}

\n Since $\chi_{|u|\leqslant k<|u_n|}$ converges a.e. to 0 in $\Omega$, we obtain
\begin{eqnarray}\label{3.46}
\lim_{n\rightarrow\infty} I_{2,2}^n=0.
\end{eqnarray}

\n Using (2.7), we obtain by Fatou's lemma and
the last estimate in the proof of Lemma 3.5, we get
\begin{eqnarray*}
0\leqslant I_{2,1}^n&\leqslant&\liminf_{m\rightarrow\infty}\int_{|u_n|\leqslant k, |u_m|\leqslant k}
\varphi \big(\nabla_A u_n-\nabla_A u_m\big)
 . \nabla (u_n-u_m) dx \nonumber\\
&\leqslant&||\varphi||_\infty \liminf_{m\rightarrow\infty}\int_{|u_n|\leqslant k, |u_m|\leqslant k}
\big(\nabla_A u_n-\nabla_A u_m\big)
 . \nabla (u_n-u_m) dx \nonumber\\
 &\leqslant&||\varphi||_\infty \liminf_{m\rightarrow\infty}||f_n-f_m||_1=||\varphi||_\infty||f_n-f||_1.
\end{eqnarray*}

\n Letting $n\rightarrow\infty$, we get
\begin{eqnarray}\label{3.47}
\lim_{n\rightarrow\infty} I_{2,1}^n=0.
\end{eqnarray}

\n Using (3.42)-(3.47), we obtain
\begin{eqnarray}\label{3.48}
\lim_{n\rightarrow\infty} \int_\Omega \varphi\nabla_A T_k(u_n). \nabla (T_k(u_n))dx= \int_\Omega \varphi\xi.\nabla T_k(u) dx.
\end{eqnarray}

\n Using (2.7) and setting $\Theta(Y)={{a(|\xi|)}\over{|\xi|}}\xi$, we obtain for every
vector function $Y$
\begin{eqnarray*}
\int_\Omega \varphi\big(\nabla_A T_k(u_n)-\Theta(Y)\big). (\nabla
T_k(u_n)-Y) dx\geqslant 0.
\end{eqnarray*}
Letting  $n\rightarrow\infty$ and taking into account Lemma 3.2 and
(3.48), we get
\begin{eqnarray*}
\int_\Omega \varphi\big(\xi-\Theta(Y)\big). (\nabla T_k(u)-Y) dx\geqslant
0.
\end{eqnarray*}
Now choosing $Y=\nabla T_k(u)-\lambda\vartheta$, where $\lambda$ is
a positive number and $\vartheta$ is an arbitrary continuous vector
function, we obtain
\begin{eqnarray*}
\int_\Omega \varphi\big(\xi-\Theta(\nabla T_k(u)-\lambda\vartheta)\big).
\vartheta dx\geqslant 0.
\end{eqnarray*}
Letting $\lambda\rightarrow 0$, we get
\begin{eqnarray*}
\int_\Omega \varphi\big(\xi-\nabla_A T_k(u)\big). \vartheta
dx\geqslant 0.
\end{eqnarray*}
Since $\vartheta$ is an arbitrary continuous vector function, it
follows that $\xi=\nabla_A T_k(u)$. As a consequence, the Lemma
follows from (3.43) and (3.48). \qed

\section{Proofs of Theorems 1.1-1.6 and Corollaries 1.1-1.2}\label{S3}

\vs 0,5cm\n This section is devoted to the proofs of Theorems
1.1-1.6 and Corollaries 1.1-1.1.

\vs 0,5cm\n\emph{Proof of Theorem 1.1.}

\vs 0,3cm\n \emph{Uniqueness of the Entropy Solution.}
Let $u$ and $v$ be two entropy solutions
of the obstacle problem associated with $(f,\psi,g)$. We argue as in
the proof of Lemma 3.5, replacing $u_n$ and
 $u_m$, respectively by $u$ and $v$,
and $f_n$, $f_m$ by $f$. We get for
$h\geqslant |g|_\infty, |\psi|_\infty$, $t=2k$ and $0<k<h$,

\begin{eqnarray*}&&\int_{\{|u|\leqslant k, |v|\leqslant k\}}\big(\nabla_A u-\nabla_A v\big)
 . \nabla (u-v)dx\nonumber\\
&&\qquad\leqslant  \int_\Omega -f(T_t(u-T_h(v))+T_t(v-T_h(u)))dx\nonumber\\
&&\qquad +\int_{\{h<|u|\leqslant h+t\}} |\nabla u| a(|\nabla u|) dx
+\int_{\{h-t<|u|\leqslant h\}} |\nabla u| a(|\nabla u|) dx\nonumber\\
&&\qquad +\int_{\{h<|v|\leqslant h+t\}} |\nabla v| a(|\nabla v|) dx
+\int_{\{h-t<|v|\leqslant h\}} |\nabla v| a(|\nabla v|) dx.
\end{eqnarray*}
Letting $h\rightarrow\infty$,  and using Lemma 3.3,
we get
$$\int_{\{|u|\leqslant k, |v|\leqslant k\}}\big(\nabla_A u-\nabla_A v\big)
 . \nabla (u-v)dx\leqslant 0.$$
Using (2.7), we obtain $\nabla u=\nabla v$ a.e. in $\{|u|\leqslant k,
|v|\leqslant k\}$. Since this holds for all $k>0$, we have $\nabla
u=\nabla v$ a.e. in $\Omega$. Therefore there exists a constant $c$
such that $u=v+c$ a.e. in $\Omega$. Using Poincar\'{e}'s inequality
for $T_k(u)-T_k(v)=T_k(u)-T_k(g)- (T_k(v)-T_k(g))\in
W^{1,A}_0(\Omega)\hookrightarrow W^{1,1}_0(\Omega)$, we obtain
\begin{eqnarray*}\int_\Omega |T_k(u)-T_k(v)|dx
&\leqslant& C\int_\Omega |\nabla T_k(u)-\nabla T_k(v)|dx
\end{eqnarray*}
from which we deduce that
\begin{eqnarray*}\int_{\{|u|\leqslant k, |v|\leqslant k\}}|u-v|dx
\leqslant C\int_{\{|u|<k<|u|+|c|\}} |\nabla u|dx+C\int_{\{|u|-|c|<k<|u|\}} |\nabla u|dx.
\end{eqnarray*}
Letting $k\rightarrow\infty$, we get
\begin{eqnarray*}\int_\Omega |u-v|dx=0.
\end{eqnarray*}
It follows that $|u-v|=0$ a.e. in $\Omega$, which leads to
$u=v$ a.e. in $\Omega$.

\vs 0,5cm\n \emph{Continuous Dependence.} This is a consequence of Proposition 3.1.

\vs 0,5cm\n\emph{Existence of an entropy solution.} Let $f_n$ be a
sequence of smooth functions converging strongly to $f$ in
$L^1(\Omega)$, with $||f_n||_1\leqslant ||f||_1+1$. We consider the sequence of approximated obstacle
problems associated with $(f_n,\psi,g)$. We know (see [25],
[21]), that there exists a unique variational solution $u_n\in
K_{\psi,g}$  of (1.7), which is also an entropy
solution.

\vs 0.2cm\n Let $v\in K_{\psi,g}$.
Taking $v$ as a test function for (1.9) associated to $(f_n,\psi,g)$,
we get
$$\int_\Omega \nabla_A u_n
 . \nabla (T_t(v-u_n))dx \geqslant \int_\Omega  -f_nT_t(v-u_n)dx.$$
Since $\{|v-u_n|<t\}\subset \{|u_n|<s\}$, with
$s=t+|v|_\infty$, the previous inequality can be
written as
\begin{equation}\label{4.1}
\int_\Omega \chi_n\nabla_A T_s(u_n)
 . \nabla v dx \geqslant \int_\Omega  -f_nT_t(v-u_n)dx+\int_\Omega \chi_n\nabla_A T_s(u_n)
 . \nabla (T_s(u_n))dx,
\end{equation}

\n where $\chi_n=\chi_{\{|v-u_n|<t\}}$. It is clear that
$\chi_n\rightharpoonup \chi$ weakly$*$ in $L^\infty(\Omega)$.
Moreover $\chi_n$ converges a.e. to $\chi_{\{|v-u|<t\}}$ in $\Omega\setminus\{|v-u|=t\}$.
It follows that
\begin{equation*}\chi=
\begin{cases}
& 1\quad
\text{ in }\quad \{|v-u|<t\}\\
 & 0\quad\text{ in }\quad \{|v-u|>t\}.
\end{cases}
\end{equation*}

\n Note that we have ${\cal L}^N(\{|v-u|=t\})=0$ for a.e. $t\in(0,\infty)$. So there exists
a measurable set $\mathcal{N}\subset (0,\infty)$ such that ${\cal L}^N(\{|v-u|=t\})=0$ for all
$t\in (0,\infty)\setminus \mathcal{N}$.

\n Assume that $t\in (0,\infty)\setminus \mathcal{N}$. Then
$\chi_n$ converges weakly$*$ in $L^\infty(\Omega)$ and a.e. in $\Omega$ to
$\chi=\chi_{\{|v-u|<t\}}$. Since $\nabla(T_s(u_n))$
converges a.e. to $\nabla(T_s(u))$ in $\Omega$ (Proposition 3.1), we obtain by Fatou's Lemma

 \begin{equation}\label{4.2}
\liminf_{n\rightarrow\infty}\int_\Omega \chi_n\nabla_A T_s(u_n)
 . \nabla T_s(u_n) dx\geqslant\int_\Omega \chi\nabla_A T_s(u)
 . \nabla T_s(u )dx.
\end{equation}
Using the strong convergence of $\nabla_A T_s(u_n)$ to $\nabla_A T_s(u)$ in
$L^1(\Omega)$ and the weak$*$ convergence of $\chi_n$ to $\chi$ in $L^\infty(\Omega)$, we obtain
 \begin{equation}\label{4.3}
\lim_{n\rightarrow\infty}\int_\Omega \chi_n\nabla_A T_s(u_n)
 . \nabla v dx=\int_\Omega \chi\nabla_A T_s(u)
 . \nabla v dx.
\end{equation}
Moreover since $f_n$ converges to $f$ in $L^1(\Omega)$ and
$T_t(v-u_n)$ converges to $T_t(v-u)$ weakly$*$ in
$L^\infty(\Omega)$, we obtain by passing to the limit in (4.1) and
taking into account (4.2)-(4.3)
\begin{equation*}
\int_\Omega \chi\nabla_A T_s(u)
 . \nabla v dx-\int_\Omega \chi\nabla_A T_s(u)
 . \nabla (T_s(u))dx \geqslant \int_\Omega  -fT_t(v-u)dx,
\end{equation*}
which can be written as
\begin{equation*}
\int_{\{|v-u|\leqslant t\}} \chi\nabla_A T_s(u)
 . (\nabla v -\nabla u)dx \geqslant \int_\Omega  -fT_t(v-u)dx,
\end{equation*}
or since $\chi=\chi_{\{|v-u|<t\}}$ and
$\nabla(T_t(v-u))=\chi_{\{|v-u|<t\}}\nabla(v-u)$
\begin{equation*}
\int_\Omega\nabla_A u
 . \nabla(T_t(v-u))dx \geqslant \int_\Omega  fT_t(v-u)dx\quad \forall t\in (0,\infty)\setminus \mathcal{N}.
\end{equation*}

\n For $t\in \mathcal{N}$, we know since $|\mathcal{N}|=0$, that
there exists a sequence $(t_k)$ of numbers in $(0,\infty)\setminus
\mathcal{N}$ such that $t_k \longrightarrow t$. Therefore we have
\begin{equation}\label{4.4}
\int_\Omega\nabla_A u
 . \nabla(T_{t_k}(v-u))dx \geqslant \int_\Omega  -fT_{t_k}(v-u)dx\quad \forall k \geqslant1.
\end{equation}

\n Since $\nabla(v-u)=0$ a.e. in $\{|v-u|=t\}$, the left hand side of (4.4) can be written as

\begin{equation*}
\int_\Omega\nabla_A u
 . \nabla(T_{t_k}(v-u))dx=\int_{\Omega\setminus|v-u|=t\}} \chi_{\{|v-u|<t_k\}}\nabla_A u. \nabla(v-u)dx.
\end{equation*}
The sequence $\chi_{\{|v-u|<t_k\}}$ converges to $\chi_{\{|v-u|<t\}}$ a.e.  in $\Omega\setminus\{|v-u|=t\}$ and therefore converges weakly$*$ in $L^\infty(\Omega\setminus\{|v-u|=t\})$. We obtain
\begin{eqnarray}\label{4.5}
\lim_{k\rightarrow\infty}\int_\Omega\nabla_A u
 . \nabla(T_{t_k}(v-u))dx&=&\int_{\Omega\setminus\{|v-u|=t\}} \chi_{\{|v-u|<t\}}\nabla_A u. \nabla(v-u)dx\nonumber\\
 &=&\int_{\Omega} \chi_{\{|v-u|<t\}}\nabla_A u. \nabla(v-u)dx\nonumber\\
 &=&\int_\Omega\nabla_A u
 . \nabla(T_{t}(v-u))dx.
\end{eqnarray}

\n For the right hand side of (4.4), we have
\begin{eqnarray}\label{4.6}
\Big|\int_\Omega  -fT_{t_k}(v-u)dx-\int_\Omega  -fT_{t}(v-u)dx\Big|\leqslant|t_k-t|.|f|_1\rightarrow~0~\text{as}~k\rightarrow\infty.
\end{eqnarray}

\n It follows from (4.4)-(4.6) that we have the inequality
\begin{equation*}
\int_\Omega\nabla_A u
 . \nabla(T_t(v-u))dx \geqslant \int_\Omega  -fT_t(v-u)dx\quad \forall t\in (0,\infty).
\end{equation*}
Hence $u$ is an entropy solution of the obstacle problem associated
with $(f,\psi,g)$.
\qed

\vs 0,5cm\n\emph{Proof of Theorem 1.2.} Let $f_n$ be
a bounded sequence of smooth functions such that

\begin{eqnarray}\label{4.7}
f_n\,\rightarrow\,f\quad\text{in} \quad L^1(\Omega).
\end{eqnarray}
We consider the unique variational solution $u_n\in K_{\psi,g}$ of
the obstacle problem associated with $(f_n,\psi,g)$.
From Proposition 2.1, we have
\begin{equation}\label{4.8}
f_n-(f_n-\Delta_A\psi)^+\leqslant \Delta_A u_n\leqslant f_n
\quad\text{a.e. in}\quad\Omega.
\end{equation}
Let us first give a proof assuming that $a_1<{N\over{N-1}}a_0$. We
deduce from (4.8) that we have for every $\varphi\in {\cal
D}(\Omega)$, $\varphi\geqslant 0$
\begin{eqnarray*}
 \int_\Omega -f_n\varphi dx\leqslant \int_\Omega \nabla_A u_n
 . \nabla \varphi dx \leqslant
 \int_\Omega (-f_n+(f_n-\Delta_A\psi)^+)\varphi dx.
\end{eqnarray*}
Using (4.7) and the strong convergence of $\nabla_A u_n$ to $\nabla_A u$ in $L^1(\Omega)$ (since
$a_1<{N\over{N-1}}a_0$, see Proposition 3.1), we get by letting
$n\rightarrow \infty$
\begin{eqnarray*}
 \int_\Omega -f\varphi dx\leqslant \int_\Omega \nabla_A u
 . \nabla \varphi dx \leqslant
 \int_\Omega (-f+(f-\Delta_A\psi)^+)\varphi dx.
\end{eqnarray*}
Hence we obtain
$$f-(f-\Delta_A\psi)^+\leqslant \Delta_A u\leqslant f\quad\text{a.e.
in}\quad\Omega.$$

\vs 0,3cm\n Let us now deal with the general case.
Choose $\epsilon>0$, $k>||\psi||_\infty$, and consider $H_\epsilon(t)$ the function
introduced in the proof of Proposition 2.1. Let $\varphi\in {\cal
D}(\{u<k\})$ with $\varphi\geqslant 0$. Denoting the functions
$-f_n+(f_n-\Delta_A\psi)^+$ and $-f+(f-\Delta_A\psi)^+$
respectively by $h_n$ and $h$, and multiplying (4.8) by $\varphi
H_\epsilon(k+\epsilon-u_n)$ and integrating, we obtain

\begin{eqnarray}\label{4.9}
I_{n,k}+\int_\Omega -f_n\varphi H_\epsilon(k+\epsilon-u_n) dx\leqslant
\int_\Omega H_\epsilon(k+\epsilon-u_n)\nabla_A (T_{k+\epsilon}(u_n))
 . \nabla\varphi dx
\end{eqnarray}

\begin{eqnarray}\label{4.10}
\int_\Omega H_\epsilon(k+\epsilon-u_n)\nabla_A (T_{k+\epsilon}(u_n))
 . \nabla\varphi dx \leqslant I_{n,k}+\int_\Omega h_n\varphi H_\epsilon(k+\epsilon-u_n)dx.
\end{eqnarray}
where
\begin{eqnarray*}
I_{n,k}&=& {1\over\epsilon}\int_{\{k<u_n<k+\epsilon\}}\varphi\nabla_A u_n . \nabla u_n dx.
\end{eqnarray*}

\n Since $u_n\geq\psi>-k$, we have

\begin{eqnarray*}
I_{n,k}&=& {1\over\epsilon}\int_{\{u_n<k+\epsilon\}} \varphi\nabla_A u_n.\nabla u_n dx-{1\over\epsilon}\int_{\{u_n<k\}} \varphi\nabla_A u_n.\nabla u_n dx\nonumber\\
&=&{1\over\epsilon}\int_{\{|u_n|<k+\epsilon\}} \varphi\nabla_A u_n.\nabla u_n dx-{1\over\epsilon}\int_{\{|u_n|<k\}} \varphi\nabla_A u_n.\nabla u_n dx\nonumber\\
&=&{1\over\epsilon}\int_\Omega \varphi\nabla_A (T_{k+\epsilon}(u_n)).\nabla
T_{k+\epsilon}(u_n)dx-{1\over\epsilon}\int_\Omega \varphi\nabla_A (T_k(u_n)).\nabla T_k(u_n)dx.
\end{eqnarray*}
Letting $n\rightarrow\infty$, we obtain by using Lemma 3.8 and the
fact that $supp(\varphi)\subset \{|u|<k\}$
\begin{eqnarray*}
\lim_{n\rightarrow\infty}I_{n,k}={1\over\epsilon}\int_\Omega \varphi\nabla_A(T_{k+\epsilon}(u)).\nabla T_{k+\epsilon}(u)dx
-{1\over\epsilon}\int_\Omega \varphi\nabla_A T_k(u).\nabla T_k(u)dx=0.
\end{eqnarray*}

\n Now letting $n\rightarrow\infty$ in (4.9)-(4.10), we obtain by
using the results of Lemmas 3.2 and 3.8
\begin{eqnarray*}
\int_\Omega f\varphi H_\epsilon(k+\epsilon-u) dx\leqslant\int_\Omega
H_\epsilon(k+\epsilon-u)\nabla_A(T_{k+\epsilon}(u))
 . \nabla\varphi dx \leqslant \int_\Omega h\varphi H_\epsilon(k+\epsilon-u)dx
\end{eqnarray*}
which can be written since $supp(\varphi)\subset \{u<k\}$ and $k+\epsilon-u>\epsilon$ in
$\{u<k\}$
\begin{eqnarray*}
\int_\Omega -f\varphi dx\leqslant\int_\Omega \nabla_A u.\nabla\varphi dx\leqslant \int_\Omega (-f+(f-\Delta_A\psi)^+\varphi dx.
\end{eqnarray*}
This means that
$$f-(f-\Delta_A\psi)^+\leqslant \Delta_A u\leqslant f\quad\text{a.e.
in}\quad\{|u|<k\},\quad\forall k>||\psi||_\infty.$$

\n Hence we have proved that
$$f-(f-\Delta_A\psi)^+\leqslant \Delta_A u\leqslant f\quad\text{a.e.
in}\quad\Omega.$$ \qed

\vs 0,5cm\n\emph{Proof of Theorem 1.3.} Let $u_1, u_2$ be two entropy
solutions of the obstacle problem associated with $(f_1,\psi, g)$
and $(f_2,\psi, g)$. Let $f^1_n, f^2_n$ be
two bounded sequences of smooth functions such that

\begin{eqnarray}\label{4.11}
f^i_n\,\rightarrow\,f^i\quad\text{in} \quad L^1(\Omega).
\end{eqnarray}

\n We consider the unique variational solutions $u^i_n\in K_{\psi,g}$ of
the obstacle problem associated with $(f^i_n,\psi,g)$.
Let $\xi_n^i=\Delta_A u_n^i-f_n^i$, $i=1, 2$.
Then we have by Proposition 2.3
\begin{equation}\label{4.12}
\int_\Omega|\xi_n^1-\xi_n^2|dx\leqslant \int_\Omega|f_n^1-f_n^2|dx.
\end{equation}

\n Let $A_k=\{|u_1|<k\}\cap \{|u_2|<k\}$ and $A_k^n=\{|u_n^1|<k\}\cap \{|u_n^2|<k\}$. We deduce from (4.12) that
\begin{equation*}
\int_{A_k^n\cap A_k}|\xi_n^1-\xi_n^2|dx\leqslant \int_\Omega|f_n^1-f_n^2|dx,
\end{equation*}

\n which can be written as
\begin{equation}\label{4.13}
\int_{A_k}|\chi_{A_k^n}f_n^2-\chi_{A_k^n}f_n^1+\chi_{A_k^n}(\Delta_A u_n^2-\Delta_A u_n^1)|dx\leqslant \int_\Omega|f_n^1-f_n^2|dx.
\end{equation}

\n Since $\chi_{A_k^n}\rightarrow 1$ a.e. in $A_k$, we have up to a subsequence
\begin{equation}\label{4.14}
\chi_{A_k^n}~~\rightarrow~~ 1\quad\text{weakly-* in} \quad L^1(A_k).
\end{equation}
\n Using (4.11) and (4.14), we obtain
\begin{equation}\label{4.15}
\chi_{A_k^n}f_n^i~~\rightarrow~~ f_i\quad\text{in} \quad L^1(A_k).
\end{equation}

\n Now we have by Proposition 2.1, for any measurable subset $A$ of
$\Omega$
\begin{eqnarray*}
\int_A |\chi_{A_k^n}\Delta_A u_n^i|dx&\leqslant &||f^i_n||_{L^1(A)}+||\Delta_A \psi||_{L^1(A)}\\
&\leqslant &||f_n^i-f_i||_{L^1(A)}+||f_i||_{L^1(A)}+||\Delta_A \psi||_{L^1(A)}.
\end{eqnarray*}

\n It follows that the sequences $(\chi_{A_k^n}\Delta_A u_n^i)$ are bounded in $L^1(\Omega)$, and
satisfy the assumptions of Dunford-Pettis Theorem. Therefore we have for a subsequence
\begin{equation}\label{4.16}
\chi_{A_k^n}\Delta_A u_n^i~~\rightharpoonup~~ {\cal U}_i\quad\text{in} \quad L^1(\Omega).
\end{equation}

\n We shall prove that ${\cal U}_i=\Delta_A u_i$ a.e. in $A_k$. To do that, we will prove that
$\chi_{A_k^n}\Delta_A u_n^i$ converges to $\Delta_A u_i$  in ${\cal D}'(A_k)$.
Indeed we have for $\varphi\in {\cal D}(A_k)$
\begin{eqnarray}\label{4.17}
\int_{A_k} \chi_{A_k^n}\Delta_A u_n^i\varphi dx&=&\int_{A_k\cap A_k^n}\Delta_A u_n^i\varphi dx\nonumber\\
&=&-\int_{A_k\cap A_k^n} \nabla_A u_n^i.\nabla\varphi dx\nonumber\\
&= &-\int_{A_k} \chi_{A_k^n}\nabla_A (T_k(u_n^i)).\nabla\varphi dx.
\end{eqnarray}

\n Using (4.14) and Lemma 3.7, and passing to the limit in (4.17), we get
\begin{eqnarray}\label{4.18}
\lim_{n\rightarrow\infty}\int_{A_k} \chi_{A_k^n}\Delta_A u_n^i\varphi dx&=&-\int_{A_k}\nabla_A (T_k(u_i)).\nabla\varphi dx=\int_{A_k} \Delta_A u_i\varphi dx.
\end{eqnarray}
\n We deduce from (4.16) and (4.18) that we have ${\cal U}_i=\Delta_A u_i$ a.e. in $A_k$.

\vs0.2cm\n It follows from (4.14) and (4.18) that $\chi_{A_k^n}f_n^2-\chi_{A_k^n}f_n^1+\chi_{A_k^n}(\Delta_A (T_k(u_n^2))-\Delta_A (T_k(u_n^1)))$ converges
weakly to $f_2-f_1+(\Delta_A (T_k(u_2))-\Delta_A (T_k(u_1)))$ in $L^1(A_k)$.
Using the semicontinuity of the norm $||.||_{L^1(A_k)}$, (4.11) and the fact that
$f_2-f_1+(\Delta_A (T_k(u_2))-\Delta_A (T_k(u_1)))=f_2-f_1+\Delta_A u_2-\Delta_A u_1=\xi_1-\xi_2$ in
$A_k$, we get from (4.13)
\begin{eqnarray*}
\int_{A_k}|\xi_1-\xi_2|dx&\leqslant &\liminf_{n\rightarrow\infty}\int_{A_k}|\chi_{A_k^n}f_n^2-\chi_{A_k^n}f_n^1+(\Delta_A (T_k(u_n^2))-\Delta_A (T_k(u_n^1)))|dx\\
&\leqslant& \liminf_{n\rightarrow\infty}\int_\Omega|f_n^1-f_n^2|dx\\
&=&\int_\Omega|f_1-f_2|dx.
\end{eqnarray*}

\n Since $k$ is arbitrary, we get
\begin{eqnarray*}
\int_\Omega|\xi_1-\xi_2|dx\leqslant \int_\Omega|f_1-f_2|dx.
\end{eqnarray*}
\qed

\vs 0,5cm\n\emph{Proof of Corollary 1.1.} Let $f_n$  be
a sequence of bounded smooth functions such that
\begin{eqnarray}\label{4.16}
&&f_n\,\rightarrow\,f\quad\text{in} \quad L^1(\Omega).
\end{eqnarray}
From Propositions 2.1 and 2.2, we have for a function $0\leqslant q_n\leqslant\chi_n=\chi_{\{u_n=\psi\}}$
\begin{equation}\label{4.17}
\Delta_A u_n-(\Delta_A\psi-f_n)q_n=f_n\quad\text{ a.e. in
}\Omega.
\end{equation}
\n Since $0\leqslant\chi_n\leqslant
1$, there exists a subsequence of $\chi_n$ and a function $\chi\in
L^\infty(\Omega)$ such that
\begin{equation}\label{4.18}
\chi_n\rightharpoonup \chi_*\quad\text{weakly* in }L^\infty(\Omega).
\end{equation}
\n Similarly, there exists a subsequence of $q_n$ and a function $q\in
L^\infty(\Omega)$ such that
\begin{equation}\label{4.19}
q_n\rightharpoonup q\quad\text{weakly* in }L^\infty(\Omega).
\end{equation}

\n Using (4.16) and (4.19) and arguing as in the proof of
Theorem 1.2, we obtain by passing to the limit in (4.17) that
\begin{equation}\label{4.20}
\Delta_A u-(\Delta_A\psi-f)q=f\quad\text{ a.e. in }\Omega.
\end{equation}
We are going to prove that $\chi=\chi_{\{u=\psi\}}$.
Since $q_n\leqslant\chi_n$, we get $q\leqslant\chi_*$ a.e. in $\Omega$.
Passing to the limit in $\displaystyle{\int_\Omega (u_n-\psi)\chi_n dx}=0$, we get
$\displaystyle{\int_\Omega (u-\psi)\chi_* dx}=0$. Since $0\leqslant \chi_*\leqslant 1$,
we obtain $\chi_*=0$ a.e. in $\{u>\psi\}$. Hence $q\leqslant\chi_*\leqslant \chi$
a.e. in $\Omega$.

\n Since $u$ satisfies (1.13), we get $(\Delta_A\psi-f)^+q=(\Delta_A\psi-f)^+\chi$,
a.e. in $\Omega$, which can be written as $(\Delta_A\psi-f)(\chi-q)$ a.e. in $\Omega$.
Because $\Delta_A\psi\neq f$ a.e. in  $\Omega$, we deduce that we have
$q=\chi=\chi_{\{u=\psi\}}$ a.e. in $\Omega$. Hence (4.18) becomes
\begin{equation*}
\chi_{\{u_n=\psi\}}\rightharpoonup \chi_{\{u=\psi\}}\quad\text{weakly*
in }L^\infty(\Omega).
\end{equation*}

\n Since $\chi_n$ and $\chi$ are characteristic functions, we obtain
\begin{equation*}
\chi_{\{u_n=\psi\}}\rightarrow\chi_{\{u=\psi\}}\quad\text{in }L^s(\Omega)~~\forall s<\infty.
\end{equation*}
\qed

\vs 0,5cm\n\emph{Proof of Corollary 1.2.} Let $u_1, u_2$ be two entropy
solutions of the obstacle problem associated with $(f_1,\psi, g)$
and $(f_2,\psi, g)$ and satisfying (1.13). Let $\omega$ be a measurable
subset of $\Omega$ such that we have, for some $\lambda>0$,
$\Delta_A\psi-f_i\leqslant -\lambda$, a.e. in
$\omega$, $i=1, 2$. Then we have by (1.13)
\begin{eqnarray*}
{\cal L}^N\big((I_1\div I_2)\cap \omega\big)&=&\int_\omega \chi_{I_1\div I_2}dx=\int_\omega |\chi_{I_1}-\chi_{I_2}|dx\\
&=&\int_\omega |\chi_{\{u=\psi_1\}}-\chi_{\{u=\psi_2\}}|dx\\
&\leqslant&{1\over\lambda}\int_\omega |(\Delta_A\psi-f_1)\chi_{\{u=\psi_1\}}-(\Delta_A\psi-f_2)\chi_{\{u=\psi_2\}}|dx\\
&=&{1\over\lambda}\int_\omega |(\Delta_Au_1-f_1)-(\Delta_Au_2-f_2)|dx\\
&=&{1\over\lambda}\int_\omega |\xi_1-\xi_2|dx\\
&\leqslant&{1\over\lambda}\int_\Omega|f_1-f_2|dx.
\end{eqnarray*}
\qed

\vs 0.5 cm
\section{Growth of the Gradient near the Free Boundary}\label{S4}

\vs 0,5cm In all what follows, we assume that the obstacle
$\psi\equiv 0$ and without loss of generality, we assume that
$x_0=0$ and $\omega_0=B_1=\{x~:~||x||<1~\}$.

\n Moreover due to the local character of the results of the
coming sections, we will restrict ourselves to the unit ball and
will consider the solutions of the following class of problems
$$ {\cal F}_A : \begin{cases}
    & u\in W^{1,A}(B_1)\cap C^{1, \alpha}_{loc}(B_1) \\
    & \Delta_A u=f\quad \text{in }\quad \{u>0\}\cap B_1,\\
    & 0\leqslant u \leqslant M_0 \quad \text{in }\quad B_1,\\
    & 0\in \partial \{u>0\},
  \end{cases}$$
where $M_0$ is a positive number. We may also assume that
there exists a positive constant
$M_1=M_1(a_0,a_1,N,a(1), M_0,\Lambda_0)$ such that

\begin{equation}\label{5.1}
|u|_{1,\alpha,\overline{B}_{3/4}} \leqslant M_1 \qquad \forall u\in
\mathcal{F}_A.
  \end{equation}

\vs 0.5 cm \n The following theorems give the growth of the elements
of the family $ \mathcal{F}_A$ and their gradients near the free boundary.
The first one was proved in \cite{[CL1]} by contradiction. We provide here a direct proof.

\vs 0,5cm
\begin{proposition}\label{p5.1} There exists a positive  constant
$C_0=C_0(a_0,a_1,N,M_0,\Lambda_0)$ such that for every $u\in
\mathcal{F}_A$, we have
$$ 0\leqslant u(x) \leqslant C_0 \widetilde{A} (|x|) \qquad \forall x\in B_1$$
\n where  $\widetilde{A}$  is the function defined by (1.21).
\end{proposition}

\vs 0,5cm
 \n \emph{Proof.} Let $k\in \mathbb{N}\cup\{0\}$ and set $\omega_k(x)= \displaystyle{ u( 2^{-k} x) \over
 2^{-k}}$ for $x\in B_1$.  We have

\begin{align*}
   & \Delta_A \omega_k(x) =  2^{-k}(\Delta_A u)(2^{-k}x),
   \qquad \| \Delta_A \omega_k\|_{\infty} \leqslant 2^{-k}\Lambda_0   \\
   & 0\leqslant \omega_k\leqslant\displaystyle {M_0\over
   2^{-k}}\quad \hbox{ in } B_1, \qquad \displaystyle {\inf_{B_{1/2}} } \omega_k
   =0.
\end{align*}

\n By Harnack's inequality (see Corollary 1.4 \cite{[L1]}), we have
for some constant $C=C(a_0,a_1,N,M_0)$,
$$\sup_{B_{1/2}} \omega_k \leqslant C\big(  \inf_{B_{1/2}} \omega_k +
a^{-1} ( 2^{-k}\Lambda_0 )  \big) = C a^{-1} ( 2^{-k}\Lambda_0).$$

\n We deduce from (2.5)-(2.6) that
$$\sup_{B_{2^{-k-1}}} u=2^{-k}\sup_{B_{1/2}}\omega_k \leqslant
2C 2^{-k-1}a^{-1} ( 2^{-k-1}2\Lambda_0)\leqslant
C'\widetilde{A}(2^{-k-1} )$$ where $C'=2C\big(1+{1\over
a_0}\big)\max\big((2\Lambda_0)^{1\over a_0},(2\Lambda_0)^{1\over
a_1}\big)$.

\n Now, let $x\in
 B_{1/2}$ and set $r=|x|$. Then there exists $k\in \mathbb{N}\cup\{0\}$ such
 that  $2^{-k-2}\leqslant r<2^{-k-1}$ and we have by (2.5)-(2.6)
 $$ u(x) \leqslant \sup_{B_{2^{-k-1}}} u \leqslant C' \widetilde{A} (2^{-k-1})
  \leqslant C'\widetilde{A} (2r)
 \leqslant C'. 2^{(1+{1\over
a_0})} \widetilde{A}(r)=2^{(1+{1\over
a_0})} C'\widetilde{A}(|x|).$$

\n If $x\in B_1\setminus B_{1/2}$, we have
$$ u(x) \leqslant M_0={M_0\over \widetilde{A}(|x|)}\widetilde{A}(|x|)
\leqslant {M_0\over \widetilde{A}(2^{-1})}\widetilde{A}(|x|).$$
 \qed

\vs 0,5cm
\begin{proposition}\label{p5.2} There exists a positive  constant
$C_1=C_1(a_0,a_1,a(1), N, M_0,M_1,\Lambda_0)$ such that for every $u\in \mathcal{F}_A$,
we have
\begin{eqnarray}\label{5.12}
&& |\nabla u(x)| \leqslant C_1 a^{-1}(|x|) \qquad \forall x\in B_{3/4}.
\end{eqnarray}
\end{proposition}

\vs 0,5cm In order to prove Proposition 5.2, we need to introduce some
notations. For a nonnegative  bounded function $u$, we define the
quantity

$$ S(r,u) = \sup_{x\in B_r(0)} u(x).$$

\n We also define for $u\in\mathcal{F}_A$ the set
$$ \mathbb{P}(u) = \{ j\in \mathbb{N}/\quad 2^{1\over{a_0}} S( 2^{-j-1},|\nabla u|)
\geqslant S( 2^{-j},|\nabla u|) \}.$$

\n Then we have

\begin{lemma}\label{l5.1} Assume that $ \mathbb{P}(u)
\neq\emptyset$. Then, there exists  a constant

\n $c_1=c_1(a_0,a_1,a(1), N, M_0, M_1,\Lambda_0)$ such that
$$ S( 2^{-j-1},|\nabla u|) \leqslant c_1 a^{-1}( 2^{-j}) \qquad
 \forall u\in \mathcal{F}_A , \quad \forall j\in \mathbb{P}(u).$$
\end{lemma}

\vs 0,5cm
 \n \emph{Proof.} We argue by contradiction.   So we assume that
  $$ \forall k\in  \mathbb{N},\quad \exists u_k \in \mathcal{F}_A,\quad
  \exists j_k\in \mathbb{P}(u_k) \quad \hbox{ such that }\quad
  S( 2^{-j_k-1},|\nabla u_k|) \geqslant k a^{-1}( 2^{-j_k}).$$
 \n
 Consider  $w_k(x)= \displaystyle{ u_k( 2^{-j_k} x) \over
 2^{-j_k}}$ for $x\in B_1$.  We have $\nabla w_k(x)=\nabla u_k( 2^{-j_k} x)$
 and by Proposition 5.1
\begin{align*}
   & \Delta_A w_k(x) =  2^{-j_k}(\Delta_A u_k )(2^{-j_k}x),
   \qquad \| \Delta_A w_k\|_{\infty} \leqslant 2^{-j_k}\Lambda_0   \\
   & 0\leqslant w_k\leqslant\displaystyle C_0{{\widetilde{A}(2^{-j_k})}\over
   {2^{-j_k}}}\quad \hbox{ in } B_1,
   \qquad
   w_k(0)=0.
\end{align*}

\n Now let  $v_k(x)= \displaystyle{  {w_k(x)} \over S( 2^{-j_k-1},
|\nabla u_k|) }  =   \displaystyle{ u_k(2^{-j_k}x) \over {2^{-j_k}S(
2^{-j_k-1}, |\nabla u_k|)} }=\displaystyle{ u_k(2^{-j_k}x) \over m_k
} $ for $x\in B_1$.

\n We introduce the functions
$$b_k(t)=\displaystyle{a(t m_k)\over a(m_k)},\qquad
B_k(t)= \int_{0}^{t} b_k(\tau) d\tau= {A(m_kt)\over {m_k a(m_k)}}.$$

\n Then it is easy to see that $v_k$ satisfies

\begin{eqnarray*}
&&\|\Delta_{B_k} v_k\|_\infty = \|\Delta_A w_k\|_\infty
\leqslant 2^{-j_k} \Lambda_0\\
&& 0 \leqslant v_k\leqslant
 {{C_0\widetilde{A} (2^{-j_k})}\over {2^{-j_k}S( 2^{-j_k-1}, u_k) }}
 \leqslant{{C_0\widetilde{A} (2^{-j_k})}\over {2^{-j_k}k a^{-1}(2^{-j_k})}}\leqslant{C_0\over k}
 \quad \hbox{ in } B_1 \quad \hbox{ by } (2.6)\\
 && \sup_{B_{1/2}} |\nabla v_k|=1\\
&&  v_k(0)=0.
\end{eqnarray*}

\n Since  $b_k$ satisfies $(1.5)$ with the same constants $a_0$ and
$a_1$, we have $v_k\in C^{1,\alpha}_{loc}(B_1)$ with
$$|v_k|_{1,\alpha,B_{3/4}} \leqslant C(a_0,a_1,a(1), N,C_0,\Lambda_0).$$

\n Hence up to a subsequence, we have $v \rightarrow v_k$  in
$C^{1,\beta}(\overline{B_{3/4}})$, for all $\beta\in(0,\alpha)$. We
deduce that $\sup_{B_{1/2}} |\nabla v_k|\rightarrow \sup_{B_{1/2}}
|\nabla v|$. Therefore we have $\sup_{B_{1/2}} |\nabla v|=1$. But
this contradicts the fact that $$\displaystyle{0\leqslant
v(x)=lim_{k\rightarrow\infty}v_k(x)\leqslant
lim_{k\rightarrow\infty}{C\over k}=0\quad \hbox{ for all } x\in B_1.}$$

\qed

\vspace {0,5cm}
 \n \emph{Proof \, of \, Proposition \,5.2.}

\n Let $S(3/4,|\nabla u|)=M_1=C(a_0,a_1,a(1),N,M_0,\Lambda_0)$. First,
we prove by induction that we have for  $c_2=\max
({M_1\over{a^{-1}(1/2)}}, c_1 2^{1\over{a_0}})$
 $$ S(2^{-j},|\nabla u|) \leqslant c_2 a^{-1}(2^{-j})\quad \forall j\in
 \mathbb{N},\quad\forall u\in {\cal F}_A.$$

\n Let $u\in {\cal F}_A.$ For $j=1$, we have
\begin{eqnarray*}
S(2^{-1},|\nabla u|)= S(1/2,|\nabla u|)\leqslant
M_1={M_1\over{a^{-1}(1/2)}}a^{-1}(1/2)\leqslant c_2 a^{-1}(2^{-1}).
\end{eqnarray*}

\n Let $j\geqslant 1$ and assume that $ S(2^{-j},u) \leqslant c_2
a^{-1}(2^{-j})$.  We distinguish two cases :

 \par -- If $j\in \mathbb{P}(u)$, we have by Lemma 5.1 and (2.5),
 $$S(2^{-j-1},|\nabla u|) \leqslant c_1
a^{-1}(2^{-j})=c_1 a^{-1}(2.2^{-j-1})\leqslant c_1 2^{1\over{a_0}}
a^{-1}(2^{-j-1})\leqslant c_2 a^{-1}(2^{-j-1}).$$

 \par -- If $j\notin \mathbb{P}(u)$, we have
$S(2^{-j-1},|\nabla u|) < {{S(2^{-j},u)}\over{2^{1\over{a_0}}}}$.
Using the induction assumption and (2.5), we get
 $$ S(2^{-j-1},|\nabla u|) \leqslant {{c_2}\over {2^{1\over{a_0}}}}
a^{-1}(2^{-j}) \leqslant {{c_2}\over {2^{1\over{a_0}}}}
2^{1\over{a_0}}a^{-1}(2^{-j-1})=c_2 a^{-1}(2^{-j-1}).$$

\n Now, let $x\in B_{1/2}$ and set $r=|x|$. Then there exists $j\in \mathbb{N}$
such that  $2^{-j-1}\leqslant r\leqslant 2^{-j}$ and we have for
$C_1=\max({M_1\over{a^{-1}(1/2)}}, c_2 2^{1\over{a_0}})$
\begin{eqnarray*}
|\nabla u(x)| &\leqslant& \sup_{y\in\overline{B_{2^{-j}}}} |\nabla
u(y)|= S(2^{-j},|\nabla u|)\leqslant c_2 a^{-1}(2^{-j})\\
&\leqslant & c_2 a^{-1}(2r)
 \leqslant c_2 2^{1\over{a_0}}  a^{-1}(r)\leqslant C_1a^{-1}(|x|).
 \end{eqnarray*}

\n If $x\in B_{3/4}\setminus B_{1/2}$, we have
$$ |\nabla u(x)| \leqslant M_1={M_1\over {a^{-1}(|x|)}}a^{-1}(|x|)
\leqslant {M_1\over a^{-1}(2^{-1})}a^{-1}(|x|)\leqslant
C_1a^{-1}(|x|).$$

\qed

\vs 0,5 cm

\section{Hausdorff Measure Estimate of the Free Boundary}\label{S6}

\vs 0.5 cm\n We shall start by establishing local
$L^2-$estimate for the second derivatives of $u$. In order to do that,
we define for each $r>0$ and each function $u\in \mathcal{F}_A$, the
quantity
\begin{eqnarray*}E(r,u) &=& {1\over{|B_r|}} \int_{B_r\cap\{\nabla u(x)\neq 0\}}
\Big[{{a(|\nabla u|}\over{|\nabla u|}} |D^2 u|\Big]^2dx\\
&=& {1\over{|B_1|}} \int_{B_1\cap\{\nabla u(rx)\neq 0\}}
\Big[{{a(|\nabla u(rx)|}\over{|\nabla u(rx)|}}D^2 u(rx) |\Big]^2dx.
\end{eqnarray*}

\n We also introduce for each $\epsilon\in(0,1)$, the unique solution of the problem
\begin{equation*}
(P_\epsilon)\begin{cases} &  u_\epsilon-u\in W_0^{1,A}(B_1)\\
&  \Delta_{A_\epsilon}u_\epsilon= f H_\epsilon(u_\epsilon) \quad \text{in }\quad
B_1,
\end{cases}
\end{equation*}
\n where $A_\epsilon$ is the function defined in Section 2.

\vs 0,5cm
\begin{lemma}\label{l6.1} Under the assumptions (1.24)-(1.25), there exists two positive constants
$C_2=C_2(a_0,a_1,N,t_*,M_1)$ and $C_3=C_3(a_0,N,t_*,M_1)$
such that we have for every $u\in {\cal F}_A$
\begin{equation}\label{6.1}
E(1/2,u)  \leqslant C_2a^2(||\nabla u||_{B_{3/4},\infty})+ C_3a(||\nabla u||_{B_{3/4},\infty})\int_{B_{3/4}}
|\nabla f|dx.
\end{equation}
\end{lemma}

\vs 0,5cm \n \emph{Proof.} \n Let $u\in {\cal F}_A$  and consider
for each $\epsilon\in(0,1)$ the unique solution of the problem
$(P_\epsilon)$.
Note that $u_\epsilon \rightarrow u$ in $W^{1,A}(B_1)$ and that
$u_\epsilon \in W_{loc}^{2,2}(B_1)$ (see \cite{[CL3]}). Moreover we
have by taking $\varphi_{x_i}$ as a test function for $(P_\epsilon)$
and integrating by parts, where $\varphi\in {\cal D}(B_1)$

\begin{eqnarray}\label{6.2}
&& \int_{B_1} \nabla_{A_\epsilon}u_\epsilon.(\nabla \varphi_{x_i})dx=
-\int_{B_1} f H_\epsilon(u_\epsilon)\varphi_{x_i}dx \nonumber\\
&& \int_{B_1} (\nabla_{A_\epsilon}u_\epsilon)_{x_i}.\nabla \varphi dx= \int_{B_1} f
H_\epsilon(u_\epsilon)\varphi_{x_i}dx.
\end{eqnarray}

\n By density  (6.2) remains valid for $\varphi=G(u_{\epsilon
x_i})\zeta^2$, where $G(t)$ is smooth enough and $\zeta\in{\cal
D}(B_{3/4})$ such that
\begin{equation*}
\begin{cases}
& 0\leqslant \zeta \leqslant 1~~
\text{ in } B_{3/4}\\
&  \zeta=1 ~~
\text{ in } B_{1/2}\\
& \displaystyle{|\nabla \zeta|\leqslant 4~~ \text{ in } B_{3/4}}.
\end{cases}
\end{equation*}

\n Setting $t_\epsilon=(\epsilon+|\nabla u_\epsilon|^2)^{1/ 2}$, the
left hand side of (6.2) becomes
\begin{eqnarray}\label{6.3}
&&I^i= \int_{B_1} \Big({{a(t_\epsilon)}\over{t_\epsilon}} \nabla
u_\epsilon\Big)_{x_i}.\nabla
(G(u_{\epsilon x_i})\zeta^2)dx\nonumber\\
&&=\int_{B_1} {{a(t_\epsilon)}\over{t_\epsilon}}\Big[\nabla
u_{\epsilon
x_i}+\Big({{a'(t_\epsilon)}\over{a(t_\epsilon)}}t_\epsilon-1\Big){{\nabla
u_{\epsilon x_i}. \nabla u_{\epsilon}}\over{\epsilon+|\nabla
u_\epsilon|^2}}\nabla u_\epsilon\Big].\big[\zeta^2G'(u_{\epsilon
x_i})\nabla u_{\epsilon x_i}+2\zeta
G(u_{\epsilon x_i})\nabla\zeta\big]\nonumber\\
&&=\int_{B_1} \zeta^2G'(u_{\epsilon
x_i}){{a(t_\epsilon)}\over{t_\epsilon}}\Big[|\nabla u_{\epsilon
x_i}|^2+\Big({{a'(t_\epsilon)}\over{a(t_\epsilon)}}t_\epsilon-1\Big){{|\nabla
u_{\epsilon x_i}. \nabla
u_{\epsilon}|^2}\over{\epsilon+|\nabla u_\epsilon|^2}}\Big]dx\nonumber\\
&&+\int_{B_1} 2\zeta G(u_{\epsilon
x_i}){{a(t_\epsilon)}\over{t_\epsilon}}\nabla\zeta\Big[\nabla
u_{\epsilon
x_i}+\Big({{a'(t_\epsilon)}\over{a(t_\epsilon)}}t_\epsilon-1\Big){{\nabla
u_{\epsilon x_i}. \nabla u_{\epsilon}}\over{\epsilon+|\nabla
u_\epsilon|^2}}\nabla
u_\epsilon\Big]dx\nonumber\\
&&=I_1^i+I_2^i.
\end{eqnarray}

\vs 0.2cm\n According to the assumption (1.25), we shall discuss two cases:

\vs 0.2cm\n \emph{\underline{$1^{st}$ Case}} :
$~~t~\rightarrow~{{a(t)}\over t}$ is non-increasing in $(0,t_*)$

\vs 0.2cm\n Let $
\displaystyle{G(t)={{a((\epsilon+t^2)^{1/2})}\over{(\epsilon+t^2)^{1/2}}}t}$.
Then we have since $a_0\leqslant 1$

\begin{eqnarray*}G'(t)&=&{{a((\epsilon+t^2)^{1/2})}\over{(\epsilon+t^2)^{1/2}}}\Big[1+
\Big({{a'((\epsilon+t^2)^{1/2})}\over{a((\epsilon+t^2)^{1/2})}}(\epsilon+t^2)^{1/2}
-1\Big){{t^2}\over{\epsilon+t^2}}\Big]\\
&\geqslant &
a_0{{a((\epsilon+t^2)^{1/2})}\over{(\epsilon+t^2)^{1/2}}}.
\end{eqnarray*}

\n Let $s_\epsilon=(\epsilon+|u_{\epsilon x_i}|^2)^{1/2}$. Using Cauchy-Schwarz
inequality and (1.5), we get
\begin{equation}\label{6.4}
I_1^i\geqslant a_0^2\int_{B_1}\zeta^2
{{a(s_\epsilon)}\over{s_\epsilon}}{{a(t_\epsilon)}\over{t_\epsilon}}
|\nabla u_{\epsilon x_i}|^2dx.
\end{equation}

\n For $I_2^i$, we have by Young's inequality, since
$|\nabla\zeta|\leqslant 4$

\begin{eqnarray}\label{6.5}
|I_2^i| &\leqslant&\int_{B_1} 8(2+a_1)\zeta
{{a(s_\epsilon)}\over{s_\epsilon}}{{a(t_\epsilon)}\over{t_\epsilon}}
|u_{\epsilon x_i}||\nabla u_{\epsilon x_i}|dx\nonumber\\
&\leqslant&{{a_0^2}\over 2}\int_{B_1}
\zeta^2{{a(s_\epsilon)}\over{s_\epsilon}}{{a(t_\epsilon)}\over{t_\epsilon}}
|\nabla u_{\epsilon x_i}|^2dx\nonumber\\
&&+{{32(2+a_1)^2}\over{a_0^2}}\int_{B_{3/4}}
{{a(s_\epsilon)}\over{s_\epsilon}}{{a(t_\epsilon)}\over{t_\epsilon}}
|u_{\epsilon x_i}|^2dx\nonumber\\
&\leqslant&{{a_0^2}\over 2}\int_{B_1}
\zeta^2{{a(s_\epsilon)}\over{s_\epsilon}}{{a(t_\epsilon)}\over{t_\epsilon}}
|\nabla u_{\epsilon x_i}|^2dx\nonumber\\
&&+{{32(2+a_1)^2}\over{a_0^2}}\int_{B_{3/4}} (a(t_\epsilon))^2dx.
\end{eqnarray}

\n On the other hand, integrating by parts and using the
monotonicity of $H_\epsilon$, we have
\begin{eqnarray}\label{6.6}
I^i&=& \int_{B_1} f H_\epsilon(u_\epsilon)(G(u_{\epsilon
x_i})\zeta^2)_{x_i}dx\nonumber\\
&=&-\int_{B_1} f H'_\epsilon(u_\epsilon)u_{\epsilon
x_i}G(u_{\epsilon x_i})\zeta^2dx -\int_{B_1} f_{x_i}
H_\epsilon(u_\epsilon)G(u_{\epsilon x_i})\zeta^2dx\nonumber\\
&\leqslant& -\int_{B_1} f_{x_i} H_\epsilon(u_\epsilon)G(u_{\epsilon
x_i})\zeta^2dx.
\end{eqnarray}

\n Hence we get from (6.2)-(6.5) that
\begin{equation}\label{6.7}
\int_{B_1}\zeta^2
{{a(s_\epsilon)}\over{s_\epsilon}}{{a(t_\epsilon)}\over{t_\epsilon}}
|\nabla u_{\epsilon x_i}|^2dx \leqslant {{64(2+a_1)^2}\over
{a_0^4}}\int_{B_{3/4}} (a(t_\epsilon))^2dx+{2\over
{a_0^2}}\int_{B_{3/4}} |\nabla f|a(t_\epsilon)dx.
\end{equation}

\n Note that we have

\begin{eqnarray}\label{6.8}
&&\int_{B_{1/2}} \Big[{{a(t_\epsilon)}\over{t_\epsilon}} |\nabla u_{\epsilon x_i}|\Big]^2dx
=\int_{B_{1/2}\cap\{s_\epsilon\leqslant t_\epsilon< t_*\}}\Big[{{a(t_\epsilon)}\over{t_\epsilon}} |\nabla u_{\epsilon x_i}|\Big]^2dx\nonumber\\
&&+\int_{B_{1/2}\cap\{s_\epsilon<t_*\leqslant t_\epsilon\}}\Big[{{a(t_\epsilon)}\over{t_\epsilon}} |\nabla u_{\epsilon x_i}|\Big]^2dx
+\int_{B_{1/2}\cap\{t_*\leqslant s_\epsilon<\leqslant t_\epsilon \}}\Big[{{a(t_\epsilon)}\over{t_\epsilon}} |\nabla u_{\epsilon x_i}|\Big]^2dx.\nonumber\\
&&
\end{eqnarray}

\n Using the monotonicity of ${{a(t)}\over t}$ in $(0,t_*)$, and the uniform
boundness of $||t_\epsilon||_{B_{3/4},\infty}$, we get

\begin{equation}\label{6.9}
\int_{B_{1/2}\cap\{s_\epsilon\leqslant t_\epsilon< t_*\}}\Big[{{a(t_\epsilon)}\over{t_\epsilon}} |\nabla u_{\epsilon x_i}|\Big]^2dx
\leqslant\int_{B_{1/2}\cap\{s_\epsilon\leqslant t_\epsilon< t_*\}}{{a(s_\epsilon)}\over{s_\epsilon}}{{a(t_\epsilon)}\over{t_\epsilon}}
|\nabla u_{\epsilon x_i}|^2dx.
\end{equation}

\begin{eqnarray}\label{6.10}
&&\int_{B_{1/2}\cap\{s_\epsilon<t_*\leqslant t_\epsilon\}}\Big[{{a(t_\epsilon)}\over{t_\epsilon}} |\nabla u_{\epsilon x_i}|\Big]^2dx
={{t_*}\over{a(t_*)}}\int_{B_{1/2}\cap\{s_\epsilon<t_*\leqslant t_\epsilon\}}{{a(t_*)}\over{t_*}}{{a(t_\epsilon)}\over{t_\epsilon}}
{{a(t_\epsilon)}\over{t_\epsilon}}|\nabla u_{\epsilon x_i}|^2 dx\nonumber\\
&&\hskip3cm\leqslant
{{a(||t_\epsilon||_{B_{3/4},\infty})}\over{a(t_*)}}\int_{B_{1/2}\cap\{s_\epsilon<t_*\leqslant t_\epsilon\}}{{a(s_\epsilon)}\over{s_\epsilon}}{{a(t_\epsilon)}\over{t_\epsilon}}
|\nabla u_{\epsilon x_i}|^2dx.
\end{eqnarray}

\begin{eqnarray}\label{6.11}
&&\int_{B_{1/2}\cap\{t_*\leqslant s_\epsilon<\leqslant t_\epsilon\}}\Big[{{a(t_\epsilon)}\over{t_\epsilon}} |\nabla u_{\epsilon x_i}|\Big]^2dx
=\int_{B_{1/2}\cap\{t_*\leqslant s_\epsilon<\leqslant t_\epsilon\}}{{a(t_\epsilon)}\over{a(s_\epsilon)}}{{s_\epsilon}\over{t_\epsilon}}{{a(s_\epsilon)}\over{s_\epsilon}}{{a(t_\epsilon)}\over{t_\epsilon}}
|\nabla u_{\epsilon x_i}|^2dx\nonumber\\
&&\hskip3cm \leqslant{{a(||t_\epsilon||_{B_{3/4},\infty})}\over{a(t_*)}}\int_{B_{1/2}\cap\{t_*\leqslant s_\epsilon<\leqslant t_\epsilon\}}{{a(s_\epsilon)}\over{s_\epsilon}}{{a(t_\epsilon)}\over{t_\epsilon}}
|\nabla u_{\epsilon x_i}|^2dx.
\end{eqnarray}

\n Taking into account (6.8)-(6.11), we get
\begin{equation}\label{6.12}
\int_{B_{1/2}}\Big[{{a(t_\epsilon)}\over{t_\epsilon}} |\nabla u_{\epsilon x_i}|\Big]^2dx
\leqslant\max\Big(1,{{a(||t_\epsilon||_{B_{3/4},\infty})}\over{a(t_*)}}\Big)
\int_{B_{1/2}}{{a(s_\epsilon)}\over{s_\epsilon}}{{a(t_\epsilon)}\over{t_\epsilon}}
|\nabla u_{\epsilon x_i}|^2dx.
\end{equation}
Using the fact that $\zeta=1$ in $B_{1/2}$, we deduce from (6.7) and (6.12), by summing up from
$i=1$ to $i=N$ that

\begin{eqnarray}\label{6.13}
&&\int_{B_{1/2}} \Big[{{a(t_\epsilon)}\over{t_\epsilon}} |D^2
u_\epsilon|\Big]^2dx \leqslant {{64N(2+a_1)^2}\over
{a_0^4}}\max\Big(1,{{a(||t_\epsilon||_{B_{3/4},\infty})}\over{a(t_*)}}\Big)\int_{B_{3/4}} (a(t_\epsilon))^2dx\nonumber\\
&&\qquad\qquad+{{2N}\over{a_0^2}}\max\Big(1,{{a(||t_\epsilon||_{B_{3/4},\infty})}\over{a(t_*)}}\Big)\int_{B_{3/4}} |\nabla f| a(t_\epsilon)dx.
\end{eqnarray}

\n It follows from (6.13) and since $\nabla u_\epsilon$ is
uniformly bounded in $B_{3/4}$, that
${{a(t_\epsilon)}\over{t_\epsilon}} D^2 u_\epsilon$ is bounded
in $L^2(B_{1/2})$. So there exists a subsequence and a function $W\in L^2(B_{1/2})$
such that
$${{a(t_\epsilon)}\over{t_\epsilon}} D^2
u_\epsilon \rightharpoonup W ~~\text{in}~L^2(B_{1/2}).$$

\n Passing to the $\displaystyle{\liminf}$ in (6.13), we obtain by
taking into account the fact that $\nabla u_\epsilon$ converges uniformly,
up to a subsequence, to $\nabla u$ in $\overline{B_{{3/4}}}$
\begin{eqnarray*}
\int_{B_{1/2}} |W|^2dx&\leqslant&\liminf_{\epsilon\rightarrow
0}\int_{B_{1/2}}\Big[{{a(t_\epsilon)}\over{t_\epsilon}} |D^2
u_\epsilon|\Big]^2dx
\nonumber\\
&\leqslant&{{64N(2+a_1)^2}\over {a_0^4}}\max\Big(1,{{a(||\nabla u||_{B_{3/4},\infty})}\over{a(t_*)}}\Big)\int_{B_{3/4}} (a(|\nabla
u|))^2dx\nonumber\\
&&+{{2N}\over{a_0^2}}\max\Big(1,{{a(||\nabla u||_{B_{3/4},\infty})}\over{a(t_*)}}\Big)\int_{B_{3/4}} |\nabla f|a(|\nabla u|)dx.
\end{eqnarray*}

\n Since $\nabla u_\epsilon$ converges uniformly to $\nabla u$ in
$\overline{B_{{3/4}}}$, we deduce that $D^2 u \in
L_{loc}^2(B_{1/2}\cap\{\nabla u\neq 0\})$. Consequently we obtain
$\displaystyle{W={a(|\nabla u|)\over  |\nabla u|} D^2 u}$ a.e. in
$B_{1/2}\cap\{\nabla u\neq 0\}$, and therefore we get
\begin{eqnarray*}
E(1/2,u)  &\leqslant& {{64N(2+a_1)^2}\over {a_0^4
}}\big({8\over 3}\big)^N \max\Big(1,{{a(M_1)}\over{a(t_*)}}\Big)a^2(|\nabla u|_{B_{3/4},\infty})\\
&&+ {{2N}\over {a_0^2 |B_{1/2}|}}\max\Big(1,{{a(M_1)}\over{a(t_*)}}\Big)a(|\nabla
u|_{B_{3/4},\infty})\int_{B_{3/4}} |\nabla f|dx \nonumber\\
&\leqslant& C_2a^2(|\nabla u|_{B_{3/4},\infty})+ C_3a(|\nabla u|_{B_{3/4},\infty})\int_{B_{3/4}}
|\nabla f|dx.
\end{eqnarray*}

\vs 0.2cm\n \emph{\underline{$2^{nd}$ Case}} :
$~~t~\rightarrow~\displaystyle{{a(t)}\over t}$ is non-decreasing in $(0,t_*)$

\vs 0.2cm\n Let $G(t)=t$. Using (1.5), we obtain

\begin{equation}\label{6.14}
I_1^i\geqslant a_0\int_{B_1}\zeta^2 {{a(t_\epsilon)}\over{t_\epsilon}}
|\nabla u_{\epsilon x_i}|^2dx.
\end{equation}

\n For $I_2^i$, we have by Young's inequality, since
$|\nabla\zeta|\leqslant 4$

\begin{eqnarray}\label{6.15}
|I_2^i| &\leqslant&\int_{B_1} 8(2+a_1)\zeta
{{a(t_\epsilon)}\over{t_\epsilon}}
|u_{\epsilon x_i}||\nabla u_{\epsilon x_i}|dx\nonumber\\
&\leqslant&{1\over 2}\int_{B_1}
\zeta^2{{a(t_\epsilon)}\over{t_\epsilon}} |\nabla u_{\epsilon
x_i}|^2dx+32(2+a_1)^2\int_{B_{3/4}}
{{a(t_\epsilon)}\over{t_\epsilon}}
|u_{\epsilon x_i}|^2dx\nonumber\\
&\leqslant&{1\over 2}\int_{B_1}
\zeta^2{{a(t_\epsilon)}\over{t_\epsilon}}
|\nabla u_{\epsilon x_i}|^2dx+32(2+a_1)^2\int_{B_{3/4}} t_\epsilon a(t_\epsilon)dx.
\end{eqnarray}

\n Integrating by parts and using the monotonicity of $H_\epsilon$,
we obtain as in the previous case
\begin{eqnarray}\label{6.16}
I^i&=&\int_{B_1} f
H_\epsilon(u_\epsilon)(u_{\epsilon x_i}\zeta^2)_{x_i}dx\nonumber\\
&=&-\int_{B_1}f H'_\epsilon(u_\epsilon)(u_{u_\epsilon
x_i})^2\zeta^2dx
-\int_{B_1}f_{x_i} H_\epsilon(u_\epsilon)u_{\epsilon x_i}\zeta^2dx\nonumber\\
&\leqslant&-\int_{B_1}f_{x_i} H_\epsilon(u_\epsilon)u_{\epsilon
x_i}\zeta^2dx.
\end{eqnarray}

\n Hence we get from (6.2) and (6.14)-(6.15) that
\begin{equation}\label{6.17}
\int_{B_{1/2}} {{a(t_\epsilon)}\over{t_\epsilon}} |D^2
u_{\epsilon}|^2dx \leqslant 64N(2+a_1)^2\int_{B_{3/4}} t_\epsilon
a(t_\epsilon)dx+2\int_{B_{3/4}} |\nabla f||\nabla u_\epsilon|dx.
\end{equation}

\n Using the monotonicity of $\displaystyle{{{a(t)}\over t}}$ in $(0,t_*)$, we
get
\begin{equation*}{{a(t_\epsilon)}\over{t_\epsilon}}\leqslant\left\{
  \begin{array}{ll}
{{a(||t_\epsilon||_{B_{3/4},\infty})}\over{||t_\epsilon||_{B_{3/4},\infty}}} & \hbox{if}\quad ||t_\epsilon||_{B_{3/4},\infty}<t_*\\
{{a(t_*)}\over{t_*}}  & \hbox{if}\quad ||t_\epsilon||_{B_{3/4},\infty}\geqslant t_*
  \end{array}
\right\}\leqslant {{a(||t_\epsilon||_{B_{3/4},\infty})}\over{||t_\epsilon||_{B_{3/4},\infty}}}\max\Big(1,{||t_\epsilon||_{B_{3/4},\infty}\over{t_*}}\Big).
\end{equation*}

\n We deduce from (6.17) that
\begin{eqnarray*}
&&\int_{B_{1/2}}\Big[{{a(t_\epsilon)}\over{t_\epsilon}} |D^2
u_\epsilon|\Big]^2dx=\int_{B_{1/2}}{{a(t_\epsilon)}\over{t_\epsilon}}.
{{a(t_\epsilon)}\over{t_\epsilon}} |D^2 u_\epsilon|^2dx\nonumber\\
&&\hskip1cm\leqslant {{a(||t_\epsilon||_{B_{3/4},\infty})}\over{||t_\epsilon||_{B_{3/4},\infty}}}\max\Big(1,{||t_\epsilon||_{B_{3/4},\infty}\over{t_*}}\Big)
\int_{B_{1/2}}{{a(t_\epsilon)}\over{t_\epsilon}} |D^2 u_\epsilon|^2dx\nonumber\\
&&\hskip1cm\leqslant 64N(2+a_1)^2{{a(||t_\epsilon||_{B_{3/4},\infty})}\over{||t_\epsilon||_{B_{3/4},\infty}}}\max\Big(1,{||t_\epsilon||_{B_{3/4},\infty}\over{t_*}}\Big)
\int_{B_{3/4}} t_\epsilon
a(t_\epsilon)dx\nonumber\\
&&\hskip1cm~~+2{{a(||t_\epsilon||_{B_{3/4},\infty})}\over{||t_\epsilon||_{B_{3/4},\infty}}}\max\Big(1,{||t_\epsilon||_{B_{3/4},\infty}\over{t_*}}\Big)
\int_{B_{3/4}} |\nabla f||\nabla u_\epsilon|dx\nonumber\\
&&\hskip1cm\leqslant 64N(2+a_1)^2|B_{3/4}|\max\Big(1,{||t_\epsilon||_{B_{3/4},\infty}\over{t_*}}\Big)a^2(||t_\epsilon||_{B_{3/4},\infty})\nonumber\\
&&\hskip1cm~~ +2\max\Big(1,{||t_\epsilon||_{B_{3/4},\infty}\over{t_*}}\Big)a(||t_\epsilon||_{B_{3/4},\infty})\int_{B_{3/4}} |\nabla f|dx.
\end{eqnarray*}

\n Passing to the limit, as in the first case, we obtain

\begin{eqnarray*}
E(1/2,u)  &\leqslant& \max\Big(1,{  {||\nabla u||_{B_{3/4},\infty} }\over{t_*}}\Big)
\Big({{64N(2+a_1)^2|B_{3/4}|}\over {|B_{1/2}|}}a^2(||\nabla
u||_{B_{3/4},\infty})\\
&&+{2\over {|B_{1/2}|}}a(||\nabla u||_{B_{3/4},\infty})\int_{B_{3/4}} |\nabla f|dx\Big)\nonumber\\
&\leqslant& C_2a^2(||\nabla u||_{B_{3/4},\infty})+ C_3 a(||\nabla u||_{B_{3/4},\infty})\int_{B_{3/4}}
|\nabla f|dx.
\end{eqnarray*}

\qed

 \vs 0,5cm \n Now we prove that $E(r,u)$ is uniformly bounded under the
assumption (1.24)-(1.25). We observe that (1.24 means that
\begin{equation}\label{6.18}
\exists K_0>0 :\quad \forall r\in(0,3/4)\quad \int_{B_r} |\nabla f|dx
\leqslant K_0 r^{N-1}.
\end{equation}
In particular, (6.18) is satisfied if $f\in C^{0,1}(\overline{B}_{3/4})$.

\begin{lemma}\label{l6.2} There exists a positive  constant
$C_4=C_4(a_0,a_1,a(1),N,M_0,M_1,\Lambda_0,K_0)$ such that we have
\begin{equation*}
E(r,u) \leqslant C_4~~~\forall u\in {\cal F}_A~~~\forall r\in
(0,{1\over 2}).
\end{equation*}

\end{lemma}

\vs 0,5cm \n \emph{Proof.} Note that it is enough to prove the lemma
for $r\in (0,{1\over 4})$. Indeed for $r\in [{1\over 4},{1\over
2})$, we have by Lemma 6.1
\begin{eqnarray*}
E(r,u)&=&{1\over{|B_r|}}\int_{B_{r}\cap\{\nabla u\neq 0\}}
\Big[{a(|\nabla u|)\over |\nabla
u|} |D^2 u|\Big]^2dx \\
&\leqslant& {1\over{|B_{1/4}|}}\int_{B_{1/2}\cap\{\nabla u\neq 0\}}
\Big[{a(|\nabla u|)\over |\nabla
u|} |D^2 u|\Big]^2dx \\
&=& {|B_{1/2}|\over{|B_{1/4}|}}E(1/2,u)\\
&\leqslant&2^N\big(C_1a^2(||\nabla u||_{B_{3/4},\infty})+ C_2a(||\nabla u||_{B_{3/4},\infty})\int_{B_{3/4}}
|\nabla f|dx\big)\\
&\leqslant& C_4=C_4(a_0,a_1,a(1),N,M_0,K_0).
\end{eqnarray*}

\n For $r\in (0,{1\over 4})$, we consider the function $v_r(x)=
\displaystyle{ {u( 2rx) }\over {2r}}$ defined  in $B_1$.  We have
by definition of $v_r$, Propositions 5.1, 5.2 and (2.5)-(2.6),
we have
\begin{align}\label{6.19-6.21}
& 0\leqslant v_r\leqslant C_0\displaystyle{ {\widetilde{A}(2r)} \over
{2r}}\leqslant C_0 2^{1+{1\over{a_0}}} a^{-1}(r)\quad \hbox{ in } B_1,
\\
&|\nabla v_r(x)|=|\nabla u(2rx)|\leqslant C_1 a^{-1}(2r)\leqslant C_12^{1\over{a_0}}  a^{-1}(r)\quad \hbox{ in } B_1,\\
&D^2 v_r(x)=2r(D^2 u)(2rx).
\end{align}
\n Using (6.21)-(6.22), we compute
\begin{eqnarray}\label{5.16}
E({1\over 2},v_r)&=&{1\over{|B_1|}}\int_{{B_1\cap\{\nabla
v_r({1\over 2}x)\neq 0\}}} \Big[{a(|\nabla v_r({1\over 2}x)|)\over
|\nabla v_r({1\over 2}x)|} \big|D^2 v_r\big({1\over
2}x\big)\big|\Big]^2 dx\nonumber
\\
&=&{1\over{|B_1|}}\int_{{B_1\cap\{\nabla u(rx)\neq 0\}}}
\Big[{a(|\nabla u(rx)|)\over |\nabla u(rx)|} 2r|D^2
u(rx)|\Big]^2dx\nonumber
\\
&=&4r^2E(r,u).
\end{eqnarray}

\n Moreover, we have
\begin{equation*}
\Delta_A v_r(x) =2r(\Delta_A
u)(2rx)=2rf(2rx)\chi_{\{u(2rx)>0\}}=f_r(x) \chi_{\{v_r(x)>0\}}.
  \end{equation*}

\n Using (5.2) and (6.18), we obtain from Lemma 6.1

\begin{eqnarray}\label{5.18}
E({1\over 2},v_r)& \leqslant& C_2a^2(||\nabla
v_r||_{B_{3/4},\infty})+ C_3a(||\nabla
v_r||_{B_{3/4},\infty})\int_{B_{3/4}}
|\nabla f_r|dx\nonumber\\
&\leqslant& C_2(a(C_12^{1\over{a_0}}a^{-1}(r)))^2+4r^2C_3a(C_12^{1\over{a_0}}
a^{-1}(r))\int_{B_{3/4}} |\nabla f(2rx)|dx\nonumber\\
&\leqslant& C_2(\max(2C_1^{a_0},2C_1^{a_1}))^2a^2( a^{-1}(r))\nonumber\\
&&+2^{2-N}r^2C_3\max(2C_1^{a_0},2C_1^{a_1})rr^{-N}\int_{B_{3r/2}} |\nabla f(y)|dy\nonumber\\
&\leqslant&4 C_2(\max(C_1^{a_0},C_1^{a_1}))^2r^2+ 2^{3-N}C_3 \max(C_1^{a_0},C_1^{a_1}) K_0
r^2\nonumber\\
&=& 4r^2\max(C_1^{a_0},C_1^{a_1}) (\max(C_1^{a_0},C_1^{a_1}) C_2+2^{1-N}C_3 K_0)=
4C_4 r^2.
\end{eqnarray}

\n Taking into account (6.22) and (6.23), we get
\begin{eqnarray*}
E(r,u)& \leqslant& C_4=C_4(N,a_0,a_1,a(1),M_0,\Lambda_0,K_0).
\end{eqnarray*}

\qed

\vs 0.5cm
\begin{lemma}\label{l6.3} We have
\begin{equation}\label{6.24}
\lambda_0^2{{(H_\epsilon(u_\epsilon))^2}\over{(2+a_1)^2}}\leqslant
\Big[{{a((\epsilon+|\nabla u_\epsilon|^2)^{1/2})}\over
{(\epsilon+|\nabla u_\epsilon|^2)^{1/2}}} |D^2
u_\epsilon|\Big]^2~~\text{a.e. in}~~B_1.
\end{equation}
\end{lemma}

\vs 0,5cm

\n \emph{Proof.} Since $u_\epsilon\in W_{loc}^{2,2}(B_1)$, we obtain
from $(P_\epsilon)$ by using (1.5)

\begin{eqnarray*}
(\lambda_0H_\epsilon(u_\epsilon))^2&=&
\Big(div\Big({{a((\epsilon+|\nabla u_\epsilon|^2)^{1/2})}\over
{(\epsilon+|\nabla u_\epsilon|^2)^{1/2}}}\nabla
u_\epsilon\Big)\Big)^2\\
&=&\Big({{a((\epsilon+|\nabla u_\epsilon|^2)^{1/2})}\over
{(\epsilon+|\nabla u_\epsilon|^2)^{1/2}}}\Big)^2 \Big[\Delta
u_{\epsilon}+\Big({{a'(t_\epsilon)}\over{a(t_\epsilon)}}t_\epsilon-1\Big)
{{D^2u_\epsilon.\nabla u_{\epsilon}}\over{\epsilon+|\nabla
u_\epsilon|^2}}\nabla u_\epsilon\Big]^2\\
&\leqslant&\Big({{a((\epsilon+|\nabla u_\epsilon|^2)^{1/2})}\over
{(\epsilon+|\nabla u_\epsilon|^2)^{1/2}}}\Big)^2
[|\Delta u_{\epsilon}|+(1+a_1)|D^2u_\epsilon|]^2\\
&\leqslant&(2+a_1)^2\Big[{{a((\epsilon+|\nabla
u_\epsilon|^2)^{1/2})}\over {(\epsilon+|\nabla
u_\epsilon|^2)^{1/2}}}|D^2u_\epsilon|\Big]^2.
\end{eqnarray*}

\n Hence (6.24) holds.\qed

\vs 0.5cm
\begin{lemma}\label{l6.4} There exists a
positive constant $C_5=C_5(a_0,a_1,a(1),N,M_0,\lambda_0,\Lambda_0,K_0)$ such that for any $0<\delta<a(t_*/2)$
and $r<1/4$ with $B_{2r}(x_0)\subset B_1$ and $x_0\in \partial
\{u>0\}\cap B_{1/2}$, we have
\begin{equation}\label{6.25}
{\cal L}^N(O_\delta\cap B_r(x_0)\cap \{u>0\})\leqslant
C_5\delta r^{N-1},
\end{equation}
\n where $O_\delta=\{|\nabla u|<a^{-1}(\delta)\}\cap B_{1/2}$.
\end{lemma}

\vs 0,5cm

\n \emph{Proof.} Let $r<1/4$, $\delta>0$ and $x_0\in
\partial \{u>0\}\cap B_{1/2}$  such that $B_{2r}(x_0)\subset B_1$. We choose a function
$\zeta\in{\cal D}(B_1)$ such that
\begin{equation*}
\begin{cases}
& 0\leqslant \zeta \leqslant 1~~
\text{ in } B_{2r}(x_0)\\
&  \zeta=1 ~~
\text{ in } B_r(x_0)\\
& \displaystyle{|\nabla \zeta|\leq {2\over r}~~ \text{ in }
B_{2r}(x_0)}.
\end{cases}
\end{equation*}

\n We denote by $u_\epsilon$ the unique solution of the problem $(P_\epsilon)$.
Let $G$ be a function such that $G(u_{\epsilon x_i})\in
W^{1,2}(B_{2r}(x_0))$. Arguing as in the proof of Lemma 5.1, we get
for $\varphi=G(u_{\epsilon x_i})\zeta$

\begin{eqnarray*}
&& \int_{B_1} \nabla_{A_\epsilon}u_\epsilon.\nabla \varphi_{x_i}dx=
-\int_{B_1}  fH_\epsilon(u_\epsilon)\varphi_{x_i}dx.
\end{eqnarray*}

\n Integrating by parts, we get by taking into account the monotonicity
of $H_\epsilon$

\begin{eqnarray}\label{6.26}
&&\int_{B_{2r}(x_0)} (\nabla_{A_\epsilon}u_\epsilon)_{x_i}.\nabla
(G(u_{\epsilon x_i})\zeta)dx=\int_{B_{2r}(x_0)}fH_\epsilon(u_\epsilon)\big(G(u_{\epsilon x_i})\zeta\big)_{x_i}dx\nonumber\\
&&\hskip 4cm=-\int_{B_{2r}(x_0)}fH_\epsilon'(u_\epsilon)u_{\epsilon
x_i}G(u_{\epsilon x_i})\zeta dx\nonumber\\
&&\hskip4.5cm-\int_{B_{2r}(x_0)}f_{x_i}H_\epsilon(u_\epsilon)G(u_{\epsilon
x_i})\zeta dx\nonumber\\
&&\hskip4cm\leqslant-\int_{B_{2r}(x_0)}f_{x_i}H_\epsilon(u_\epsilon)G(u_{\epsilon
x_i})\zeta dx.
\end{eqnarray}

\n Note that
\begin{eqnarray*}
&&(\nabla_{A_\epsilon}u_\epsilon)_{x_i}={{a(t_\epsilon)}\over{t_\epsilon}}\Big[\nabla
u_{\epsilon
x_i}+\Big({{a'(t_\epsilon)}\over{a(t_\epsilon)}}t_\epsilon-1\Big){{\nabla
u_{\epsilon x_i}. \nabla u_{\epsilon}}\over{\epsilon+|\nabla
u_\epsilon|^2}}\nabla u_\epsilon\Big].
\end{eqnarray*}

\n Then (6.26) becomes
\begin{eqnarray}\label{6.27}
&&\int_{B_{2r}(x_0)} \zeta G'(u_{\epsilon
x_i}){{a(t_\epsilon)}\over{t_\epsilon}}\Big[|\nabla u_{\epsilon x_i}|^2+\Big({{a'(t_\epsilon)}\over{a(t_\epsilon)}}t_\epsilon-1\Big){{|\nabla
u_{\epsilon x_i}. \nabla u_{\epsilon}|^2}\over{\epsilon+|\nabla
u_\epsilon|^2}}\Big]dx\nonumber\\
&&~~\leqslant-\int_{B_{2r}(x_0)} G(u_{\epsilon
x_i}){{a(t_\epsilon)}\over{t_\epsilon}}\Big[\nabla u_{\epsilon
x_i}+\Big({{a'(t_\epsilon)}\over{a(t_\epsilon)}}t_\epsilon-1\Big){{\nabla
u_{\epsilon x_i}. \nabla u_{\epsilon}}\over{\epsilon+|\nabla
u_\epsilon|^2}}\nabla u_\epsilon\Big].\nabla\zeta dx\nonumber\\
&&\qquad-\int_{B_{2r}(x_0)}f_{x_i}H_\epsilon(u_\epsilon)G(u_{\epsilon
x_i})\zeta dx\nonumber\\
&&\qquad\leqslant {2\over r} a_1\int_{B_{2r}(x_0)}
|G(u_{\epsilon x_i})|{{a(t_\epsilon)}\over{t_\epsilon}}|\nabla
u_{\epsilon x_i}|
dx-\int_{B_{2r}(x_0)}f_{x_i}H_\epsilon(u_\epsilon)G(u_{\epsilon
x_i})\zeta dx.\nonumber\\
&&
\end{eqnarray}

\n Now we have
\begin{eqnarray}\label{6.28}
&&\int_{B_{2r}(x_0)} |G(u_{\epsilon
x_i})|{{a(t_\epsilon)}\over{t_\epsilon}}|\nabla u_{\epsilon x_i}| dx
\leqslant  \int_{B_{2r}(x_0)} |G(u_{\epsilon x_i})- G(u_{ x_i})|
{{a(t_\epsilon)}\over{t_\epsilon}}|\nabla u_{\epsilon x_i}| dx \nonumber\\
&& ~~~+ \int_{B_{2r}(x_0)} |G(u_{
x_i})|{{a(t_\epsilon)}\over{t_\epsilon}}|\nabla
u_{\epsilon x_i}| dx \nonumber\\
&&~~~ = J^i_1 + J^i_2.
\end{eqnarray}

\vs 0.2cm\n According to the assumption (1.25), we shall discuss two cases:

\vs 0.2cm\n \emph{\underline{$1^{st}$ Case}} :
$~~\displaystyle{{{a(t)}\over t}}$ is non-increasing in $(0,t_*)$

\vs 0.2cm\n \n For each $\epsilon>0$ and
$\eta=a(2a^{-1}(\delta))$, we consider the function
\begin{equation*}
G(t)=\begin{cases} \displaystyle{
{{a((\epsilon+(a^{-1}(\eta))^2)^{1/2})}\over
{(\epsilon+(a^{-1}(\eta))^2)^{1/2}}}a^{-1}(\eta)} ~~
&\text{ if }~ t> a^{-1}(\eta)\\
 \displaystyle{ {{a((\epsilon+t^2)^{1/2})}\over {(\epsilon+t^2)^{1/2}}}t} ~~
&\text{ if } ~|t|\leqslant a^{-1}(\eta)\\
\displaystyle{ -{{a((\epsilon+(a^{-1}(\eta))^2)^{1/2})}\over
{(\epsilon+(a^{-1}(\eta))^2)^{1/2}}}a^{-1}(\eta)} ~~ &\text{ if }~
t<-a^{-1}(\eta).
\end{cases}
\end{equation*}

\n $G$ is Lipschitz continuous and we have

\begin{eqnarray*}G'(t)&=&{{a((\epsilon+t^2)^{1/2})}\over{(\epsilon+t^2)^{1/2}}}\Big[1+
\Big({{a'((\epsilon+t^2)^{1/2})}\over{a((\epsilon+t^2)^{1/2})}}(\epsilon+t^2)^{1/2}
-1\Big){{t^2}\over{\epsilon+t^2}}\Big]\chi_{\{|t|<
a^{-1}(\eta)\}}
\end{eqnarray*}

$$a_0 {{a((\epsilon+t^2)^{1/2})}\over{(\epsilon+t^2)^{1/2}}}\chi_{\{|t|<
a^{-1}(\eta)\}} \leqslant G'(t)\leqslant a_1
{{a((\epsilon+t^2)^{1/2})}\over{(\epsilon+t^2)^{1/2}}}\chi_{\{|t|<
a^{-1}(\eta)\}}.$$

\n Let $t_\epsilon=(\epsilon+|\nabla u_\epsilon|^2)^{1/
2}$. Since $\zeta=1$ in $B_{r/2}$, and $\{|\nabla u_{\epsilon}|<a^{-1}(\eta)\}\subset\{|u_{\epsilon x_i}|<a^{-1}(\eta)\}$,
we obtain from (6.27)-(6.28) and the monotonicity of ${{a(t)}\over t}$, since
$a^{-1}(\eta)=2a^{-1}(2\delta)<t_*$

\begin{eqnarray}\label{6.29}
&&a_0^2\int_{B_{r}\cap\{|\nabla u_{\epsilon}|<a^{-1}(\eta)\}}
\Big({{a(t_\epsilon)}\over{t_\epsilon}}\Big)^2 |\nabla u_{\epsilon
x_i}|^2dx\leqslant a_0^2\int_{B_{r}\cap\{|u_{\epsilon x_i}|<a^{-1}(\eta)\}}
\Big({{a(t_\epsilon)}\over{t_\epsilon}}\Big)^2 |\nabla u_{\epsilon
x_i}|^2dx\nonumber\\
&&\hskip 1cm \leqslant a_0^2\int_{B_{r}\cap\{|u_{\epsilon
x_i}|<a^{-1}(\eta)\}} {{a((\epsilon+|u_{\epsilon
x_i}|^2)^{1/2})}\over{(\epsilon+|u_{\epsilon
x_i}|^2)^{1/2}}}{{a(t_\epsilon)}\over{t_\epsilon}} |\nabla
u_{\epsilon
x_i}|^2dx\nonumber\\
&&\hskip 1cm \leqslant{2\over r} a_1( J^i_1 +
J^i_2)+\int_{B_{2r}(x_0)}|\nabla f||G(u_{\epsilon x_i})| dx.
\end{eqnarray}

\n Since $|G(t)|\leqslant a((\epsilon+(a^{-1}(\eta))^2)^{1/2}$, we
obtain by taking into account (6.18)
\begin{eqnarray}\label{6.30}
\int_{B_{2r}(x_0)}|\nabla f||G(u_{\epsilon x_i})| dx
&\leqslant&a((\epsilon+(a^{-1}(\eta))^2)^{1/2}\int_{B_{2r}(x_0)}|\nabla
f|dx\nonumber\\
&\leqslant& K_0a((\epsilon+(a^{-1}(\eta))^2)^{1/2}(2r)^{N-1}.
\end{eqnarray}

\n We claim that $|G(u_{\epsilon x_i})- G(u_{ x_i})| \longrightarrow
0$ in $L^2(B_{2r}(x_0))$, as $\epsilon\rightarrow 0$.

\vs 0.2cm \n Indeed, first we have $|G(u_{\epsilon x_i})- G(u_{
x_i})|\leqslant \eta$ $\forall x\in B_{2r}(x_0)$.

\n Let $x\in B_{2r}(x_0)$. We discuss two cases :

\n $\ast$ If $u_{ x_i}(x) =0$, then there exists by the uniform convergence of
$\nabla u_\epsilon$ to $\nabla u$ in $\overline{B}_{1/2}$, $\epsilon_0>0$ such
that $\forall\epsilon\in(0,\epsilon_0)$, $|u_{\epsilon x_i}(x)| <
a^{-1}(\eta)$. So
 $ \displaystyle{G(u_{\epsilon x_i})={{a((\epsilon+|u_{
x_i}|^2)^{1/2})}\over {(\epsilon+|u_{ x_i}|^2)^{1/2}}} u_{\epsilon
x_i}}$, $G(u_{ x_i}(x))=0$ and

\begin{eqnarray*}
&&|G(u_{\epsilon x_i})- G(u_{\epsilon x_i})| =
\Big|{{a((\epsilon+|u_{\epsilon x_i}|^2)^{1/2})}\over
{(\epsilon+|u_{\epsilon x_i}(x)|^2)^{1/2}}} u_{\epsilon
x_i}(x) -0\Big|\\
&&\qquad\qquad\qquad\qquad\leqslant a((\epsilon+|u_{\epsilon
x_i}|^2)^{1/2})\longrightarrow a(|u_{ x_i}|)=0\\
&&\qquad\qquad\qquad\qquad\qquad\epsilon\rightarrow 0.
\end{eqnarray*}

\n $\ast$ If $u_{ x_i}(x) \neq 0$, then $|u_{ x_i}(x)|>0$ and there
exists an $\epsilon_0>0$ such that $\forall\epsilon\in(0,\epsilon_0)$,
$|u_{ x_i}(x)|/2 < |u_{\epsilon x_i}(x)| < 3|u_{ x_i}(x)|/2$. So

\begin{eqnarray*}
|G(u_{\epsilon x_i})- G(u_{ x_i})| &=& | G^\prime(\theta_\epsilon)|.
| u_{\epsilon x_i}(x)- u_{ x_i}(x)|  \qquad \hbox{ with }
\theta_\epsilon
\longrightarrow u_{ x_i}(x)\\
&\leqslant &a_1 {{a(\theta_\epsilon)}\over
{\theta_\epsilon}} | u_{\epsilon x_i}(x)- u_{ x_i}(x)|\\
&\longrightarrow & a_1 {{a(|u_{ x_i}(x)|)}\over {|u_{
x_i}(x)|}} | u_{ x_i}(x)- u_{ x_i}(x)|=0 .
\end{eqnarray*}

\n On the other hand, we know (see proof of Lemma 6.1) that

$$\int_{B_{2r}(x_0)} \Big[ {{a((\epsilon+|\nabla u_\epsilon|^2)^{1/
2})}\over{(\epsilon+|\nabla u_\epsilon|^2)^{1/ 2}}}|\nabla
u_{\epsilon x_i}|\Big]^2 dx\leqslant C.$$

\n It follows that up to a subsequence

$${{a((\epsilon+|\nabla u_\epsilon|^2)^{1/
2})}\over{(\epsilon+|\nabla u_\epsilon|^2)^{1/ 2}}}\nabla
u_{\epsilon x_i}\rightharpoonup W_i \qquad \hbox{in }
L^2(B_{2r}(x_0)) .$$

\n Therefore
\begin{equation}\label{6.31}
J^i_1=\int_{B_{2r}(x_0)} |G(u_{\epsilon x_i})- G(u_{ x_i})|
{{a((\epsilon+|\nabla u_\epsilon|^2)^{1/
2})}\over{(\epsilon+|\nabla u_\epsilon|^2)^{1/ 2}}}|\nabla
u_{\epsilon x_i}| dx \longrightarrow
0,~~\text{as}~~\epsilon\rightarrow 0.
\end{equation}

\n For $J^i_2$, we have

\begin{eqnarray}\label{6.32}
&&J^i_2= \int_{B_{2r}\cap \{\nabla u \neq 0\}} |G(u_{ x_i})|
.{{a((\epsilon+|\nabla u_\epsilon|^2)^{1/
2})}\over{(\epsilon+|\nabla u_\epsilon|^2)^{1/ 2}}}|\nabla
u_{\epsilon x_i}|  \qquad \hbox{ since  } G(0)=0\nonumber
\\
&& \leqslant a((\epsilon+(a^{-1}(\eta))^2)^{1/2}) \int_{B_{2r}(x_0)\cap \{\nabla u \neq 0\}}
{{a((\epsilon+|\nabla u_\epsilon|^2)^{1/ 2})}\over{(\epsilon+|\nabla
u_\epsilon|^2)^{1/ 2}}}|\nabla
u_{\epsilon x_i}|\nonumber\\
&& \leqslant a((\epsilon+(a^{-1}(\eta))^2)^{1/2}) |B_{2r}(x_0)|^{1/2}.
 \Big( \int_{B_{2r}(x_0)\cap \{\nabla u
\neq 0\}} \Big[{{a((\epsilon+|\nabla u_\epsilon|^2)^{1/
2}))}\over{(\epsilon+|\nabla u_\epsilon|^2)^{1/ 2}}}|D^2
u_\epsilon|\Big]^2\Big)^{1/2}.\nonumber\\
&&
\end{eqnarray}

\vs 0,5cm \n We claim that $O_{\delta}\subset \{|\nabla u_\epsilon|<a^{-1}(\eta)\}$.
Indeed since $\nabla u_\epsilon$ converges uniformly to $\nabla u$
in $\overline{B}_{1/2}$, there exists $\epsilon_0>0$ such that
$$\forall \epsilon\in(0,\epsilon_0),~~~
|\nabla u_\epsilon-\nabla
u|_{\infty,\overline{B}_{1/2}}<a^{-1}(\delta)/2.$$

\n We deduce that for $x\in \cap O_{\delta}$,

\begin{eqnarray*}
\forall \epsilon\in(0,\epsilon_0),~~~ |\nabla
u_\epsilon(x)|&\leqslant&|\nabla u_\epsilon(x)-\nabla u(x)|+|\nabla
u(x)|\\
&<&a^{-1}(\delta)/2+a^{-1}(\delta)\\
&<&2a^{-1}(\delta)=a^{-1}(\eta).
\end{eqnarray*}

\n We obtain from (6.29)

\begin{eqnarray}\label{6.33}
a_0^2\int_{B_r(x_0)\cap O_\delta}
\Big({{a(t_\epsilon)}\over{t_\epsilon}}\Big)^2 |\nabla u_{\epsilon
x_i}|^2dx&\leqslant& {2\over r}
a_1(J_1^i+J_2^i)\nonumber\\
&&+K_02^{N-1}a((\epsilon+(a^{-1}(\eta))^2)^{1/2}r^{N-1}.
\end{eqnarray}

\n Summing up from $i=1$ to $N$ in (6.33) and
using (6.1) and (6.32), we get
\begin{eqnarray}\label{6.34}
&&{{a_0^2}\over{(2+a_1)^2}}\int_{B_r(x_0)\cap O_\delta}
(\lambda_0 H_\epsilon(u_\epsilon)^2dx\leqslant {2a_1\over
r}\Big(\sum_{i=1}^{i=N}J_1^i\nonumber\\
&&+Na((\epsilon+(a^{-1}(\eta))^2)^{1/2}) |B_{2r}(x_0)|^{1/2}. \Big(
\int_{B_{2r}(x_0)\cap \{\nabla u \neq 0\}}
\Big[{{a(t_\epsilon)}\over{t_\epsilon}} |D^2
u_\epsilon|\Big]^2dx\Big)^{1/2}\Big)\nonumber\\
&&+NK_02^{N-1}a((\epsilon+(a^{-1}(\eta))^2)^{1/2})r^{N-1}.
\end{eqnarray}

\n Letting $\epsilon\rightarrow 0$ in (6.34), and using (6.31) and
Lemma 6.2, we obtain
\begin{eqnarray*}
&&{{\lambda_0^2a_0^2}\over{(2+a_1)^2}}\int_{B_r(x_0)\cap
O_\delta\cap\{u>0\}} dx\nonumber\\
&&~~~~\leqslant {2a_1\over r}N\eta |B_{2r}(x_0)|^{1/2}.
 \Big( \int_{B_{2r}(x_0)\cap \{\nabla u
\neq 0\}} \Big[{{a(|\nabla u|)}\over{|\nabla u|}} |D^2
u|\Big]^2dx\Big)^{1/2}\nonumber\\
&&~~~~~+NK_02^{N-1}\eta
r^{N-1}\nonumber\\
&&~~~~\leqslant N\eta {2\over r}
|B_{2r}(x_0)|^{1/2}.(C(a_0,a_1,a(1),N,M_0,\Lambda_0,K_0)|B_{2r}(x_0)|)^{1/2}\nonumber\\
&&~~~~~~+NK_02^{N-1}\eta
r^{N-1},
\end{eqnarray*}

\n which can be written
\begin{eqnarray*}
{\cal L}^N(O_\delta\cap B_r(x_0)\cap \{u>0\})\leqslant
C(a_0,a_1,a(1),N,M_0,\lambda_0,\Lambda_0,K_0)\delta r^{N-1}.
\end{eqnarray*}

\vs 0.2cm\n \emph{\underline{$2^{nd}$ Case}} :
$~~\displaystyle{{{a(t)}\over t}}$ is nondecreasing in $(0,t_*)$

\vs 0.2cm\n \n For each $\epsilon>0$ and
$\eta=a(2a^{-1}(\delta))$, we consider the function
\begin{equation*}
G(t)=\begin{cases} \displaystyle{
{{a((\epsilon+(a^{-1}(\eta))^2)^{1/2})}\over
{(\epsilon+(a^{-1}(\eta))^2)^{1/2}}}a^{-1}(\eta)} ~~
&\text{ if } ~~t> a^{-1}(\eta)\\
\displaystyle{{{a((\epsilon+(a^{-1}(\eta))^2)^{1/2})}\over
{(\epsilon+(a^{-1}(\eta))^2)^{1/2}}}t} ~~
&\text{ if } ~~|t|\leqslant a^{-1}(\eta)\\
 \displaystyle{-{{a((\epsilon+(a^{-1}(\eta))^2)^{1/2})}\over
{(\epsilon+(a^{-1}(\eta))^2)^{1/2}}}a^{-1}(\eta)} ~~ &\text{ if } ~~
t<-a^{-1}(\eta).
\end{cases}
\end{equation*}

\n $G$ is Lipschitz continuous and we have

$$0 \leqslant
G'(t)={{a((\epsilon+(a^{-1}(\eta))^2)^{1/2})}\over
{(\epsilon+(a^{-1}(\eta))^2)^{1/2}}}\chi_{\{|t|<
a^{-1}(\eta)\}}.$$

\n Using (1.5) and the fact that $\zeta=1$ in $B_{r}$, and since $|G(t)|\leqslant
a((\epsilon+(a^{-1}(\eta))^2)^{1/2}$, we obtain from (6.27)-(6.28) by taking into
account (6.18), we obtain
\begin{eqnarray*}
&&\int_{B_{r}\cap\{|u_{\epsilon x_i}|<a^{-1}(\eta)\}}
{{a((\epsilon+(a^{-1}(\eta))^2)^{1/2})}\over
{(\epsilon+(a^{-1}(\eta))^2)^{1/2}}}{{a(t_\epsilon)}\over{t_\epsilon}}
|\nabla u_{\epsilon
x_i}|^2dx\nonumber\\
&&\hskip 0.5cm \leqslant {2\over r} a_1 \int_{B_{2r}(x_0)}
|G(u_{\epsilon x_i})|{{a(t_\epsilon)}\over{t_\epsilon}}|\nabla
u_{\epsilon x_i}|
dx+K_0a((\epsilon+(a^{-1}(\eta))^2)^{1/2}(2r)^{N-1}.\nonumber\\
&&
\end{eqnarray*}

\n Since $\{|\nabla u_{\epsilon}|<a^{-1}(\eta)\}\subset \{|u_{\epsilon x_i}|<a^{-1}(\eta)\}$,
we obtain by using the monotonicity of $\displaystyle{{{a(t)}\over
t}}$ in $(0,t_*)$

\begin{eqnarray}\label{6.35}
&&\int_{B_{r}\cap \{|\nabla u_{\epsilon}|<a^{-1}(\eta)\}}
\Big({{a(t_\epsilon)}\over{t_\epsilon}}\Big)^2 |\nabla u_{\epsilon
x_i}|^2dx\nonumber\\
&&\hskip 1cm \leqslant \int_{B_{r}\cap\{|u_{\epsilon
x_i}|<a^{-1}(\eta)\}} {{a((\epsilon+(a^{-1}(\eta))^2)^{1/2})}\over
{(\epsilon+(a^{-1}(\eta))^2)^{1/2}}}{{a(t_\epsilon)}\over{t_\epsilon}}
|\nabla u_{\epsilon
x_i}|^2dx\nonumber\\
&&\hskip 1cm \leqslant{2\over r} a_1( J^i_1 +
J^i_2)+K_02^{N-1}a((\epsilon+(a^{-1}(\eta))^2)^{1/2}r^{N-1}.
\end{eqnarray}

\n Since
$$|G(u_{\epsilon x_i})- G(u_{ x_i})|\leqslant
{{a((\epsilon+(a^{-1}(\eta))^2)^{1/2})}\over
{(\epsilon+(a^{-1}(\eta))^2)^{1/2}}}|u_{\epsilon x_i}- u_{
x_i}|,~~\forall x\in B_{2r}(x_0)$$

\n we deduce that $|G(u_{\epsilon x_i})- G(u_{ x_i})|
\longrightarrow 0$ in $L^2(B_{2r}(x_0))$, as $\epsilon\rightarrow 0$.

\n As in the first case, we have up to a subsequence

$${{a((\epsilon+|\nabla u_\epsilon|^2)^{1/
2})}\over{(\epsilon+|\nabla u_\epsilon|^2)^{1/ 2}}}\nabla
u_{\epsilon x_i}\rightharpoonup W_i \qquad \hbox{in }
L^2(B_{2r}(x_0)) .$$

\n Therefore
\begin{equation}\label{6.36}
J^i_1=\int_{B_{2r}(x_0)} |G(u_{\epsilon x_i})- G(u_{ x_i})|
|\big({{a((\epsilon+|\nabla u_\epsilon|^2)^{1/
2})}\over{(\epsilon+|\nabla u_\epsilon|^2)^{1/ 2}}}|\nabla
u_{\epsilon x_i}| dx \longrightarrow
0,~~\text{as}~~\epsilon\rightarrow 0.
\end{equation}

\n For $J^i_2$, we have

\begin{eqnarray}\label{6.37}
&&J^i_2= \int_{B_{2r}(x_0)\cap \{\nabla u \neq 0\}} |G(u_{ x_i})|
.{{a((\epsilon+|\nabla u_\epsilon|^2)^{1/
2})}\over{(\epsilon+|\nabla u_\epsilon|^2)^{1/ 2}}}|\nabla
u_{\epsilon x_i}|  \qquad \hbox{ since  } G(0)=0\nonumber
\\
&& \leqslant a((\epsilon+(a^{-1}(\eta))^2)^{1/2})
\int_{B_{2r}(x_0)\cap \{\nabla u \neq 0\}} {{a((\epsilon+|\nabla
u_\epsilon|^2)^{1/ 2})}\over{(\epsilon+|\nabla u_\epsilon|^2)^{1/
2}}}|\nabla
u_{\epsilon x_i}|\nonumber\\
&& \leqslant a((\epsilon+(a^{-1}(\eta))^2)^{1/2})
|B_{2r}(x_0)|^{1/2}.
 \Big( \int_{B_{2r}(x_0)\cap \{\nabla u
\neq 0\}} \Big[{{a((\epsilon+|\nabla u_\epsilon|^2)^{1/
2})}\over{(\epsilon+|\nabla u_\epsilon|^2)^{1/ 2}}}|D^2
u|\Big]^2\Big)^{1/2}.\nonumber\\
&&
\end{eqnarray}

\vs 0,5cm \n As in the first case, we have $
O_{\delta}\subset \overline{B}_{1/2}\cap \{|\nabla u_{\epsilon}|<a^{-1}(\eta)\}$.
Since we have also $\{|\nabla u_\epsilon|<a^{-1}(\eta)\}\subset \{|
u_{\epsilon x_i}|<a^{-1}(\eta)\}$, we obtain from (6.35)

\begin{eqnarray}\label{6.38}
&&\int_{B_r(x_0)\cap O_\delta\cap\{u>0\}}
\Big({{a(t_\epsilon)}\over{t_\epsilon}}\Big)^2 |\nabla u_{\epsilon
x_i}|^2dx\leqslant {2\over r} a_1(J_1^i+J_2^i)\nonumber\\
&&\hskip 4cm +K_02^{N-1}a((\epsilon+(a^{-1}(\eta))^2)^{1/2}r^{N-1}.
\end{eqnarray}

\n Using (6.1), (6.37), we get from (6.38)
\begin{eqnarray}\label{6.39}
&&{1\over{(2+a_1)^2}}\int_{B_r(x_0)\cap O_\delta}
(\lambda_0H_\epsilon(u_\epsilon)^2dx\leqslant {2a_1\over
r}\Big(\sum_{i=1}^{i=N}J_1^i\nonumber\\
&& +N a((\epsilon+(a^{-1}(\eta))^2)^{1/2}) {2\over
r}|B_{2r}(x_0)|^{1/2}. \Big( \int_{B_{2r}(x_0)\cap \{\nabla u \neq
0\}} \Big[{{a(t_\epsilon)}\over{t_\epsilon}} |D^2
u_\epsilon|\Big]^2dx\Big)^{1/2}\Big)\nonumber\\
&&+NK_02^{N-1}a((\epsilon+(a^{-1}(\eta))^2)^{1/2}r^{N-1}.
\end{eqnarray}

\n Letting $\epsilon\rightarrow 0$ in (6.39), and using (6.36) together with
Lemma 6.2, we obtain
\begin{eqnarray*}
&&{\lambda_0^2\over{(2+a_1)^2}}\int_{B_r(x_0)\cap
O_\delta\cap\{u>0\}}dx\\
&&~~\leqslant {2a_1\over r}N\eta |B_{2r}(x_0)|^{1/2}.
 \Big( \int_{B_{2r}(x_0)\cap \{\nabla u
\neq 0\}} \Big[{{a(|\nabla u|)}\over{|\nabla u|}} |D^2
u|\Big]^2dx\Big)^{1/2}+NK_02^{N-1}\eta
r^{N-1}\nonumber\\
&&~~\leqslant N\eta {2\over r}
|B_{2r}(x_0)|^{1/2}.(C_4(a_0,a_1,a(1),N,M_0,\Lambda_0, K_0)|B_{2r}(x_0)|)^{1/2}+NK_02^{N-1}\eta
r^{N-1},
\end{eqnarray*}

\n which can be written
\begin{eqnarray*}
{\cal L}^N(O_\delta\cap B_r(x_0)\cap \{u>0\})\leqslant C_5=
C_5(a_0,a_1,a(1),\lambda_0,N,M_0,\Lambda_0, K_0)\delta r^{N-1}.
\end{eqnarray*}

\qed

\vs 0.5cm\n Now we give the proof of Theorem 1.5.

\vs 0,3cm\n \emph{Proof of Theorem 1.5.} Let $r\in(0,1/4)$, $B_r(x_0)\subset B_1$ with
$x_0\in
\partial \{u>0\}\cap B_{1/2}$  and $\delta>0$. First we recall the
definition of Hausdorff measure.  Let $E$ be a subset of
$\mathbb{R}^N$ and $s\in [0,\infty)$. The $s-$dimensional Hausdorff
measure of $E$ is defined by
$$H^s(E)=\lim_{\delta\rightarrow 0}H_\delta^s(E)=\sup_{\delta>0}H_\delta^s(E),$$

\n where for $\delta>0$,
$$H_\delta^s(E)=\inf\Big\{~\sum_{j=1}^\infty\alpha(s)\Big({{diam (C_j)}\over
2}\Big)^s~|~ E\subset \displaystyle{\bigcup_{j=1}^\infty
}C_j,~diam (C_j)\leqslant \delta~\Big\}$$

\n $\alpha(s)=\displaystyle{{\pi^{s/2}}\over {\Gamma(s/2+1)}}$,
$\Gamma(s)=\displaystyle{\int_0^\infty e^{-t}t^{s-1}dt}$, for
$s>0$ is the usual Gamma function.

\vs 0.2cm\n Let $E=\partial \{u>0\}\cap B_r(x_0)$ and denote by
$(B_\delta(x_i))_{i\in I}$ a finite covering of $E$, with $x_i\in
\partial \{u>0\}$ and $P(N)$ maximum overlapping.

\n From the proof of Theorem 3.1 of \cite{[CL1]}, we know that there
exists a constant $c_0=c(a_0,a_1,N,M_0,\lambda_0, \Lambda_0)$ such that
\begin{equation}\label{6.17}
\forall i\in I~~\exists y_i\in B_\delta(x_i)~: ~~
B_{c_0\delta}(y_i)\subset B_\delta(x_i)\cap \{u>0\}\cap O_\delta.
\end{equation}

\n We deduce from (6.17) and Lemma 6.2 that
\begin{eqnarray*}
\sum_{i\in I}{\cal L}^N(B_1)c_0^N\delta^N&=&\sum_{i\in I}{\cal L}^N(B_{c_0\delta}(y_i))\leqslant
\sum_{i\in
I}{\cal L}^N(B_\delta(x_i)\cap \{u>0\}\cap O_\delta)\nonumber\\
&&\leqslant P(N){\cal L}^N(B_r(x_0)\cap \{u>0\}\cap
O_\delta)\leqslant P(N)C\delta r^{N-1}
\end{eqnarray*}

\n where $C=C(N,a_0,a_1,a(1),M,\lambda_0,\Lambda_0,K_0)$. This leads to
\begin{eqnarray*}
\sum_{i\in I}\alpha(N-1)\Big({{diam B_\delta(x_i)}\over
2}\Big)^{N-1}&\leqslant& {{\alpha(N-1)}\over{{\cal
L}^n(B_1)c_0(N,a_0,a_1,M_0,\lambda_0, \Lambda_0)}}P(N)C
r^{N-1}\\
&&=\overline{C}=\overline{C}(N,a_0,a_1,a(1),M_0,\lambda_0,\Lambda_0,K_0) r^{N-1}.
\end{eqnarray*}

\n It follows  that $H_\delta^{N-1}(\partial \{u>0\}\cap
B_r(x_0))\leqslant \overline{C} r^{N-1}$. Letting
$\delta\rightarrow 0$, we get
$$H^{N-1}(\partial \{u>0\}\cap
B_r(x_0))\leqslant \overline{C}r^{N-1}.$$ \qed

\vskip 0,5cm\n Finally we give the proof of the stability result in Corollary 1.1..
We will use the same notation $I(u)$ to denote the local coincidence set $\{u=0\}\cap B_{1/2}$
for each $u\in {\cal F}_A$.
We also set $E_{-\epsilon}=\{x\in E~:~d(x,E^c)>\epsilon~\}$, where $E^c$ is the complement of the set $E$. Then we have

\vs 0,3cm\n \emph{Proof of Corollary 1.1.} We set $E_{-\epsilon}=\{x\in E~:~d(x,E^c)>\epsilon~\}$, where $E^c$ is the complement of the set $E$.

\vs 0,3cm\n $i)$ First note that if $x\in I(u_1)$, then $u_1(x)=0$ and then
$$u_2(x)=|u_1(x)-u_2(x)|\leqslant |u_1-u_2|_\infty < \widetilde{A}(\epsilon).$$
We deduce that $I(u_1)\subset \{u_2<\widetilde{A}(\epsilon)\}\cap B_{1/2}$ and in particular
$$I(u_1)\setminus I(u_2) \subset \{0<u_2<\widetilde{A}(\epsilon)\}\cap B_{1/2}.$$
Since moreover
$\{0<u_2<\widetilde{A}(\epsilon)\}\cap B_{1/2}\subset O_{\epsilon}\cap B_{1/2}$, it follows then by Lemma 6.1 that we have
$${\cal L}^N(I(u_1)\setminus I(u_2))\leqslant {\cal L}^N(O_{\epsilon}\cap B_{1/2})\leqslant C(N,a_0,a_1,a(1),M_0,\lambda_0)\epsilon.$$
Similarly we prove that
$${\cal L}^N(I(u_2)\setminus I(u_1))\leqslant  C(N,a_0,a_1,a(1),M_0,\lambda_0)\epsilon.$$
It follows that
${\cal L}^N(
I(u_1)\div I(u_2))={\cal L}^N(I(u_1)\setminus I(u_2))+{\cal L}^N(I(u_2)\setminus I(u_1))\leqslant C\epsilon$.

\vskip 0,3cm \n
$ii)$ We have already proved that $I(u_1)\subset \{u_2<\widetilde{A}(\epsilon)\}\cap B_{1/2}$. Let us show that $(I(u_2))_{-C\epsilon}\subset I(u_1)$. So let $x\in (I(u_2))_{-C\epsilon}$. Then we have $x\in I(u_2)$ and
$d(x,(I(u_2))^c)>C\epsilon$. Now assume that $u_1(x)>0$, we obtain by using Lemma 3.1 of \cite{[CL1]}, for $C$ large enough
$$\sup_{\partial B_{C\epsilon}(x)}u_1 \geqslant {N\over \lambda_0}\widetilde{A}\Big({\lambda_0\over N}C\epsilon\Big)+u_1(x)>\widetilde{A}(\epsilon).$$
We deduce that there exists $x'\in \partial B_{C\epsilon}(x)$ such that $u_1(x')>\widetilde{A}(\epsilon)$. If $u_2(x')=0$, we obtain
$u_1(x')\leqslant\widetilde{A}(\epsilon)$. Hence we have necessarily $u_2(x')>0$. This leads to
$d(x,(I(u_2))^c)\leqslant |x-x'|=C\epsilon$, which contradicts the fact that $d(x,(I(u_2))^c)>C\epsilon$. We conclude that $u_1(x)=0$, i.e. $x\in I(u_1)$.\qed

\vskip 0,5cm \n\emph{Acknowledgments.} The first two authors are
grateful for the facilities and excellent research conditions at the
Fields Institute. This work was completed while the second author
was visiting the CMAF, and while the third author was visiting the
Fields Institute. They would like to thank these institutions
for the hospitality during their visits.

\vs 0,5cm

\end{document}